\documentclass[amsart,11pt]{article}

\usepackage[hscale=0.7,centering]{geometry}
\usepackage{graphicx}
\usepackage{float}
\usepackage{amssymb,amsmath,amsfonts,amssymb,url}
\usepackage{times}
\usepackage{color}
\usepackage{amsfonts}[11pt]
\usepackage{amsmath}[11pt]
\usepackage{pdfsync}
\numberwithin{equation}{section}

\def\qed{ \hfill $\Box$}
\def\eqd{\,{\buildrel d \over =}\,}

\def\tendp{\,{\buildrel P \over \rightarrow}\,}

\def\tendl2{\,{\buildrel L_2 \over \rightarrow}\,}

\makeatletter

\renewcommand{\section}{
  \@startsection
  {section}
  {1}
  {0pt}
  {1.1\baselineskip}
  {0.2\baselineskip}
  {\sc \centering}
}
\makeatother

\newcommand{\Z}{\mathbb Z}

\definecolor{db}{rgb}{0.1,0,0.75}
\definecolor{lm}{cmyk}{0 ,1,0,0}

\newtheorem{lem}{Lemma}

\newtheorem{thm}{Theorem}

\newtheorem{cor}{Corollary}[thm]

\renewcommand{\P}{\text{P}}

\newcommand{\barphi}{\overline{\Phi}}
\newcommand{\baru}{\overline{u}}

\begin{document}
\def\shorttitle{\textbf{Conga Line}}
\title{\textbf{\Large \sc The Brownian Conga Line}}
\author{\sc Sayan Banerjee\\\\
 Department of Statistics \\ University of Warwick \\ Sayan.Banerjee@warwick.ac.uk}

\date{\today}


\maketitle
\begin{abstract}
We introduce a new model called the Brownian Conga Line. It is a random curve evolving in time, generated when a particle performing a two dimensional Gaussian random walk leads a long chain of particles connected to each other by cohesive forces. 
We approximate the discrete Conga line in some sense by a smooth random curve and subsequently study the properties of this smooth curve.
\end{abstract}
\tableofcontents

\section{\textbf{Introduction}}
The Conga Line is a Cuban carnival march that has become popular in many cultures over time. It consists of a queue of people, each one holding onto the person in front of him. The person at the front of the line can move as he will, and the person holding onto him from behind follows him. The third person in the queue follows the second, and so on. Often people keep on joining the line over time by attaching themselves to the last person in the line. As the Conga Line grows in time, it displays interesting motion patterns where the randomness in the motion of the first person propagates down the line, diminishing in influence as it moves further down. The Conga line also appears in molecular biology as a model for various long polymers built of smaller monomers and is closely linked to (a discrete version of) the \textit{curve shortening problem} (see \cite{chouzhu}). In this article, we devise a mathematical formulation of the Conga Line and study its properties.

The formulation is as follows.

Let $Z_k, \ k \ge 1$, be $i.i.d$ standard $2$-dimensional normal random variables. Fix some $\alpha \in (0,1)$. 
Let $X_1(0)=0$ and $X_1(n)=\sum_{i=1}^nZ_i$ for $n \ge 1$. This denotes the leading particle, or the \textit{tip} of the Conga line.

\begin{figure}[H]\label{figure:Congasteps}
\begin{center}
\includegraphics[width=120mm,height=90mm]{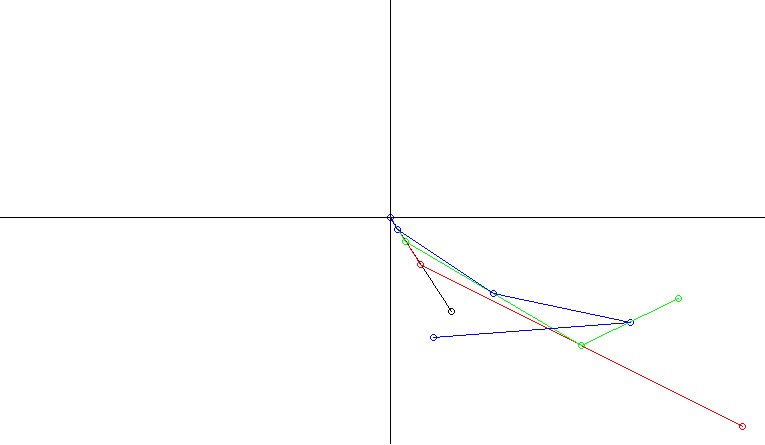}
\caption{\footnotesize The Conga line for $n=1$ (black), $n=2$ (red), $n=3$ (green) and $n=4$ (blue) when $\alpha=0.5$. (courtesy: Mary Solbrig)}
 \end{center}
 \end{figure}

Now, we define processes $X_k$ inductively as follows. Suppose that $\{X_k(n), \ n \ge 0\}$ have already been defined for $1 \le k \le j$. Then we let $X_{j+1}(0)=0$ and
\begin{equation}\label{equation:recursion}
X_{j+1}(n)=(1-\alpha)X_{j+1}(n-1) + \alpha X_{j}(n-1)
\end{equation}
for $n \ge 1$. Here, the process $X_k$ denotes the motion over time of the particle at distance $k$ from the leading particle. The relation (\ref{equation:recursion}) describes the manner in which a particle $X_{j+1}$ follows the preceding particle $X_j$. It is easy to check from (\ref{equation:recursion}) that $X_j(n)=0$ for all $j > n$. These represent the particles at rest at the origin at time $n$. Note that the $j$-th particle $X_j$ joins the Conga Line at time $j$. See Figure \ref{figure:Congasteps} for the construction of the Conga line for $n=1,2,3,4$.

 The Conga line at time $n$ is defined as the collection of random variables $\{X_k(n), \ k \le n\}$.
 
One can also think of this model as a discrete version of a long string or molecule whose tip is moving randomly under the effect of an erratic force and the rest of it performs a constrained motion governed by the tip together with the cohesive forces. Burdzy and Pal \cite{burdzypal} performed some simulations (see Figure 2) which led them to make the following observations:
\begin{enumerate}
\item For a fixed large $n$, the locations of the particles $\{X_k(n), \ k \ge 1\}$ sufficiently away from the tip look like a `smooth' curve, and the smoothness increases as we move away from the tip.
\item For $k$ significantly larger than $1$, there is very little variability in the location of the particles over short periods of time.
\item The small loops in the curve tend to die out over time. Just before death, they look `elongated' and their death site forms a cusp.
\item The particles near the origin seem to \textit{freeze} showing very little movement over time. Moreover, the \textit{direction} in which the particles come out of the origin seems to stabilise over time.
\end{enumerate}
All the above observations need precise mathematical formulations. Once the rigorous foundations are established, we can ask the correct questions and try to answer them. This, broadly, is the goal of the article.

We give a brief outline of the content of each section.
In Section \ref{section:discrete}, we try to make mathematical sense of the statement `the process looks like a smooth curve'. This is the toughest challenge as the Conga line, unlike most known stochastic processes which can be approximated by continuous models, does not seem to have an interesting scaling limit under \textit{uniform scaling}. This is because if we look at the Conga Line for any fixed $n$, the distance between the particles decays as we move away from the tip, i.e., increase $k$, as suggested by Figure 2. But it is precisely why this is a novel model, which exhibits \textit{particles moving in different scales `in the same picture'}. The particles near the tip are wider spaced and their paths mostly resemble a Gaussian random walk, but those for large $k$ are more closely packed and the Conga Line looks very smooth in this region (see Figure 2). To circumvent this problem, we describe a \textit{coupling} between our discrete process $\{X_k(n), \ k \le n\}$ and a smooth random process $\{u(x,t): (x,t) \in \mathbb{R}^{+2}\}$ such that, when observed sufficiently away from the tip, more precisely for $k \ge n^{\epsilon}$ for any fixed $\epsilon>0$ and large $n$, the points $X_{k+1}(n)$ are uniformly close to the points $u(k,n)$. Thus, $u$ serves as a \textit{smooth approximation} to the discrete process $X$ in a suitable sense. The $x$ variable of $u$ represents distance from the tip and the $t$ variable represents time. Future references to the Conga line refer to this smooth version $u$. We close the section by presenting another smooth process $\baru$ that also serves as an approximation in the same sense, and is more intuitive when considering the motion of individual particles, i.e. trajectories of the form $\{X_k(n): n \ge k\}$ for fixed $k$. It is also used in Section \ref{section:freezing} to study the phenomenon of freezing of the Conga line near the origin.\\\\
\textbf{Note: }In the following, we will be investigating the properties of the continuous \textit{two dimensional Conga line} (which is the smooth approximation to our original discrete model) as well as the process corresponding to its $x$ (equivalently $y$) coordinate, which we will call the (continuous) \textit{one dimensional Conga line}. To save additional notation, we will denote both of these by $u$, but the dimension will be clear from the context.\\\\
In Section \ref{section:continuous}, we study the properties of the continuous, one dimensional Conga line $u$. First, we investigate the phenomenon of the particles at different distances from the tip moving in `different scales' suggested by their different order of variances. The particles near the tip move in the same scale as the leading Gaussian random walk as indicated by their variance being $O(t)$, while those far away from the tip show very little movement away from the origin, as indicated by exponentially decaying variances. Furthermore, there exists a \textit{cutoff} near $x=\alpha t$, where the \textit{variance shows a sharp transition from `very large to very small'}. We identify this and study the fine changes in variance around this point. 
\\
\indent Next, using the scaling properties of Brownian motion, we show that for fixed $t$, the Conga line (in both one and two dimensions) can be \textit{scaled} so that the space variable $x$ runs in $[0,1]$. We \textit{call this scaled version} $u_t$ and study its analytical properties. Upper bounds on the growth rate of the derivatives show that $u_t$ is real analytic. We also make a detailed study of the covariance structure of the derivatives. This turns out to be a major tool in studying the subsequent properties like critical points, length, loops, etc.\\
\indent With the basic framework of the Conga line established, we set out to investigate its finer properties. We investigate the distribution of critical points of the scaled one dimensional Conga line $u_t$, i.e., points at which the derivative vanishes. The number of critical points in an interval serves as a measure of how wiggly the Conga line looks on that interval. The critical points are distributed as a point process on the real line and we show using an expectation meta-theorem for smooth Gaussian fields (see \cite{adlertaylor} p. 263) that its first intensity at $x$ (for a large time $t$) is approximately of the form $\sqrt{t} x ^{-1/2}$. This shows that, though the typical number of these points in a given interval is $O(\sqrt{t})$ for large $t$, the proportion of critical points around $x$ decreases as $\displaystyle{x^{-\frac{1}{2}}}$ as we go farther away from the tip. We also show subsequently using second moment estimates that the critical points are reasonably \textit{well-spaced} and they do not tend to crowd around any point. Furthermore, we show that the \textit{first intensity is a good estimate of the point process itself} as for a given interval $I$ sufficiently away from the ends $x=0$ and $x=1$, the ratio $\displaystyle{\frac{N_t(I)}{\mathbb{E}N_t(I)}}$ goes to one in probability as $t$ grows large.

In Section \ref{section:twodimensional} we give a fluctuation estimate for the continuous two dimensional Conga line. This, along with the approximation result on Section \ref{section:discrete}, tells us that for any $\delta \in (0,\alpha)$ and for large $n$, the \textit{linear interpolation} of the discrete Conga line, along with its derivative (which exists everywhere except at integer points), comes uniformly close to $(u(\cdot, n), \partial_xu(\cdot, n))$ on the entire interval $[\delta n, \alpha n]$. This explains why the discrete Conga line \textit{looks smooth} away from the tip. We then study properties of the scaled two dimensional Conga line, like length and number of loops. We also investigate a strange phenomenon. Although the mechanism of subsequent particles following the preceding ones and `cutting corners' results in progressively smoothing out the roughness of the Gaussian random walk of the tip, we see that as $t$ increases, the \textit{scaled Conga line looks more and more like Brownian motion} in that the sup-norm distance between them on $[0,1]$ is roughly of order $t^{-1/4}$ (with a log correction term). This can be explained by the fact that the noticeable smoothing of the paths of the unscaled Conga line takes place in a window of width $\sqrt{t}$ around each point, which translates into a window of width $t^{-1/2}$ as we scale space and time by $t$. Thus in the scaled version, the \textit{smoothing window becomes smaller with time}, resulting in this phenomenon. Consequently, the scaled Conga line $u_t$ for large $t$ serves as a \textit{smooth approximation to Brownian motion} which smooths out microscopic irregularities but retains its macroscopic characteristics.

In Section \ref{section:loopevolution} we study the evolution of loops in the family of two dimensional paths that the particles at successively larger distances from the tip trace out. We study this evolution under a metric similar to the Skorohod metric. It turns out that with probability one, every singularity, i.e a point where the speed of the curve becomes zero, in a particle path is a \textit{cusp singularity} (looks like the graph of $y=x^{2/3}$ in a local co-ordinate frame). Furthermore, there is a \textit{bijection between dying loops and cusp singularities} in the sense that small loops die (i.e. the loop shrinks to a point) creating cusp singularities, and conversely, if such a singularity appears in the path of some particle, we can find a loop in the path of the immediately preceding particles, and it dies creating the singularity.

Finally, in Section \ref{section:freezing}, we investigate the phenomenon of \textit{freezing near the origin}. We work with the smooth approximation $\baru$, and show that for an appropriate choice of a sequence $x_t$ of distances from the tip such that the particles at these distances remain sufficiently close to the origin, $\baru(x_t,t)$ converges almost surely and in $L^2$, and find the limiting function.
\\\\
\textbf{Notation: }Before we proceed, we clarify the notation that we will be using here:
\begin{enumerate}
\item[(i)] If  $g$ is a random function and $V$ is a random variable with distribution function $F$ and independent of $g$, then $$\mathbb{E}_Vg(V)=\int g(v)dF(v)$$ denotes the expectation with respect to $V$ for a fixed realisation of $g$.
\item[(ii)] $\Phi$ denotes the normal distribution function and $\barphi=1-\Phi$. We denote the corresponding density function by $\phi$.
\item[(iii)] For any function $f$ of several variables, $\partial^k_x f$ denotes the partial derivative of $f$ with respect to the variable $x$ taken $k$ times.
\item[(iv)] For functions $f,g: [0,\infty) \rightarrow \mathbb{R}^+$, $f(t) \sim g(t)$ means that there $f$ and $g$ have the same growth rate in $t$, i.e., there exists a constant $C$ such that $$\displaystyle{\frac{f(t)}{g(t)}\vee \frac{g(t)}{f(t)} \le C}$$ for all sufficiently large $t$.
\item[(v)] For a family of real-valued functions $\{f_t: t \in (0, \infty)\}$ defined on a compact set $I \subseteq \mathbb{R}^k$ and a function $a: (0,\infty) \rightarrow [0,\infty)$, we say that $$f_t=O^{\infty}\left(a(t)\right) \ \text{on} \ I$$ if $$\sup_{t \in (0,\infty)}\frac{\sup_{x \in I}|f_t(x)|}{a(t)} \le C$$ for some constant $C<\infty$. Sometimes, (by abuse of notation) we will write $$f_t(x)=O^{\infty}\left(a(t)\right) \ \text{for} \ x \in I$$ to denote the same.
\end{enumerate}
\section{\textbf{The discrete Conga line}}\label{section:discrete}
We set out by finding a useful way to express $X_{k}(n)$ in terms of $X_1(n)$. This turns out to be the starting point  in our approximation procedure of the discrete Conga line $\{X_k(n): k \ge 1\}$ for sufficiently large $k$ by a smooth curve.\\\\
Let $T_1,T_2,\dots$ be i.i.d Geom($\alpha$) and let $$\Theta_j=\sum_{i=1}^jT_i.$$ Then $\Theta_j \sim NB(j,\alpha)$, where $NB(a,b)$ represents the Negative Binomial distribution with parameters $a$ and $b$. It follows from the recursion relation (\ref{equation:recursion}) that one can write $$X_k(n)=\mathbb{E}_{T_1}X_{k-1}(n-T_1),$$ where we set $X_k(n)=0$ for all $n \le 0$, for each $k$.\\\\
By induction, we get
\begin{equation}\label{equation:expression}
X_k(n)= \mathbb{E}_{\Theta_{k-1}} X_1(n-\Theta_{k-1})=\sum_{m=0}^{n-k+1} {m+k-2\choose m} (1-\alpha)^m\alpha^{k-1}X_1(n-k+1-m).
\end{equation}

\subsection{Approximation by a smooth process}\label{subsection:approximation}
Here we show that for any fixed $\epsilon >0$, the discrete Conga line $\{X_k(n)\}$ can be approximated uniformly in $k$, for $n^{\epsilon}\le k \le n$, for large $n$, by a smooth process $\{u(x,t): (x,t) \in \mathbb{R}^+ \times \mathbb{R}^+\}$ evaluated at integer points $(k,n)$. This process arises as a smoothing kernel acting on Brownian motion. Furthermore, the increments $X_{k+1}(n)-X_k(n)$ (which correspond to the `discrete derivative' of $X$ for fixed $n$) can also be uniformly approximated by the derivative of $u(\cdot, n)$ at $k$.\\\\
Let $B_l\sim \operatorname{Bin}(l,\alpha)$. From (\ref{equation:expression}) and the fact that $$P(\Theta_{k-1}\le l)=P(B_l\ge k-1),$$ we get
\begin{eqnarray}\label{eqnarray:movone}
X_k(n)&=&\mathbb{E}_{\Theta_{k-1}}X_1(n-\Theta_{k-1})=\sum_{j=k-1}^{n}P(\Theta_{k-1}=j)\sum_{l=1}^{n-j}Z_l\nonumber\\
&=&\sum_{l=1}^{n-k+1}Z_l\sum_{j=k-1}^{n-l}P(\Theta_{k-1}=j)=\sum_{l=k-1}^{n-1} P(\Theta_{k-1}\le l)Z_{n-l}\nonumber\\
&=&\sum_{l=k-1}^{n-1} P(B_l\ge k-1)Z_{n-l}
\end{eqnarray}
The above expression also yields
\begin{eqnarray}\label{eqnarray:movtwo}
X_{k+1}(n)-X_k(n)= -\sum_{l=k-1}^nP(B_l= k-1)Z_{n-l}.
\end{eqnarray}
The next step is the key to the approximation. We obtain a \textit{coupling} between a Brownian motion and our process $X$. Let $(\Omega, \mathcal{F}, P)$ be a probability space supporting a Brownian motion $\{W(t): t \ge 0\}$, where we set $W(t)=0$ for all $t \le 0$. Then $$X_k(n)= \sum_{l=k-1}^{n-1}P(B_l\ge k-1)(W(n-l)-W(n-l-1))$$ gives the desired coupling on this space. Note that we can write
\begin{eqnarray*}
X_{k+1}(n)= \int_0^n g(k,z)dW^n_z
\end{eqnarray*}
where $$g(k,z)= P(B_{\lfloor z \rfloor} \ge k).$$
Here $\lfloor z \rfloor$ denotes the largest integer less than or equal to $z$, and $W^t_z= W(t)-W(t-z)$, $0\le z\le t$, is the time reversed Brownian motion from time $t$.\\\\
Let $\sigma=\sqrt{\alpha(1-\alpha)}$. Consider the ``space-time" process
\begin{eqnarray}\label{eqnarray:SIP}
 u(x,t)&=& \int_0^t\barphi\left(\frac{x-\alpha z}{\sigma \sqrt{z}}\right)dW^t_z\nonumber\\
 &=& \int_0^tW(t-z)(\sqrt{2\pi})^{-1}\left(\frac{x+\alpha z}{2\sigma z^{3/2}}\right)\exp\left(-\frac{(x-\alpha z)^2}{2 \sigma^2z}\right)dz
\end{eqnarray}
(We obtain the second expression from the first by an application of \textit{stochastic integration by parts}, see \cite{revuzyor}). Note that $u$ is smooth in $x$ and its first derivative is given by $$\partial_xu(x,t)= -\int_0^t\frac{1}{\sigma\sqrt{z}}\phi\left(\frac{x-\alpha z}{\sigma \sqrt{z}}\right)dW^t_z.$$
 We prove in what follows that for large $n$, and for $n^{\epsilon}\le k \le n$ for any fixed $\epsilon >0$, the points $X_{k+1}(n)$ are ``uniformly close" to the points $u(k,n)$ ($u$ evaluated at integer points) for the given range of $k$. Similarly, the increments $X_{k+2}(n)-X_{k+1}(n)$ are also uniformly close to $\partial_x u(k,n)$ for $k$ in the given range. Our strategy is to first consider a discretized version of the process $u(x,t)$ and its derivative $\partial_x u(x,t)$, respectively given by
\begin{eqnarray*}
\hat{u}(k,n)&=&\sum_{l=0}^{n}\barphi\left(\frac{k-\alpha l}{\sigma \sqrt{l}}\right)(W(n-l)-W(n-l-1)),\\
\widehat{\partial_x u}(k,n)&=&-\sum_{l=0}^{n}\frac{1}{\sigma\sqrt{l}}\phi\left(\frac{k-\alpha l}{\sigma \sqrt{l}}\right)(W(n-l)-W(n-l-1)).
\end{eqnarray*}
In Lemma 1, we give a bound on the $L^2$ distance between $X_{k+1}(n)$ and $\hat{u}(k,n)$ and between $X_{k+2}(n)-X_{k+1}(n)$ and $\widehat{\partial_x u}(k,n)$ for large $n$ when $n^{\epsilon} \le k \le n$. In Lemma 2, a similar bound is achieved for the $L^2$ distance between $\hat{u}(k,n)$ and $u(k,n)$ (and between $\widehat{\partial_x u}(k,n)$ and $\partial_x u(k,n)$). In Theorem \ref{thm:conga}, we prove using a Borel Cantelli argument that for large $n$ the two processes $X$ and $u$ (evaluated at integer points) come uniformly close on $n^{\epsilon} \le k \le n$. A similar result holds for the increments $X_{k+2}(n)-X_{k+1}(n)$ and the partial derivative $\partial_x u(k,n)$.\\\\
In the following, $C_1, \ C_2, \dots$ represent absolute constants, $C_{\epsilon}$, $C'_{\epsilon}$ denote constants that depend only on $\epsilon$, $C_p$ denotes a constant depending only on $p$ and $D_{\epsilon,p}$, $D'_{\epsilon,p}$ denote constants depending upon both $\epsilon$ and $p$.
\begin{lem}\label{lem:l2dist1}
Fix $\epsilon>0$. For $n^{\epsilon}\le k \le n$, 
\begin{eqnarray}\label{eqnarray:ltwo}
\displaystyle{\sum_{l=0}^{n}\left[P(B_l\ge k)-\barphi\left(\frac{k-\alpha l}{\sigma\sqrt{l}}\right)\right]^2 \le C_{\epsilon}\frac{\sqrt{\log k}}{\sqrt{k}}}\nonumber\\
\displaystyle{\sum_{l=0}^{n}\left[P(B_l = k)-\frac{1}{\sigma\sqrt{l}}\phi\left(\frac{k-\alpha l}{\sigma\sqrt{l}}\right)\right]^2 \le C_{\epsilon}\frac{(\log k)^{7/2}}{k^{3/2}}}
\end{eqnarray}
where $\sigma=\sqrt{\alpha(1-\alpha)}$. Consequently, $$\displaystyle{\mathbb{E}|X_{k+1}(n)-\hat{u}(k,n)|^2 \le C_{\epsilon}\frac{\sqrt{\log k}}{\sqrt{k}}}$$ and $$\displaystyle{\mathbb{E}|(X_{k+2}(n)-X_{k+1}(n))-\widehat{\partial_x u}(k,n)|^2 \le C_{\epsilon}\frac{(\log k)^{7/2}}{k^{3/2}}}$$ uniformly on $n^{\epsilon}\le k \le n$.
\end{lem}
\textbf{Proof: }Choose $C>0$ such that $\epsilon(C-1)\ge 2$ and $C\ge \frac{3}{2}$. Take $L_k=\lfloor\alpha^{-1}\sqrt{Ck \log k}\rfloor$. Then, we can write
\begin{eqnarray*}
\sum_{l=0}^{n}\left[P(B_l\ge k)-\barphi\left(\frac{k-\alpha l}{\sigma\sqrt{l}}\right)\right]^2&\le &\sum_{l=0}^{\lfloor\frac{k}{\alpha}\rfloor-L_k}\left[P(B_l\ge k)-\barphi\left(\frac{k-\alpha l}{\sigma\sqrt{l}}\right)\right]^2\\
&+&\sum_{l=\lfloor\frac{k}{\alpha}\rfloor-L_k}^{\lfloor\frac{k}{\alpha}\rfloor+L_k}\left[P(B_l\ge k)-\barphi\left(\frac{k-\alpha l}{\sigma\sqrt{l}}\right)\right]^2\\
&+&\sum_{l=\lfloor\frac{k}{\alpha}\rfloor+L_k}^n\left[P(B_l\ge k)-\barphi\left(\frac{k-\alpha l}{\sigma\sqrt{l}}\right)\right]^2\\\\
&=& S_1^{(k)} + S_2^{(k)} + S_3^{(k)}.
\end{eqnarray*}
Here, $S_1^{(k)}$ and $S_3^{(k)}$ correspond to the \textit{tails} of the distribution functions, and we shall show that they are negligible compared to $S_2^{(k)}$. To this end, note that
\begin{eqnarray*}
S_1^{(k)} \le 2k\alpha^{-1}\left({P^2(B_{\lfloor\frac{k}{\alpha}\rfloor-L_k}} \ge k) + \barphi^2\left(\frac{\alpha L_k}{\sigma\sqrt{\frac{k}{\alpha}}}\right)\right).
\end{eqnarray*}
Now, by Bernstein's Inequality, 
\begin{eqnarray*}
P(B_{\lfloor\frac{k}{\alpha}\rfloor-L_k} \ge k) \le \exp\left(-\frac{\alpha^2L_k^2/2}{\alpha(\lfloor\frac{k}{\alpha}\rfloor-L_k)+\alpha L_k/3}\right) \le \exp \left(-\frac{\alpha^2L_k^2}{2k}\right).
\end{eqnarray*} 
We also have
\begin{eqnarray*}
\barphi\left(\frac{\alpha L_k}{\sigma\sqrt{\frac{k}{\alpha}}}\right) \le \frac{\sigma\sqrt{k}}{\alpha^{3/2}L_k}\exp\left(-\frac{\alpha^{3}L_k^2}{2\sigma^2 k}\right).
\end{eqnarray*} 
Therefore, for large $k$,
\begin{eqnarray*}
S_1^{(k)} \le \frac{4k}{\alpha}\exp\left(-\frac{\alpha^2 L_k^2}{k}\right) \le \frac{C_1}{k^{C-1}} \le \frac{C_1}{\sqrt{k}}.
\end{eqnarray*}
Similarly, for $S_3^{(k)}$, we get
\begin{eqnarray*}
S_3^{(k)} &\le & 2(n-\lfloor k\alpha^{-1}\rfloor)\left({P^2(B_{\lfloor\frac{k}{\alpha}\rfloor+L_k}} < k) + \Phi^2\left(\frac{-\alpha L_k}{\sigma\sqrt{\frac{2k}{\alpha}}}\right)\right)\\
&\le & 4(n-\lfloor k\alpha^{-1}\rfloor)\exp\left(-\frac{\alpha^2L_k^2}{2k}\right) \le C_3(n-\lfloor k\alpha^{-1}\rfloor)k^{-C/2}\\
&\le & \frac{C_3(n-\lfloor k\alpha^{-1}\rfloor)}{\sqrt{k}n^{\epsilon(C-1)/2}} \le \frac{C_3}{\sqrt{k}}.
\end{eqnarray*}
Now, for $S_2^{(k)}$, we use the \textit{Berry Esseen Theorem} (see \cite{feller}).
\begin{eqnarray*}
S_2^{(k)} \le \sum_{l=\lfloor\frac{k}{\alpha}\rfloor-L_k}^{\lfloor\frac{k}{\alpha}\rfloor+L_k}\frac{1}{l}\le C_4\frac{\sqrt{C \log k}}{\sqrt{k}}.
\end{eqnarray*}
The second inequality in (\ref{eqnarray:ltwo}) follows similarly, except that we use the \textit{De Moivre-Laplace Theorem} (see \cite{DL}) in place of the Berry Esseen Theorem to bound the second sum. The De Moivre-Laplace Theorem gives a moderate deviations bound: it says that if $A_N$ is non-decreasing and satisfies $\displaystyle{\frac{A_N}{N^{1/6}} \rightarrow 0}$ as $N \rightarrow \infty$, then there is a constant $C$ depending only on $\alpha$ such that $$\sup_{|u-N\alpha| \le A_N \sqrt{N}}\left|\frac{\sigma\sqrt{N}P(B_N = k)}{\phi\left(\frac{k-\alpha N}{\sigma\sqrt{N}}\right)}-1\right| \le \frac{CA_N^3}{\sqrt{N}}.$$ Putting $A_N = \alpha^{-1}\sqrt{C\log N}$ in this expression yields the required bound. This completes the proof of the lemma. \qed
\begin{lem}\label{lem:l2dist2}
$\displaystyle{\mathbb{E}|\hat{u}(k,n)-u(k,n)|^2 \le C_{\epsilon}\frac{\sqrt{\log k}}{\sqrt{k}}}$ and $\displaystyle{\mathbb{E}|\widehat{\partial_x u}(k,n)-\partial_xu(k,n)|^2 \le C_{\epsilon}\frac{\sqrt{\log k}}{k^{3/2}}}$ \ uniformly on $n^{\epsilon}\le k \le n$.
\end{lem}
\textbf{Proof: }Write $\displaystyle{\hat{u}(k,n)=\int_0^n\hat{f}(k,z)dW_z^n}$ and $\displaystyle{u(k,n)=\int_0^nf(k,z)dW_z^n}$ where\\\\ $\displaystyle{\hat{f}(k,z)= \barphi\left(\frac{k-\alpha \lfloor z \rfloor}{\sigma\sqrt{\lfloor z \rfloor}}\right)}$ and $\displaystyle{f(k,z)=\barphi\left(\frac{k-\alpha z}{\sigma\sqrt{z}}\right)}$.\\\\
Then, we can decompose $\displaystyle{\mathbb{E}|\hat{u}(k,n)-u(k,n)|^2}$ as in the proof of Lemma \ref{lem:l2dist1} as follows:
\begin{eqnarray*}
\mathbb{E}|\hat{u}(k,n)-u(k,n)|^2&=&\int_0^n\left[\barphi\left(\frac{k-\alpha \lfloor z \rfloor}{\sigma\sqrt{\lfloor z \rfloor}}\right)-\barphi\left(\frac{k-\alpha j}{\sigma\sqrt{j}}\right)\right]^2dz\\
&\le & \int_0^{\lfloor\frac{k}{\alpha}\rfloor-L_k}\left[\barphi\left(\frac{k-\alpha \lfloor z \rfloor}{\sigma\sqrt{\lfloor z \rfloor}}\right)-\barphi\left(\frac{k-\alpha j}{\sigma\sqrt{j}}\right)\right]^2dz\\
&+& \int_{\lfloor\frac{k}{\alpha}\rfloor-L_k}^{\lfloor\frac{k}{\alpha}\rfloor+L_k}\left[\barphi\left(\frac{k-\alpha \lfloor z \rfloor}{\sigma\sqrt{\lfloor z \rfloor}}\right)-\barphi\left(\frac{k-\alpha j}{\sigma\sqrt{j}}\right)\right]^2dz\\
&+&\int_{\lfloor\frac{k}{\alpha}\rfloor+L_k}^{n}\left[\barphi\left(\frac{k-\alpha \lfloor z \rfloor}{\sigma\sqrt{\lfloor z \rfloor}}\right)-\barphi\left(\frac{k-\alpha j}{\sigma\sqrt{j}}\right)\right]^2 dz\\\\
&=&I_1^{(k)}+I_2^{(k)}+I_3^{(k)}.
\end{eqnarray*}
Now,
\begin{eqnarray}\label{eqnarray:partial}
\frac{\partial f}{\partial z}(k,z)=(\sqrt{2\pi})^{-1}\left(\frac{k+\alpha z}{2\sigma z^{3/2}}\right)\exp\left(-\frac{(k-\alpha z)^2}{2 \sigma^2z}\right).
\end{eqnarray}
We can follow the same argument as in the proof of Lemma \ref{lem:l2dist1} and verify that $I_1^{(k)}$ and $I_3^{(k)}$  are bounded above by $C_5(\sqrt{k})^{-1}$. To handle the second term, note that by (\ref{eqnarray:partial}), we have
\begin{eqnarray*}
\left\lvert\frac{\partial f}{\partial z}(k,z)\right\rvert \le C_6(\sqrt{k})^{-1}.
\end{eqnarray*}
 on $\lfloor\frac{k}{\alpha}\rfloor-L_k \le z \le \lfloor\frac{k}{\alpha}\rfloor+L_k$. So,
 \begin{eqnarray*}
I_2^{(k)} \le \int_{\lfloor\frac{k}{\alpha}\rfloor-L_k}^{\lfloor\frac{k}{\alpha}\rfloor+L_k}\left(C_6(\sqrt{k})^{-1}\right)^2dz = 2C_6\frac{\sqrt{\log k}}{\sqrt{k}}.
\end{eqnarray*}
This gives the first bound claimed in the lemma. To prove the second bound, we proceed exactly the same way, but now we use $$\left|\frac{\partial}{\partial z}\frac{1}{\sigma \sqrt{z}}\phi\left(\frac{k-\alpha z}{\sigma\sqrt{z}}\right)\right| \le C_7 k^{-1}$$ on $\lfloor\frac{k}{\alpha}\rfloor-L_k \le z \le \lfloor\frac{k}{\alpha}\rfloor+L_k$, in place of the derivative bound on $f$.\qed\\\\
So, by the preceding lemmas, we have proved that $$\mathbb{E}|X_{k+1}(n)-u(k,n)|^2 \le C_{\epsilon}'\frac{\sqrt{\log k}}{\sqrt{k}}$$ and $$\mathbb{E}|(X_{k+2}(n)-X_{k+1}(n))-\partial_x u(k,n)|^2 \le C_{\epsilon}\frac{(\log k)^{7/2}}{k^{3/2}}$$ uniformly on $n^{\epsilon}\le k \le n$.\\\\
Now, $X_{k+1}(n)- u(k,n) = \int_0^n(g(k,z)-f(k,z))dW^n_z$. As this is a centred Gaussian random variable,
\begin{eqnarray}\label{eqnarray:momentbound}
\mathbb{E}|X_{k+1}(n)- u(k,n)|^{2p} &\le& C_p \left[\int_0^n(g(k,z)-f(k,z))^2dz\right]^p \nonumber\\
&\le & C_p  C_{\epsilon}'^p\left(\frac{\log k}{k}\right)^{p/2},
\end{eqnarray}
and similarly,
\begin{eqnarray}\label{eqnarray:momentboundder}
\mathbb{E}|(X_{k+2}(n)-X_{k+1}(n))-\partial_x u(k,n)|^{2p} \le C_p  C_{\epsilon}'^p\left(\frac{(\log k)^7}{k^3}\right)^{p/2}.
\end{eqnarray}
We use this to obtain the following theorem.
\begin{thm}\label{thm:conga}
For any $\mu > 0$, $\epsilon > 0$ and $ \eta >0$, define the following events:
$$\displaystyle{A_{k,n}=\left\lbrace|X_{k+1}(n)-u(k,n)| > \mu k^{-\left(\frac{1}{4}- \eta\right)} \text{ or } |(X_{k+2}(n)-X_{k+1}(n))-\partial_x u(k,n)|> \mu k^{-\left(\frac{3}{4}- \eta\right)}\right\rbrace}$$ and $$\displaystyle{B_n^{\epsilon}= \bigcup_{k=n^{\epsilon}}^n A_{k,n}}$$
Then $P(\limsup_nB_n^{\epsilon})=0$.
\end{thm}
\textbf{Proof: }By Chebychev type Inequality (for $2p$-th moment) and (\ref{eqnarray:momentbound}), we get for any $p \ge 1$,
\begin{eqnarray*}
P(A_{k,n})&\le & (\mu k^{-(\frac{1}{4}- \eta)})^{-2p}\mathbb{E}|X_{k+1}(n)-u(k,n)|^{2p}\\
&\quad & + (\mu k^{-(\frac{3}{4}- \eta)})^{-2p}\mathbb{E}|(X_{k+2}(n)-X_{k+1}(n))-\partial_x u(k,n)|^{2p}\\
&\le &  \frac{C_p  C_{\epsilon}'^p}{\mu^{2p}}k^{-2\eta p}\left[k^{p/2}\left(\frac{\log k}{k}\right)^{p/2} + k^{3p/2}\left(\frac{(\log k)^7}{k^3}\right)^{p/2}\right] \le D_{\epsilon ,p} \ k^{-(2\eta p-1)}.
\end{eqnarray*}
Hence,
\begin{eqnarray*}
P(B_n^{\epsilon}) \le D_{\epsilon ,p} \sum_{k=n^{\epsilon}}^nk^{-(2\eta p-1)} \le \frac{D'_{\epsilon ,p}}{2\eta p -2} \ n^{-\epsilon(2\eta p -2)}.
\end{eqnarray*} 
Now, choose $p$ large enough such that $P(B_n^{\epsilon})\le C n^{-2}$ for some constant $C$. The result now follows by the Borel Cantelli lemma. \qed\\\\
\textbf{Note: }The above theorem suggests that the distance between the points $X_{k+1}(n)$ and $u(k,n)$ decreases as we increase $k$. This is exactly what is suggested by Figure 2. Furthermore, the difference between $X_{k+2}(n)-X_{k+1}(n)$ and $\partial_x u(k,n)$ decreases even more rapidly on increasing $k$.

After we have investigated the path properties of the continuous Conga line $u$ in subsequent sections, we will prove a fluctuation estimate for $u$ in Lemma \ref{lem:fluctuation} which can be combined with Theorem \ref{thm:conga} to prove Theorem \ref{thm:congafl} which states that the linear interpolation of $X(\cdot,n)$ and its derivative (which exists everywhere except integer points) come uniformly close to $(u(\cdot,n), \partial_x u(\cdot,n))$. This will formally explain why the discrete Conga line \textit{looks smooth} when observed sufficiently far away from the tip.
\subsection{Another smooth approximation}\label{subsection:related}
Here we give another smooth approximation $\overline{u}$ to $X$ given by
\begin{equation}\label{equation:baru}
\displaystyle{\baru(x,t)=\mathbb{E}_{Z_{\rho^2x}}W\left(t-\frac{x}{\alpha}-Z_{\rho^2x}\right)},
\end{equation}
where $\rho=\sigma/\alpha^{3/2}$, $Z_{\rho^2x}$ is a normal random variable with mean zero and variance $\rho^2x$ and $\mathbb{E}_{Z_{\rho^2x}}$ represents expectation taken with respect to $Z_{\rho^2x}$ for a fixed realisation of $W$ (see Notation (i)).

This approximation is more convenient and intuitive for investigating the paths of individual particles and studying the phenomenon of freezing near the origin. Note that $\baru$ has the following properties:\\
$\bullet$ The curve $\{\baru(x,t): x \in (0,t]\}$ for fixed $t$ corresponds to the \textit{Conga line}. Increasing $x$ (i.e. moving away from the tip) results in \textit{increasing smoothness along the same curve} $\baru(\cdot,t)$, indicated by the increasing variance of the $Z_.$ variable. (See Figure 2)\\
$\bullet$ The curve $\{\baru(x,t): t \in [x,\infty)\}$ for fixed $x$ corresponds to the \textit{path of the particle} at distance $x$ from the tip. As successive particles `cut corners', the \textit{family of curves $\baru(x,\cdot)$ become progressively smoother with $x$}. (See Figure 5)
\begin{thm}
The result of Theorem \ref{thm:conga} holds with $u$ replaced by $\baru$.
\end{thm}
\textbf{Proof: }Let $\rho=\sigma/\alpha^{3/2}$. Consider another continuous process $u^*$ given by $$u^*(x,t)=\int_0^t\barphi\left(\frac{\frac{x}{\alpha}-s}{\rho\sqrt{x}}\right)dW^t_s.$$ Rewrite $u^*$ as follows:
\begin{eqnarray*}
u^*(x,t)&=&\int_0^t\int_{-\infty}^s\frac{1}{\sqrt{ x}\rho}\phi\left(\frac{\frac{x}{\alpha}-w}{\rho\sqrt{x}}\right)dwdW^t_s\\\\
&=&\int_{-\infty}^t\int_{\max\{w,0\}}^tdW^t_s\frac{1}{\sqrt{x}\rho}\phi\left(\frac{\frac{x}{\alpha}-w}{\rho\sqrt{x}}\right)dw\\\\
&=&\int_0^t\frac{1}{\sqrt{x}\rho}\phi\left(\frac{\frac{x}{\alpha}-w}{\rho\sqrt{x}}\right)W(t-w)dw+W(t)\barphi\left(\frac{\sqrt{x}}{\alpha\rho}\right)\\\\
&=&\int_{-\infty}^t\frac{1}{\sqrt{x}\rho}\phi\left(\frac{\frac{x}{\alpha}-w}{\rho\sqrt{x}}\right)W(t-w)dw\\\\
&-&\int_{-\infty}^0\frac{1}{\sqrt{x}\rho}\phi\left(\frac{\frac{x}{\alpha}-w}{\rho\sqrt{x}}\right)W(t-w)dw
+W(t)\barphi\left(\frac{\sqrt{x}}{\alpha\rho}\right)\\\\
&=&\baru(x,t)+ +e_1(x,t)+e_2(x,t).
\end{eqnarray*}
Clearly, $$\displaystyle{\mathbb{E}(|e_2(k,n)|^2) \le \frac{\alpha^2\rho^2n}{k}\exp\left\lbrace-\frac{k}{\alpha^2\rho^2}\right\rbrace}.$$ So, on $n^{\epsilon} \le k \le n$, $$\mathbb{E}(|e_2(k,n)|^2) \le \alpha^2\rho^2n\exp\left\lbrace-\frac{n^{\epsilon}}{\alpha^2\rho^2}\right\rbrace.$$
Furthermore, by using stochastic integration by parts to express $e_1(x,t)$ as a stochastic integral and then computing the second moment, it can be verified that
\begin{eqnarray*}
\mathbb{E}(|e_1(x,t)|^2)=t\barphi^2\left(\frac{\sqrt{x}}{\alpha\rho}\right)+\int_0^{\infty}\barphi^2\left(\frac{s+\frac{x}{\alpha}}{\rho\sqrt{x}}\right)ds.
\end{eqnarray*} 
So, on $n^{\epsilon} \le k \le n$,
\begin{eqnarray*}
\mathbb{E}(|e_1(k,n)|^2)&\le& n\barphi^2\left(\frac{n^{\epsilon/2}}{\rho \alpha}\right) + \int_0^{\infty}\barphi^2\left(\frac{s}{\rho\sqrt{n}}+\frac{n^{\epsilon/2}}{\rho \alpha}\right)ds\\
&\le &\frac{\alpha^2\rho^2}{n^{\epsilon}}\exp\left\lbrace-\frac{n^{\epsilon}}{\alpha^2\rho^2}\right\rbrace+\rho\sqrt{n}\int_{\frac{n^{\epsilon/2}}{\rho \alpha}}^{\infty}\frac{1}{y^2}e^{-y^2}dy\\
&\le &\frac{\alpha^2\rho^2}{n^{\epsilon}}\exp\left\lbrace-\frac{n^{\epsilon}}{\alpha^2\rho^2}\right\rbrace +\frac{\rho\sqrt{n}}{2\left(\frac{n^{\epsilon/2}}{\rho \alpha}\right)^3}\exp\left\lbrace\frac{-n^{\epsilon}}{\rho^2\alpha^2}\right\rbrace.
\end{eqnarray*}
By calculations similar to those in the proof of Lemma \ref{lem:l2dist2},
\begin{eqnarray*}
\mathbb{E}|u(k,n)-u^*(k,n)|^2&=&\int_0^n\left[\barphi\left(\frac{k-\alpha s}{\sigma\sqrt{s}}\right)-\barphi\left(\frac{k-\alpha s}{\alpha\rho\sqrt{k}}\right)\right]^2ds\\
&\le & C\frac{(\log k)^{\frac{5}{2}}}{\sqrt{k}}.
\end{eqnarray*}
It is routine to check that a similar analysis yields the required bound for $\mathbb{E}|\partial_xu(k,n)-\partial_x\baru(k,n)|^2$.
Now, proceeding exactly as in the proof of Theorem \ref{thm:conga}, we get the result.\qed
\section{\textbf{The continuous one dimensional Conga Line}}\label{section:continuous}
Here we investigate properties of the continuous one dimensional Conga line $u$, the \textit{coordinate process} of the random smooth curve obtained in the previous section as an approximation to the discrete Conga line $X$.
\subsection{Particles moving at different scales}
It is not hard to observe by estimating $$\displaystyle{\operatorname{Var}\left(u(x,t)\right)=\int_0^t\barphi^2\left(\frac{x-\alpha y}{\sigma \sqrt{y}}\right)dy}$$ that particles at distances $ct$ from the leading particle have variance $O(t)$ if $c < \alpha$ and $o(1)$ (in fact, the variance decays exponentially with $t$) if $c > \alpha$. Also in a window of width $c\sqrt{t}$ about $\alpha t$, the variance is $O(\sqrt{t})$. In particular, this indicates that there is a window of the form $[\alpha t, \alpha t + c_t]$, with $\frac{c_t}{t} \rightarrow 0$ and $\frac{c_t}{\sqrt{t}} \rightarrow \infty$, where the variance changes from being 'very large to very small'. Furthermore, we will show that there is a \textit{cut-off} around which the variance shows a sharp transition: below it the variance grows to infinity with $t$, and above it the variance decreases to zero.
\begin{thm}\label{thm:variance}
\hspace{0.01cm}
\begin{enumerate}
\item[(i)] For $\lambda > 0$ and $1/2 \le \beta < 1$, $\operatorname{Var}(u(\alpha t - \lambda t^{\beta},t))\sim t^{\beta}$.
\item[(ii)] $\operatorname{Var}(u(\alpha t + \sigma\sqrt{\lambda t \log t},t)) \sim \displaystyle{\frac{t^{1/2-\lambda}}{(\log t)^{3/2}}}$.
\end{enumerate}
\end{thm}
\textbf{Proof: }(i) Take any $c>\lambda/\alpha$. Then, decomposing the variance,
\begin{eqnarray*}
\operatorname{Var}(u(\alpha t - \lambda t^{\beta},t))=\int_0^t\barphi^2\left(\frac{\alpha t - \lambda t^{\beta}-\alpha y}{\sigma \sqrt{y}}\right)dy &=& \int_0^{t-ct^{\beta}}\barphi^2\left(\frac{\alpha t - \lambda t^{\beta}-\alpha y}{\sigma \sqrt{y}}\right)dy\\\\
&+&\int_{t-ct^{\beta}}^t\barphi^2\left(\frac{\alpha t - \lambda t^{\beta}-\alpha y}{\sigma \sqrt{y}}\right)dy.
\end{eqnarray*}
The first integral satisfies
\begin{eqnarray*}
\int_0^{t-ct^{\beta}}\barphi^2\left(\frac{\alpha t - \lambda t^{\beta}-\alpha y}{\sigma \sqrt{y}}\right)dy &\le &Ct^{1/2}\int_0^{t-ct^{\beta}}\barphi^2\left(\frac{\alpha t - \lambda t^{\beta}-\alpha y}{\sigma \sqrt{y}}\right)\left(\frac{\alpha t - \lambda t^{\beta}+\alpha y}{2\sigma y\sqrt{y}}\right)dy\\
&\le & Ct^{1/2}\int_{\frac{\alpha c - \lambda}{\sigma}t^{\beta-\frac{1}{2}}}^{\infty}\barphi^2(z)dz\\
&\le & Ct^{1/2}\exp\left(-\frac{(\alpha c -\lambda)^2t^{2\beta-1}}{\sigma^2}\right),
\end{eqnarray*}
where the second step above follows from a change of variables.\\\\
It is easy to check that $$\int_{t-ct^{\beta}}^t\barphi^2\left(\frac{\alpha t - \lambda t^{\beta}-\alpha y}{\sigma \sqrt{y}}\right)dy \sim t^{\beta}$$ proving part (i).\\\\
(ii) We decompose the variance as $$\operatorname{Var}(u(\alpha t + \sigma\sqrt{\lambda t \log t},t))=\int_0^{t/2}\barphi^2\left(\frac{\alpha t + \sigma\sqrt{\lambda t \log t}-\alpha y}{\sigma \sqrt{y}}\right)dy+\int_{t/2}^t\barphi^2\left(\frac{\alpha t + \sigma\sqrt{\lambda t \log t}-\alpha y}{\sigma \sqrt{y}}\right)dy.$$ The first integral decays like $e^{-Ct}$ for some constant $C$. For the second integral, we make a change of variables similar to (i) and standard estimates for the normal c.d.f. to get
\begin{eqnarray*}
\int_{t/2}^t\barphi^2\left(\frac{\alpha t + \sigma\sqrt{\lambda t \log t}-\alpha y}{\sigma \sqrt{y}}\right)dy &\sim &  t^{1/2}\int_{\sqrt{\lambda \log t}}^{\infty}\barphi^2\left(z\right)dz \sim  t^{1/2}\int_{\sqrt{\lambda \log t}}^{\infty}z^{-2}\exp(-z^2)dz\\
&\sim & \frac{t^{1/2}}{(\log t)^{3/2}}\int_{\sqrt{\lambda \log t}}^{\infty}2z\exp(-z^2)dz \sim \frac{t^{1/2-\lambda}}{(\log t)^{3/2}},
\end{eqnarray*}
proving (ii).
\qed\\\\
Part (ii) of the above theorem has the following interesting consequence, demonstrating a \textit{cut-off phenomenon} for the variance of $u(x,t)$ in the vicinity of $\displaystyle{x=\alpha t + \sigma\sqrt{\frac{1}{2}t \log t}}$.
\begin{cor}\label{cor:vartrans}
As $t \rightarrow \infty$\\\\
 $(i)$ $\operatorname{Var}(u(\alpha t + \sigma\sqrt{\lambda t \log t},t)) \rightarrow 0$ if $\lambda \ge 1/2$.\\\\
 $(ii)$ $\operatorname{Var}(u(\alpha t + \sigma\sqrt{\lambda t \log t},t)) \rightarrow \infty$ if $\lambda < 1/2$.\\\\
 $(iii)$ For $0 < \delta < \infty$, $\operatorname{Var}(u(\alpha t + \sigma\sqrt{t((1/2) \log t - (3/2)\log \log t - \log \delta)},t)\sim \delta$.\\\\
 So, the variance exhibits a sharp transition around $\displaystyle{\alpha t + \sigma\sqrt{t((1/2) \log t - (3/2)\log \log t)}}$.
 \end{cor}
 The proof follows easily from part (ii) of Theorem \ref{thm:variance}.
 \subsection{Analyticity of the scaled Conga Line}
For a fixed time $t$, by a change of variables in (\ref{eqnarray:SIP}), we have:
\begin{equation}\label{equation:scalingone}
t^{-\frac{1}{2}}u(tx,t)=\int_0^1W^{(t)}(1-z)(\sqrt{2\pi})^{-1}\left(\frac{x+\alpha z}{2\sigma_t z^{3/2}}\right)\exp\left(-\frac{(x-\alpha z)^2}{2 \sigma_t^2z}\right)dz
\end{equation}
where $\displaystyle{\sigma_t=\frac{\sigma}{\sqrt{t}}}$ and $W^{(t)}(z)=t^{-\frac{1}{2}}W(tz), 0 \le z \le 1$. 
So, to study the Conga line for fixed $t$, we study the scaled process $$\displaystyle{u^W_t(x)=\int_0^1W(1-z)(\sqrt{2\pi})^{-1}\left(\frac{x+\alpha z}{2\sigma_t z^{3/2}}\right)\exp\left(-\frac{(x-\alpha z)^2}{2 \sigma_t^2z}\right)dz=\int_0^1\barphi\left(\frac{x-\alpha s}{\sigma_t \sqrt{s}}\right)dW^1_s}$$ for $0 \le x \le 1$. From (\ref{equation:scalingone}), note that $u$ and $u^W_t$ are connected by the exact equality
\begin{equation}\label{equation:exeq}
t^{-\frac{1}{2}}u(tx,t)=u^{W^{(t)}}_t(x)
\end{equation}
for $0 \le x \le 1$. In particular, for fixed $t$,
\begin{equation}\label{equation:scaling}
\{t^{-\frac{1}{2}}u(xt,t): 0 \le x \le 1\} \eqd \{u^W_t(x): 0 \le x \le 1\}.
\end{equation}
When the driving Brownian motion $W$ is clear from the context, we will suppress the superscript $W$ and write $u_t$ for $u^W_t$.\\\\
Now, we take a look at the derivatives of $u_t$. It is easy to check that we can differentiate under the integral. Thus,
\begin{eqnarray*}
\partial_xu_t(x) &=& -\int_0^1\frac{1}{\sigma_t \sqrt{s}}\phi\left(\frac{x-\alpha s}{\sigma_t \sqrt{s}}\right)dW^1_s,\\
\partial_x^{2}u_t(x) &=& -\int_0^1\frac{1}{\sigma_t^2 s}\phi'\left(\frac{x-\alpha s}{\sigma_t \sqrt{s}}\right)dW^1_s\\
&=&\int_0^1\frac{1}{\sigma_t^2 s}\left(\frac{x-\alpha s}{\sigma_t \sqrt{s}}\right)\phi\left(\frac{x-\alpha s}{\sigma_t \sqrt{s}}\right)dW^1_s.
\end{eqnarray*}
 
In general, the $(n+1)$th derivative takes the following form:
\begin{eqnarray*}
\partial_x^{n+1}u_t(x) &=& \int_0^1(\sigma_t \sqrt{s})^{-(n+1)}(-1)^{n+1}\operatorname{He}_n\left(\frac{x-\alpha s}{\sigma_t \sqrt{s}}\right)\phi\left(\frac{x-\alpha s}{\sigma_t \sqrt{s}}\right)dW^1_s,
\end{eqnarray*}
where $\operatorname{He}_n$ is the $n$-th \textit{Hermite polynomial (probabilist version)} given by $$\operatorname{He}_n(x)=(-1)^ne^{x^2/2}\operatorname{\frac{d^n}{dx^n}}e^{-x^2/2}.$$
In the following lemma, we give an upper bound on the growth rate of the derivatives. Using this, we will prove that, for fixed $t$, $u_t$ is \textit{real analytic} on the interval $(0,1)$, and the radius of convergence around $x_0$ is comparable to $|x_0|$. This is natural as the Conga line gets smoother as we move away from the tip. We start off with the following lemma.
\begin{lem}\label{lem:derivative}
For $\displaystyle{0<\epsilon < \frac{x}{\alpha}}$,
\begin{eqnarray*}
\left|\partial_x^{n+1}u_t(x)\right| &\le & (2\pi)^{1/4}\lVert W \rVert\left\lbrace\left(\frac{\sqrt{2}}{x-\alpha\epsilon}\right)^{n+1}+ \frac{1}{x}\left(\frac{\sqrt{2}}{x-\alpha\epsilon}\right)^{n}\right\rbrace(n+1)!\\
&+& \lVert W \rVert\left\lbrace\left(\frac{1}{\sigma_t\sqrt{\epsilon}}\right)^{n+1}+\frac{n+1}{x}\left(\frac{1}{\sigma_t\sqrt{\epsilon}}\right)^n\right\rbrace\sqrt{(n+1)!},
\end{eqnarray*}
 where $\displaystyle{\lVert W \rVert = \sup_{0\le s\le 1}|W_s|}$.
\end{lem}
\textbf{Proof: } In the proof, we consider $C$ as a generic positive constant whose value might change in between steps.

Let
\begin{equation}\label{equation:K}
K_t^n(x,s)=(\sigma_t \sqrt{s})^{-(n+1)}(-1)^{n+1}\operatorname{He}_n\left(\frac{x-\alpha s}{\sigma_t \sqrt{s}}\right)\phi\left(\frac{x-\alpha s}{\sigma_t \sqrt{s}}\right)
\end{equation}
 Then 
\begin{eqnarray}\label{eqnarray:derivativeK}
\partial_x^{n+1}u_t(x)&=&\int_0^1(\sigma_t \sqrt{s})^{-(n+1)}(-1)^{n+1}\operatorname{He}_n\left(\frac{x-\alpha s}{\sigma_t \sqrt{s}}\right)\phi\left(\frac{x-\alpha s}{\sigma_t \sqrt{s}}\right)dW^1_s\nonumber\\
&=& \int_0^1W(1-s)\partial_sK_t^n(x,s)ds.
\end{eqnarray}
So, $\displaystyle{\left|\partial_x^{n+1}u_t(x)\right| \le ||W||\int_0^1\left|\partial_sK_t^n(x,s)\right|ds}$. So, we have to estimate the integral $\displaystyle{\int_0^1\left|\partial_sK_t^n(x,s)\right|ds}$. Now,
\begin{eqnarray*}
\partial_sK_t^n(x,s) &=& (-1)^{n}\frac{n+1}{2\sigma_t^{n+1}s^{(n+3)/2}}\operatorname{He}_n\left(\frac{x-\alpha s}{\sigma_t \sqrt{s}}\right)\phi\left(\frac{x-\alpha s}{\sigma_t \sqrt{s}}\right)\\
&\quad & +(-1)^{n+1}\ \frac{x+\alpha s}{2\sigma_t s^{3/2}}(\sigma_t \sqrt{s})^{-(n+1)}\operatorname{He}_{n+1}\left(\frac{x-\alpha s}{\sigma_t \sqrt{s}}\right)\phi\left(\frac{x-\alpha s}{\sigma_t \sqrt{s}}\right).
\end{eqnarray*}
So,
\begin{eqnarray*}
\int_0^1 \left|\partial_sK_t^n(x,s)\right|ds &\le & \frac{n+1}{x}\int_0^1(\sigma_t \sqrt{s})^{-n}\frac{x+\alpha s}{2\sigma_t s^{3/2}}\left|\operatorname{He}_n\left(\frac{x-\alpha s}{\sigma_t \sqrt{s}}\right)\phi\left(\frac{x-\alpha s}{\sigma_t \sqrt{s}}\right)\right|ds \\
&\quad & + \int_0^1(\sigma_t \sqrt{s})^{-(n+1)}\frac{x+\alpha s}{2\sigma_t s^{3/2}}\left|\operatorname{He}_{n+1}\left(\frac{x-\alpha s}{\sigma_t \sqrt{s}}\right)\phi\left(\frac{x-\alpha s}{\sigma_t \sqrt{s}}\right)\right|ds.
\end{eqnarray*}
From the above, it is clear that estimating the second integral suffices.
\begin{equation*}
 \int_0^1(\sigma_t \sqrt{s})^{-(n+1)}\frac{x+\alpha s}{2\sigma_t s^{3/2}}\left|\operatorname{He}_{n+1}\left(\frac{x-\alpha s}{\sigma_t \sqrt{s}}\right)\phi\left(\frac{x-\alpha s}{\sigma_t \sqrt{s}}\right)\right|ds \le I_t^x + J_t^x,
\end{equation*}
where
\begin{eqnarray*}
I_t^x &=& \int_0^{\epsilon}(\sigma_t \sqrt{s})^{-(n+1)}\frac{x+\alpha s}{2\sigma_t s^{3/2}}\left|\operatorname{He}_{n+1}\left(\frac{x-\alpha s}{\sigma_t \sqrt{s}}\right)\phi\left(\frac{x-\alpha s}{\sigma_t \sqrt{s}}\right)\right|ds \\
&=& \int_0^{\epsilon}(x-\alpha s)^{-(n+1)}\left(\frac{x-\alpha s}{\sigma_t \sqrt{s}}\right)^{n+1}\frac{x+\alpha s}{2\sigma_t s^{3/2}}\left|\operatorname{He}_{n+1}\left(\frac{x-\alpha s}{\sigma_t \sqrt{s}}\right)\phi\left(\frac{x-\alpha s}{\sigma_t \sqrt{s}}\right)\right|ds \\
&\le &(x-\alpha\epsilon)^{-(n+1)}\int_0^{\epsilon}\left(\frac{x-\alpha s}{\sigma_t \sqrt{s}}\right)^{n+1}\frac{x+\alpha s}{2\sigma_t s^{3/2}}\left|\operatorname{He}_{n+1}\left(\frac{x-\alpha s}{\sigma_t \sqrt{s}}\right)\phi\left(\frac{x-\alpha s}{\sigma_t \sqrt{s}}\right)\right|ds \\
&= & (x-\alpha\epsilon)^{-(n+1)}\int_{\frac{x-\alpha \epsilon}{\sigma_t \sqrt{\epsilon}}}^{\infty}s^{n+1}\left|\operatorname{He}_{n+1}(s)\right|\phi(s)ds \\
&\le & (x-\alpha\epsilon)^{-(n+1)}\left(\int_{\frac{x-\alpha \epsilon}{\sigma_t \sqrt{\epsilon}}}^{\infty}s^{2n+2}\phi(s)ds\right)^{1/2}\left(\int_{\frac{x-\alpha \epsilon}{\sigma_t \sqrt{\epsilon}}}^{\infty}\operatorname{He}_{n+1}^2(s)\phi(s)ds \right)^{1/2}\\
&\le & (x-\alpha\epsilon)^{-(n+1)}(2\pi)^{1/4}\sqrt{(n+1)!}\left(\int_0^{\infty}s^{2n+2}\phi(s)ds\right)^{1/2}\\
&\le & (2\pi)^{1/4}\left(\frac{\sqrt{2}}{x-\alpha\epsilon}\right)^{n+1}(n+1)!.
\end{eqnarray*}
Here we use the facts that $$\int_{-\infty}^{\infty}\operatorname{He}_n^2(s)\phi(s)ds=n!$$ and $$\displaystyle{\int_0^{\infty}s^{2n+2}\phi(s)ds = \frac{(2n+2)!}{2^{n+1}(n+1)!} < 2^{n+1}(n+1)!}.$$ Similarly,
\begin{eqnarray*}
J_t^x &=& \int_{\epsilon}^1(\sigma_t \sqrt{s})^{-(n+1)}\frac{x+\alpha s}{2\sigma_t s^{3/2}}\left|\operatorname{He}_{n+1}\left(\frac{x-\alpha s}{\sigma_t \sqrt{s}}\right)\phi\left(\frac{x-\alpha s}{\sigma_t \sqrt{s}}\right)\right|ds \\
&\le & (\sigma_t\sqrt{\epsilon})^{-(n+1)}\int_{-\infty}^{\infty}|\operatorname{He}_{n+1}(s)|\phi(s)ds\\
&\le & (\sigma_t\sqrt{\epsilon})^{-(n+1)}\left(\int_{-\infty}^{\infty}\operatorname{He}_{n+1}^2(s)\phi(s)ds\right)^{1/2}\left(\int_{-\infty}^{\infty}\phi(s)ds\right)^{1/2}\\
&= & \sqrt{(n+1)!}(\sigma_t\sqrt{\epsilon})^{-(n+1)}.
\end{eqnarray*}
The lemma follows from the above. \qed\\\\
From this lemma, it is not too hard to see that $u_t$ is real analytic on $(0,1)$. Let $x_0$ be any point in this interval. For $0<\delta < x_0$, define $$\displaystyle{\Delta_t^n(x_0,\delta) = \left(\sup_{x \in (x_0-\delta, x_0+\delta)}|\partial_x^{n+1}u_t(x)|\right)\frac{\delta^{n+1}}{(n+1)!}}.$$ The n-th order Taylor polynomial based at $x_0$ is given by $\displaystyle{T_t^n(x)= \sum_{i=0}^n\frac{\partial_x^i u_t(x_0)}{i!}(x-x_0)^i}$. By Taylor's Inequality, $$\displaystyle{\sup_{x \in (x_0-\delta, x_0+\delta)}|u_t(x)-T_t^n(x)| \le \Delta_t^n(x_0,\delta)}.$$ From the above lemma, we know that, for $\displaystyle{\epsilon < \frac{x_0}{\alpha}}$,
\begin{eqnarray*}\label{eqnarray:error}
\Delta_t^n(x_0,\delta) &\le &(2\pi)^{1/4} \lVert W \rVert\left\lbrace \left(\frac{\sqrt{2}\delta}{x_0-\delta-\alpha\epsilon}\right)^{n+1}+ \frac{\delta}{x_0-\delta}\left(\frac{\sqrt{2}\delta}{x_0-\delta-\alpha\epsilon}\right)^{n}\right\rbrace \\
&\quad &  + \lVert W \rVert\left\lbrace\left(\frac{\delta}{\sigma_t\sqrt{\epsilon}}\right)^{n+1}+\frac{(n+1)\delta}{x_0-\delta}\left(\frac{1}{\sigma_t\sqrt{\epsilon}}\right)^n\right\rbrace\frac{1}{\sqrt{(n+1)!}}.
\end{eqnarray*} 
The above error will go to zero only when $\displaystyle{\frac{\sqrt{2}\delta}{x_0-\delta-\alpha\epsilon} < 1}$, i.e., $ \displaystyle{\delta < \frac{x_0-\alpha\epsilon}{1+\sqrt{2}}}$. We can make $\epsilon$ arbitrarily small to get the following:
\begin{cor}\label{cor:ana}
With probability one, the scaled Conga line $u_t$ is real analytic on $(0,1)$. The power series expansion of $u_t$ around $x_0 \in (0,1)$ converges in $\displaystyle{\left(\frac{\sqrt{2}x_0}{1+\sqrt{2}}, \min\{\sqrt{2}x_0,1\}\right)}$.
\end{cor}
We are going to use this property of the Conga line multiple times in this article.
\subsection{Covariance structure of the derivatives}\label{subsection:covstruct}
In the following sections, we will analyse the finer properties of the Conga line like distribution of critical points, length and shape and number of loops. For all of these, fine estimates on the covariance structure of the derivatives are of utmost importance. This section is devoted to finding these estimates for the one dimensional scaled Conga line $u_t$.

To get uniform estimates for the covariance structure of derivatives, we will look at the scaled Conga line sufficiently away from the tip $x=0$. More precisely, we will consider $\{u_t(x): \delta \le x \le \alpha\}$ for an arbitrary $\delta \in (0, \alpha)$. This amounts to analysing the unscaled version $u(\cdot,t)$ in the region $\delta t \le x \le \alpha t$.

The next lemma is about the covariance between the first derivatives at two points. In what follows, we write $\displaystyle{L_t^x(M)=\alpha^{-1}\sqrt{-M\frac{x}{t}\log \frac{x}{t}}}$.
\begin{lem}\label{lem:derivativecov}
For $\delta \le x,y \le \alpha$, $\displaystyle{\operatorname{Cov}(u_t'(x), u_t'(y))  \ge 0}$ and satisfies
\begin{eqnarray*}
\operatorname{Cov}(u_t'(x), u_t'(y)) = \exp\left\lbrace-\frac{2 \alpha t}{\sigma^2}\left(\sqrt{\frac{x^2+y^2}{2}}-\frac{x+y}{2}\right)\right\rbrace\left(\frac{\sqrt{t}}{2\sqrt{\pi \alpha}\sigma\left(\frac{x^2+y^2}{2}\right)^{1/4}}\right)\left(1 + O^{\infty}\left(\sqrt{\frac{\log t}{t}}\right)\right).
\end{eqnarray*}
Consequently, the correlation function $\displaystyle{\operatorname{\rho}_t(x,y)=\operatorname{Corr}(u_t'(x), u_t'(y))}$ is always non-negative and has the following decay rate
\begin{eqnarray*}
\operatorname{\rho}_t(x,y) \ge \exp\left\lbrace-C_1t(x-y)^2\right\rbrace\left(\frac{(xy)^{1/4}}{\left(\frac{x^2+y^2}{2}\right)^{1/4}}\right)\left(1 + O^{\infty}\left(\sqrt{\frac{\log t}{t}}\right)\right),\\
\operatorname{\rho}_t(x,y) \le \exp\left\lbrace-C_2t(x-y)^2\right\rbrace\left(\frac{(xy)^{1/4}}{\left(\frac{x^2+y^2}{2}\right)^{1/4}}\right)\left(1 + O^{\infty}\left(\sqrt{\frac{\log t}{t}}\right)\right),
\end{eqnarray*}
where constants $C_1,C_2$ depend only on $\delta$ and $\alpha$.
\end{lem}
\textbf{Proof: }To prove this lemma, note that, by completing squares in the exponent, we get
\begin{eqnarray*}
\operatorname{Cov}(u_t'(x), u_t'(y))&=& \int_0^1\frac{1}{\sigma_t^2 s} \phi\left(\frac{x - \alpha s}{\sigma_t\sqrt{s}}\right)\phi\left(\frac{y - \alpha s}{\sigma_t\sqrt{s}}\right)ds\\
&=& \exp\left\lbrace-\frac{2 \alpha}{\sigma_t^2}\left(\sqrt{\frac{x^2+y^2}{2}}-\frac{x+y}{2}\right)\right\rbrace\int_0^1\frac{1}{\sigma_t^2s}\phi^2\left(\frac{\sqrt{\frac{x^2+y^2}{2}} - \alpha s}{\sigma_t\sqrt{s}}\right)ds.
\end{eqnarray*}
Now we want to estimate the integral $\displaystyle{\int_0^1\frac{1}{\sigma_t^2s}\phi^2\left(\frac{x - \alpha s}{\sigma_t\sqrt{s}}\right)ds}$ where $\delta \le x \le \alpha$. By choosing $M$ large enough, we can ensure that $$\int_0^1\frac{1}{\sigma_t^2s}\phi^2\left(\frac{x - \alpha s}{\sigma_t\sqrt{s}}\right)ds=\int_{\frac{x}{\alpha}-L_t^x(M)}^{\frac{x}{\alpha}+L_t^x(M)}\frac{1}{\sigma_t^2s}\phi^2\left(\frac{x - \alpha s}{\sigma_t\sqrt{s}}\right)ds + O^{\infty}\left(\frac{1}{t}\right).$$ Notice that\\
\begin{eqnarray*}
\int_{\frac{x}{\alpha}-L_t^x(M)}^{\frac{x}{\alpha}+L_t^x(M)}\frac{1}{\sigma_t^2s}\phi^2\left(\frac{x - \alpha s}{\sigma_t\sqrt{s}}\right)ds &=& \frac{1}{2\pi}\int_{\frac{x}{\alpha}-L_t^x(M)}^{\frac{x}{\alpha}+L_t^x(M)}\frac{x + \alpha s}{2\sigma_ts^{3/2}}\frac{2\sqrt{s}}{\sigma_t(x + \alpha s)}\exp\left\lbrace-\left(\frac{x - \alpha s}{\sigma_t\sqrt{s}}\right)^2\right\rbrace ds\\
&=&\frac{1}{2\pi\sigma_t\sqrt{\alpha x}}\int_{\frac{x}{\alpha}-L_t^x(M)}^{\frac{x}{\alpha}+L_t^x(M)}\frac{x + \alpha s}{2\sigma_ts^{3/2}}\exp\left\lbrace-\left(\frac{x - \alpha s}{\sigma_t\sqrt{s}}\right)^2\right\rbrace ds\\
&\quad & + \ O^{\infty}(\sqrt{\log t})\\
&=&\frac{1}{2\pi\sigma_t\sqrt{\alpha x}}\int_{-\infty}^{\infty}e^{-s^2}ds + O^{\infty}(\sqrt{\log t})\\
&=&\frac{\sqrt{t}}{2\sigma\sqrt{\pi\alpha x}} + O^{\infty}(\sqrt{\log t}).
\end{eqnarray*}
Substituting $\sqrt{\frac{x^2 +y^2}{2}}$ in place of $x$ proves the lemma.\qed\\\\
\begin{lem}\label{lem:seconddervar}
For $\delta \le x \le \alpha$,
\begin{eqnarray*}
\operatorname{Var}(u_t''(x))
 = \frac{\sqrt{\alpha} \ t^{3/2}}{4\sqrt{\pi}\sigma^3x^{3/2}}\left(1 + O^{\infty}\left(\sqrt{\frac{\log t}{t}}\right)\right).
\end{eqnarray*}
\end{lem}
\textbf{Proof: }Follows along the same lines as the proof of Lemma \ref{lem:derivativecov}.\qed
\begin{lem}\label{lem:onetwoderivative}
For $\delta \le x \le \alpha$,
\begin{eqnarray*}
\operatorname{Cov}(u_t'(x), u_t''(x))= O^{\infty}(\sqrt{t\log t}).
\end{eqnarray*}
\end{lem}
\textbf{Proof: }
\begin{eqnarray*}
\operatorname{Cov}(u_t'(x), u_t''(x))&=& \int_0^1\frac{1}{(\sigma_t\sqrt{s})^3}\left(\frac{x - \alpha s}{\sigma_t\sqrt{s}}\right)\phi^2\left(\frac{x - \alpha s}{\sigma_t\sqrt{s}}\right)ds\\\\
&=&\frac{1}{\pi}\int_{\frac{x}{\alpha}-L_t^x(M)}^{\frac{x}{\alpha}+L_t^x(M)}\frac{x + \alpha s}{2\sigma_ts^{3/2}}\frac{1}{\sigma_t^2(x + \alpha s)}\left(\frac{x - \alpha s}{\sigma_t\sqrt{s}}\right)\exp\left\lbrace-\left(\frac{x - \alpha s}{\sigma_t\sqrt{s}}\right)^2\right\rbrace ds\\
&\quad & + O^{\infty}(t^{-1})\\\\
&=& \frac{1}{\pi}\int_{\frac{x}{\alpha}-L_t^x(M)}^{\frac{x}{\alpha}+L_t^x(M)}\frac{x + \alpha s}{2\sigma_ts^{3/2}}\frac{1}{2\sigma_t^2 x}\left(1+f_t(x,s)\right)\\
&\quad & \hspace{4cm}\times\left(\frac{x - \alpha s}{\sigma_t\sqrt{s}}\right)\exp\left\lbrace-\left(\frac{x - \alpha s}{\sigma_t\sqrt{s}}\right)^2\right\rbrace ds\\
&\quad & + O^{\infty}(t^{-1}),
\end{eqnarray*}
where $\displaystyle{f_t(x,s)=\frac{2x}{x+\alpha s}-1}$.\\\\
Using the fact that $\displaystyle{\int_{-\infty}^{\infty}s \exp(-s^2)ds=0}$, we get
\begin{eqnarray}\label{eqnarray:doobyone}
\frac{1}{\pi}\int_{\frac{x}{\alpha}-L_t^x(M)}^{\frac{x}{\alpha}+L_t^x(M)}\frac{x + \alpha s}{2\sigma_ts^{3/2}}\frac{1}{2\sigma_t^2 x}\left(\frac{x - \alpha s}{\sigma_t\sqrt{s}}\right)\exp\left\lbrace-\left(\frac{x - \alpha s}{\sigma_t\sqrt{s}}\right)^2\right\rbrace ds = O^{\infty}(t^{-1})
\end{eqnarray}
choosing sufficiently large $M$.\\\\
Also note that $$f_t(x,s)=O^{\infty}\left(\sqrt{\frac{\log t}{t}}\right)$$ for $\delta \le x \le \alpha, \ \frac{x}{\alpha}-L_t^x(M) \le s \le \frac{x}{\alpha}+L_t^x(M)$, which yields
\begin{eqnarray}\label{eqnarray:doobytwo}
\frac{1}{\pi}\int_{\frac{x}{\alpha}-L_t^x(M)}^{\frac{x}{\alpha}+L_t^x(M)}\frac{x + \alpha s}{2\sigma_ts^{3/2}}\frac{f_t(x,s)}{2\sigma_t^2 x}
\left(\frac{x - \alpha s}{\sigma_t\sqrt{s}}\right)\exp\left\lbrace-\left(\frac{x - \alpha s}{\sigma_t\sqrt{s}}\right)^2\right\rbrace ds =O^{\infty}(\sqrt{t\log t}).
\end{eqnarray}
(\ref{eqnarray:doobyone}) and (\ref{eqnarray:doobytwo}) prove the lemma.\qed
\begin{cor}\label{cor:correlation}
For $\delta \le x \le \alpha$,
$$\operatorname{Corr}(u_t'(x), u_t''(x))= O^{\infty}\left(\sqrt{\frac{\log t}{t}}\right).$$
\end{cor}
\textbf{Proof: }This follows from Lemmas \ref{lem:derivativecov}, \ref{lem:seconddervar} and \ref{lem:onetwoderivative}.\qed\\\\
Let $\Sigma_t(x,y)$ be the covariance matrix of $(u_t'(x),u_t'(y))$. We need the following technical lemma to estimate the determinant of the matrix. It turns out to be crucial in certain second moment computations in Subsection \ref{subsection:critical}.
\begin{lem}\label{lem:det}
There exist constants $C^*,C_1,C_2$ such that, for $\delta \le x,y \le \alpha$ with $\displaystyle{|x-y| \le \frac{C^*}{\sqrt{t}}}$,
\begin{eqnarray*}
C_1t^2(y-x)^2 \le \det\Sigma_t(x,y)\le C_2t^2(y-x)^2
\end{eqnarray*}
\end{lem}
\textbf{Proof: }We fix $x \in [\delta,\alpha]$ and consider the function $\Psi_{t,x}(y)=\det\Sigma_t(x,y)$. 
Consider the function $\displaystyle{g_t(y)= \operatorname{Var}{u_t'(y)} = \int_0^1\frac{1}{\sigma_t^2 s} \phi^2\left(\frac{y - \alpha s}{\sigma_t\sqrt{s}}\right)ds}$. Let $\operatorname{H}_n$ denote the $n$-th \textit{Hermite polynomial (physicist version)} given by $$\operatorname{H}_n(x)=(-1)^ne^{x^2}\operatorname{\frac{d^n}{dx^n}}e^{-x^2}=2^{n/2}\operatorname{He}_n(\sqrt{2}x).$$ Then we can write the n-th derivative of $g_t$ as
\begin{eqnarray*}
g_t^{(n)}(x) = (-1)^n\int_0^1\frac{1}{\sigma_t^2 s} \frac{1}{(\sigma_t\sqrt{s})^n}\operatorname{H}_n\left(\frac{x - \alpha s}{\sigma_t\sqrt{s}}\right)\phi^2\left(\frac{x - \alpha s}{\sigma_t\sqrt{s}}\right)ds.
\end{eqnarray*}
Using the fact that $\int_{-\infty}^{\infty}H_n(s)\exp\{-s^2\}ds=0$ and the same technique as the proof of Lemma \ref{lem:onetwoderivative}, one can show that, for $n\ge 1$,
\begin{eqnarray}\label{eqnarray:growthder}
g_t^{(n)}(\cdot)=O^{\infty}(t^{n/2}\sqrt{\log t})
\end{eqnarray}
Let $\displaystyle{\eta=\sqrt{\frac{x^2+y^2}{2}}}$.\\\\
Consider the functions $$E_{t,x}(y)=g_t(x)g_t(y)-g_t^2(\eta)$$ and $$F_{t,x}(y)=\left[1-\exp\left\lbrace-\frac{4 \alpha t}{\sigma^2}\left(\sqrt{\frac{x^2+y^2}{2}}-\frac{x+y}{2}\right)\right\rbrace\right]g_t^2(\eta).$$ Then, writing down $\operatorname{Cov}(u_t'(x),u_t'(y))$ as in the proof of Lemma \ref{lem:derivativecov}, we get
\begin{eqnarray}\label{eqnarray:detform}
\Psi_{t,x}(y)=F_{t,x}(y)+E_{t,x}(y).
\end{eqnarray}
It is easy to check that
\begin{eqnarray}\label{eqnarray:onetwozero}
E_{t,x}(x)=E_{t,x}'(x)=0.
\end{eqnarray}
The double derivative of $E_{t,x}$ takes the form $$E_{t,x}''(y)=g_t(x)g_t''(y)-2(g_t'(\eta))^2(\partial_y \eta)^2-2(g_t(\eta))(g_t''(\eta))(\partial_y \eta)^2-2(g_t(\eta))(g_t'(\eta))(\partial_y^2 \eta).$$
Using (\ref{eqnarray:growthder}) we deduce $$E_{t,x}''(\cdot)=O^{\infty}(t^{3/2}\sqrt{\log t}),$$
which, along with (\ref{eqnarray:onetwozero}) yields
\begin{eqnarray}\label{eqnarray:scoobyone}
|E_{t,x}(y)| \le Ct^{3/2}\sqrt{\log t} \ (y-x)^2
\end{eqnarray}
for some constant $C < \infty$.

Now, to estimate $F_{t,x}$, note that in the region $\delta \le x,y \le \alpha$,
\begin{eqnarray}\label{eqnarray:scoobytwo}
\exp\left\lbrace-\frac{\alpha t(y-x)^2}{2\delta\sigma^2}\right\rbrace\le\exp\left\lbrace-\frac{4 \alpha t}{\sigma^2}\left(\sqrt{\frac{x^2+y^2}{2}}-\frac{x+y}{2}\right)\right\rbrace \le \exp\left\lbrace-\frac{t(y-x)^2}{2\sigma^2}\right\rbrace.
\end{eqnarray}
Using (\ref{eqnarray:scoobytwo}) along with the fact that $e^{-C}x \le 1-e^{-x} \le x$ on $0\le x \le C$, and Lemma \ref{lem:derivativecov}, we get
\begin{eqnarray}\label{eqnarray:scoobythree}
C^*t^2(y-x)^2 \le F_{t,x}(y) \le Ct^2(y-x)^2,
\end{eqnarray}
where $C, C^*$ are positive, finite constants.

(\ref{eqnarray:scoobyone}) and (\ref{eqnarray:scoobythree}) together prove the lemma.\qed
\subsection{Analyzing the distribution of critical points}\label{subsection:critical}
Let $\displaystyle{N_t(I)}$ denote the number of critical points of the scaled one dimensional Conga line $u_t$ in an interval $I\subseteq [\delta, \alpha]$. Then $N_t$ defines a \textit{simple point process} on $[\delta, \alpha]$. Our first goal is to find out the \textit{first intensity} of this process. For this, we use the \textit{Expectation meta-theorem for smooth Gaussian fields} (see \cite{adlertaylor} p. 263), which implies the following:
\begin{eqnarray}\label{eqnarray:expmeta}
\mathbb{E}\left(N_t(I)\right)= \int_I \mathbb{E}\left(|u_t''(y)| \ \middle| u_t'(y)=0\right)p_t^y(0)dy,
\end{eqnarray}
where $p_t^y$ is the density of $u_t'(y)$. Before we go further, we remark that the meta-theorem from \cite{adlertaylor} mentioned above is a very general theorem which holds in a much wider set-up under a set of assumptions. In our case, it is easy to check that all the assumptions hold. Now, we utilize (\ref{eqnarray:expmeta}) and the developments in subsection \ref{subsection:covstruct} to derive a nice expression for the first intensity density.
\begin{lem}\label{lem:firstint}
The first intensity density $\rho_t$ for $N_t$ satisfies
\begin{eqnarray}
\rho_t(x)=\frac{\sqrt{\alpha t}}{\pi\sigma\sqrt{2x}}\left(1+O^{\infty}\left(\sqrt{\frac{\log t}{t}}\right)\right)
\end{eqnarray}
for $\delta \le x \le \alpha$.
\end{lem}
\textbf{Proof: }
By standard formulae for normal conditional densities and the lemmas proved in Subsection \ref{subsection:covstruct}, we manipulate (\ref{eqnarray:expmeta}) as follows: 
\begin{eqnarray*}
\mathbb{E}\left(N_t(I)\right)&=& \sqrt{\frac{2}{\pi}}\int_I\left[\frac{\operatorname{Var}(u_t'(y)).\operatorname{Var}(u_t''(y))-\operatorname{Cov}^2(u_t'(y), u_t''(y))}{\operatorname{Var}(u_t'(y))}\right]^{\frac{1}{2}}\frac{1}{\sqrt{2\pi\operatorname{Var}(u_t'(y))}}dy\\\\
&=&\sqrt{\frac{2}{\pi}}\int_I\frac{\left(1+O^{\infty}\left(\frac{\log t}{t}\right)\right)}{\sqrt{2\pi}}\left[\frac{\operatorname{Var}(u_t''(y))}{\operatorname{Var}(u_t'(y))}\right]^{\frac{1}{2}}dy\\\\
&=& \int_I\frac{\left(1+O^{\infty}\left(\sqrt{\frac{\log t}{t}}\right)\right)}{\pi}\left[\frac{\alpha t}{2\sigma^2y}\right]^{\frac{1}{2}}dy\\\\
&=&\int_I\left(1+O^{\infty}\left(\sqrt{\frac{\log t}{t}}\right)\right)\frac{\sqrt{\alpha t}}{\sqrt{2}\pi\sigma}\frac{1}{\sqrt{y}}dy.
\end{eqnarray*}
Here, we use Corollary \ref{cor:correlation} for to get the second equality and the estimates proved in Lemma \ref{lem:derivativecov} and Lemma \ref{lem:seconddervar} to get the third equality. This proves the lemma.\qed\\\\
Thus $\displaystyle{\hat{\rho}_t(x)=\frac{\sqrt{\alpha t}}{\pi\sigma\sqrt{2x}}}$ gives us the approximate first intensity for $N_t$. From this, we see that the expected number of critical points in a small interval $[x,x+h]$ is approximately $\displaystyle{\frac{\sqrt{\alpha t}h}{\pi\sigma\sqrt{2x}}}$.\\\\
Now that we know the first intensity reasonably accurately, we can ask finer questions about the distribution of critical points, such as\\\\
(i) What can we say about the spacings of the critical points? Are there points in $[\delta, \alpha]$ around which there is a large concentration of critical points, or are they more or less well-spaced?\\\\
(ii) Given an interval $I \in [\delta, \alpha]$, how good is $\mathbb{E}N_t(I)$ as an estimate of $N_t(I)$?\\\\
 The next lemma answers (i) by estimating the second intensity of $N_t$. First we present a formula for the second intensity of $N_t$ taken from \cite{adlertaylor}.
\begin{eqnarray}\label{eqnarray:secint}
\mathbb{E}\left(N_t(I)^2-N_t(I)\right)=\int_{I \times I}\mathbb{E}\left(|u_t''(y)||u_t''(z)| \ \middle| u_t'(y)=0,u_t'(z)=0\right)p_t^{y,z}(0)dydz,
\end{eqnarray}
where $\displaystyle{p_t^{y,z}}$ is the joint density of $(u_t'(y),u_t'(z))$.\\\\
In the following, $C^*$ represents a positive constant.
\begin{lem}\label{lem:condsec}
For $\displaystyle{t > \frac{4(1+\sqrt{2})^2}{2\delta^2}}$, and $\delta \le y,z \le \alpha$ with $\displaystyle{|y-z| \le \frac{C^*}{\sqrt{t}}}$, 
\begin{eqnarray*}
\mathbb{E}\left(|u_t''(y)||u_t''(z)|\mid u_t'(y)=0,u_t'(z)=0\right) \le C(y-z)^2t^{5/2}
\end{eqnarray*}
for some constant $C >0$.
\end{lem}
\textbf{Proof: }The hypothesis of the lemma tells us that $y$ and $z$ lie in the region of analyticity of each other, i.e. we can write
\begin{eqnarray*}
u_t'(y)=\sum_{n=0}^{\infty}u_t^{(n+1)}(z)\frac{(y-z)^n}{n!},
\end{eqnarray*}
and the same holds with $y$ and $z$ interchanged. If we know that $u_t'(y)=0$ and $u_t'(z)=0$, the above equation becomes
\begin{eqnarray*}
0=\sum_{n=1}^{\infty}u_t^{(n+1)}(z)\frac{(y-z)^n}{n!}.
\end{eqnarray*}
From this, we can solve for $u_t''(z)$ to get
\begin{eqnarray*}
u_t''(z)=-\sum_{n=2}^{\infty}u_t^{(n+1)}(z)\frac{(y-z)^{n-1}}{n!},
\end{eqnarray*}
and the same holds with $y$ and $z$ interchanged.
Thus, the conditional expectation in (\ref{eqnarray:secint}) becomes\\\\
$\displaystyle{\mathbb{E}\left[\left|-\sum_{n=1}^{\infty}u_t^{(n+2)}(z)\frac{(y-z)^n}{(n+1)!}\right|\left|-\sum_{n=1}^{\infty}u_t^{(n+2)}(y)\frac{(y-z)^n}{(n+1)!}\right| \ \bigg| u_t'(y)=0,u_t'(z)=0 \right]}\\\\
=\displaystyle{(y-z)^2\mathbb{E}\left[\left|-\sum_{n=0}^{\infty}u_t^{(n+3)}(z)\frac{(y-z)^{n}}{(n+2)!}\right|\left|-\sum_{n=0}^{\infty}u_t^{(n+3)}(y)\frac{(y-z)^n}{(n+2)!}\right| \ \bigg| u_t'(y)=0,u_t'(z)=0 \right]}.$\\\\
Now, by the Cauchy-Schwarz inequality and the fact that the conditional variance is bounded above by the total variance, we have\\\\
$\mathbb{E}\left(|u_t^{(m)}(y)||u_t^{(n)}(z)| \ \big| u_t'(y)=0,u_t'(z)=0\right)\\\\
\hspace{2cm}\le\sqrt{\mathbb{E}\left(u_t^{(m)}(y)^2 \ \big| u_t'(y)=0,u_t'(z)=0\right)}\sqrt{\mathbb{E}\left(u_t^{(n)}(z)^2 \ \big| u_t'(y)=0,u_t'(z)=0\right)}\\\\
\hspace{2cm}\le\sqrt{\mathbb{E}\left(u_t^{(m)}(y)\right)^2}\sqrt{\mathbb{E}\left(u_t^{(n)}(z)\right)^2}.$\\\\
We know that
\begin{eqnarray*}
u_t^{(m+3)}(y) &=& \int_0^1(\sigma_t \sqrt{s})^{-(m+3)}(-1)^{m+3}\operatorname{He}_{m+2}\left(\frac{y-\alpha s}{\sigma_t \sqrt{s}}\right)\phi\left(\frac{y-\alpha s}{\sigma_t \sqrt{s}}\right)dW^1_s.
\end{eqnarray*}
So, using the same techniques as in the proof of Lemma \ref{lem:derivative}, for $\displaystyle{0<\epsilon<\frac{\delta}{2\alpha}}$, we estimate the variance as
\begin{eqnarray*}
\mathbb{E}\left(u_t^{(m+3)}(y)\right)^2&=&\int_0^1(\sigma_t \sqrt{s})^{-2(m+3)}\operatorname{He}_{m+2}^2\left(\frac{y-\alpha s}{\sigma_t \sqrt{s}}\right)\phi^2\left(\frac{y-\alpha s}{\sigma_t \sqrt{s}}\right)ds\\
&\le & \frac{C_1t}{\sigma^2(y-\alpha\epsilon)^{2m+4}}\int_0^{\infty}s^{2m+3}\operatorname{He}_{m+2}^2(s)\phi^2(s)ds + \frac{C_2}{(\sigma_t\sqrt{\epsilon})^{2m+5}}\int_0^{\infty}\operatorname{He}_{m+2}^2(s)\phi^2(s)ds.
\end{eqnarray*}
To estimate the first integral, we note that the function $$\displaystyle{g_m(s)=s^{2m+3}\exp\{-s^2/2\}}$$ is maximised at $s=\sqrt{2m+3}$. So, $\displaystyle{g_m(s) \le (2m+3)^{(2m+3)/2}\exp\{-(2m+3)/2\}}$.\\\\
By Stirling's Formula,
\begin{eqnarray*}
\Gamma(n)=\sqrt{\frac{2\pi}{n}}\left(\frac{n}{e}\right)^n\left(1+O(n^{-1})\right),
\end{eqnarray*}
where $\Gamma(\cdot)$ is the Gamma function.
Using this, we get $\displaystyle{g_m(s) \le C\sqrt{m}2^{m}(m+1)!}$. So,
\begin{eqnarray*}
\int_0^{\infty}s^{2m+3}\operatorname{He}_{m+2}^2(s)\phi^2(s)ds \le C\sqrt{m}2^{m}(m+1)!\int_{-\infty}^{\infty}\operatorname{He}_{m+2}^2(s)\phi(s)ds \le C2^m\{(m+2)!\}^2.
\end{eqnarray*}
The second integral is easier to estimate. Finally, we get
\begin{eqnarray}\label{eqnarray:higherdervar}
\mathbb{E}\left(u_t^{(m+3)}(y)\right)^2 \le \frac{C_1t}{\sigma^2(y-\alpha\epsilon)^{2m+4}}2^m\{(m+2)!\}^2 + \frac{C_2}{(\sigma_t\sqrt{\epsilon})^{2m+5}}\{(m+2)!\}.
\end{eqnarray}
Therefore,
\begin{eqnarray*}
\frac{\sqrt{\mathbb{E}\left(u_t^{(m+3)}(y)\right)^2}}{(m+2)!}|y-z|^m &\le & C \left[t\left(\frac{\sqrt{2}|y-z|}{y-\alpha \epsilon}\right)^{2m}\left(\frac{1}{y-\alpha \epsilon}\right)^4\right.\\
&\quad & \left. + \ \frac{1}{(m+2)!}\left(\frac{\sqrt{t}|y-z|}{\sigma\sqrt{\epsilon}}\right)^{2m}\left(\frac{\sqrt{t}}{\sigma\sqrt{\epsilon}}\right)^5\right]^{1/2}\\
&\le & Ct^{5/4}a_m(t,y,z),
\end{eqnarray*} 
where, by the assumptions of the lemma, $\displaystyle{\sum_{m=0}^{\infty}a_m(t,y,z) \le C}$, where $C$ is a constant that does not depend on $t,y,z$. Thus, we have
\begin{eqnarray*}
\mathbb{E}\left(|u_t''(y)||u_t''(z)| \ \big| u_t'(y)=0,u_t'(z)=0\right)&\le &(y-z)^2\sum_{m,n=0}^{\infty}\frac{\sqrt{\mathbb{E}\left(u_t^{(m+3)}(y)\right)^2}}{(m+2)!}|y-z|^m\\
&\quad & \hspace{2.5cm} \times\frac{\sqrt{\mathbb{E}\left(u_t^{(n+3)}(y)\right)^2}}{(n+2)!}|y-z|^n\\\\
&\le &C(y-z)^2t^{5/2}\left(\sum_{m=0}^{\infty}a_m(t,y,z)\right)^2,
\end{eqnarray*}
which proves the lemma.\qed\\\\
We know that $\displaystyle{p_t^{y,z}(0)=\frac{1}{2\pi\sqrt{\det\Sigma_t(y,z)}}}$. Using Lemmas \ref{lem:det} and \ref{lem:condsec}, we get
\begin{lem}\label{lem:secondintensity}
For $\displaystyle{t > \frac{4(1+\sqrt{2})^2}{2\delta^2}}$, and $h \le C^* t ^{-1/2}$,
\begin{eqnarray*}
\mathbb{E}(N_t^2([x,x+h])-N_t([x,x+h])) \le C(\delta)h^3t^{3/2}.
\end{eqnarray*}
In particular, we get
\begin{eqnarray*}
\mathbb{E}(N_t([x,x+h])\mathbb{I}(N_t([x,x+h]) \ge 2)) \le C(\delta)h^3t^{3/2}.
\end{eqnarray*}
\end{lem}
Using this lemma, we can deduce that if we divide $[\delta,\alpha]$ into subintervals of sufficiently small length, the number of critical points in any of these should not exceed one. The following corollary makes this precise.
\begin{cor}
Let $\{a_t\}$ be any sequence such that $a_t=o(t^{-1/4})$. Divide the interval $[\delta,\alpha]$ into subintervals $I_1,...,I_{[\sqrt{t}/a_t]+1}$ of length at most $\displaystyle{\frac{a_t}{\sqrt{t}}}$. Then
\begin{eqnarray*}
P\left(\max_{1 \le j \le [\sqrt{t}/a_t]+1}N_t(I_j) \ge 2\right) \le C(\delta)a_t^2t^{1/2} \rightarrow 0
\end{eqnarray*}
as $t \rightarrow \infty$.
\end{cor}
This follows easily from Lemma \ref{lem:secondintensity} using the union bound.\\\\
Now, we answer (ii).\\\\
Note that for a Poisson point process, the first intensity determines the whole process. The Conga line lacks the Markov property. We can think of it as a process that `\textit{gains smoothness at the cost of Markov property}'. But Lemma \ref{lem:derivativecov} tells us that there is \textit{exponential decorrelation}, i.e. pieces of the Conga line that are reasonably far apart are almost independent. Thus, we expect that the first intensity of $N_t$ should give us a lot of information about the process $N_t$ itself. We conclude this section on critical points by giving a basis to this intuition. We show the following:
\begin{lem}\label{lem:critprob}
Let $I \subseteq (0, \alpha]$ be a closed interval. Then
\begin{eqnarray*}
\frac{N_t(I)}{\mathbb{E}N_t(I)}\tendp 1
\end{eqnarray*}
as $t \rightarrow \infty$.
\end{lem}
\textbf{Proof: }We prove the result for $I=[\delta, \alpha]$ for an arbitrary $\delta>0$, although the same proof carries over to a general closed interval contained in $(0,\alpha]$.

Consider a collection of intervals $$\mathcal{C}= \{I_j: 1 \le j \le C[\sqrt{t}/r]\}$$ where each interval is of length $\displaystyle{\frac{1}{\sqrt{t}}}$ in $[\delta,\alpha]$, and $\displaystyle{d(I_j,I_k) \ge \frac{r}{\sqrt{t}}}$ for a sufficiently large $r$ (which can be a function of $t$), whose optimal choice will be made later, and $d(A,B)$ represents the usual distance between sets $A$ and $B$. Using the long range independence of the Conga line (see Lemma \ref{lem:derivativecov}), we will prove that $\operatorname{Var}\left(N_t\left(\bigcup_{j=1}^{C[\sqrt{t}/r]}I_j\right)\right)$ is very small compared to $\mathbb{E}\left(N_t\left(\bigcup_{j=1}^{C[\sqrt{t}/r]}I_j\right)\right)^2$. The proof is completed by covering $[\delta,\alpha]$ with $[r]$ translates $\mathcal{C}_1,\cdots,\mathcal{C}_{[r]}$ of such collections and an application of Chebychev Inequality.

Note that all the constants used in this proof depend on $\delta$.\\\\
We begin by computing $\mathbb{E}\left(N_t(I_1)N_t(I_2)\right)$ using an analogue of the Expectation meta-theorem (which can also be derived from the second intensity formula (\ref{eqnarray:secint})) as follows:
\begin{eqnarray}\label{eqnarray:cross}
\mathbb{E}\left(N_t(I_1)N_t(I_2)\right)=\int_{I_1\times I_2}\mathbb{E}\left(|u_t''(y)u_t''(z)|\mid u_t'(y)=0, u_t'(z)=0\right)\frac{1}{2\pi\sqrt{\det\Sigma_t(y,z)}}dydz
\end{eqnarray}
where $\Sigma_t(y,z)$ is the covariance matrix for $(u_t'(y),u_t'(z))$. We know that if $$\left(u_t''(y),u_t''(z),u_t'(y),u_t'(z)\right)\sim N(0,\Sigma),$$ then $$\left(u_t''(y),u_t''(z) \ \big| u_t'(y)=0,u_t'(z)=0\right) \sim N(0,\Sigma^*),$$ where $\Sigma=\begin{bmatrix}
\Sigma_{11}&\Sigma_{12}\\
\Sigma_{21}&\Sigma_{22}
\end{bmatrix}
$ and $\Sigma^*=\Sigma_{11}-\Sigma_{12}\Sigma_{22}^{-1}\Sigma_{21}.$ Now
\begin{eqnarray}\label{eqnarray:cov}
\operatorname{Cov}\left(|u_t''(y)|,|u_t''(z)| \ \big| u_t'(y)=0, u_t'(z)=0\right)&=&\frac{2}{\pi}\sqrt{\sigma^*_{11}\sigma^*_{22}}\left(\rho^*_{12}\arcsin \rho^*_{12} + \sqrt{1-\rho_{12}^{*2}} -1\right)\nonumber\\
& \le & \frac{2}{\pi}\sqrt{\sigma^*_{11}\sigma^*_{22}}\left(\rho^*_{12}\arcsin \rho^*_{12}\right)\nonumber\\
&\le & \sigma^*_{12}=\sigma_{12}-\left(\Sigma_{12}\Sigma_{22}^{-1}\Sigma_{21}\right)_{12}.
\end{eqnarray}
Take $(y,z) \in I_j \times I_k$, where $I_j, I_k \in \mathcal{C}$ with $j\neq k$.

The proof of Lemma \ref{lem:derivativecov} shows that for $\eta=\sqrt{\frac{y^2+z^2}{2}}$, $$\displaystyle{\operatorname{Cov}\left(u_t'(y),u_t'(z)\right)\le \exp\{-Ct(y-z)^2\}\operatorname{Var}(u_t'(\eta))}$$ and $$\operatorname{Var}(u_t'(\eta)) \le C_1\sqrt{t}.$$ So, as $\displaystyle{|y-z|\ge \frac{r}{\sqrt{t}}}$,
\begin{eqnarray}\label{eqnarray:oneone}
\operatorname{Cov}\left(u_t'(y),u_t'(z)\right)&\le & \exp\{-Cr^2\}\operatorname{Var}(u_t'(\eta))\nonumber \\
&\le & C_1\sqrt{t}\exp\{-C_2r^2\}.
\end{eqnarray}
Calculations similar to those in the proof of Lemma \ref{lem:derivativecov} show
\begin{eqnarray*}
\operatorname{Cov}\left(u_t'(y),u_t''(z)\right) &=&-\int_0^1\frac{1}{(\sigma_t\sqrt{s})^3}\left(\frac{z-\alpha s}{\sigma_t\sqrt{s}}\right)\phi\left(\frac{y-\alpha s}{\sigma_t\sqrt{s}}\right)\phi\left(\frac{z-\alpha s}{\sigma_t\sqrt{s}}\right)ds\\
&=&-\exp\left\lbrace-\frac{2\alpha}{\sigma_t^2}\left(\sqrt{\frac{y^2+z^2}{2}}-\frac{y+z}{2}\right)\right\rbrace\\
&\quad & \hspace{1cm}\times\int_0^1\frac{1}{(\sigma_t\sqrt{s})^3}\left(\frac{z-\alpha s}{\sigma_t\sqrt{s}}\right)\phi^2\left(\frac{\sqrt{\frac{y^2+z^2}{2}}-\alpha s}{\sigma_t\sqrt{s}}\right)ds.
\end{eqnarray*}
and thus
\begin{eqnarray}\label{eqnarray:onetwo}
\left|\operatorname{Cov}\left(u_t'(y),u_t''(z)\right)\right| &\le & C_1\exp\{-C_2r^2\}\left(\int_0^1|z-\eta|\frac{1}{(\sigma_t\sqrt{s})^4}\phi^2\left(\frac{\eta-\alpha s}{\sigma_t\sqrt{s}}\right)ds\right.\nonumber\\
&&\quad\left.+\int_0^1\frac{|\eta-\alpha s|}{(\sigma_t\sqrt{s})^4}\phi^2\left(\frac{\eta-\alpha s}{\sigma_t\sqrt{s}}\right)ds\right)\nonumber\\ \nonumber\\
&\le &C_1\exp\{-C_2r^2\}t^{3/2}.
\end{eqnarray}
Also, from Lemma \ref{lem:onetwoderivative},
\begin{eqnarray}
\left|\operatorname{Cov}\left(u_t'(y),u_t''(y)\right)\right| \le C\sqrt{t\log t}.
\end{eqnarray}
Similar calculations also show
\begin{eqnarray}\label{eqnarray:twotwo}
\left|\operatorname{Cov}\left(u_t''(y),u_t''(z)\right)\right| &\le& C_1\exp\{-C_2r^2\}\int_0^1\frac{1}{(\sigma_t\sqrt{s})^4}\left|\frac{y-\alpha s}{\sigma_t\sqrt{s}}\right|\left|\frac{z-\alpha s}{\sigma_t\sqrt{s}}\right|\phi^2\left(\frac{\sqrt{\frac{y^2+z^2}{2}}-\alpha s}{\sigma_t\sqrt{s}}\right)ds \nonumber\\
&\le& C_1\exp\{-C_2r^2\}t^{5/2},
\end{eqnarray}
writing $\displaystyle{\frac{y-\alpha s}{\sigma_t\sqrt{s}}}$ as $\displaystyle{\frac{y-\eta}{\sigma_t\sqrt{s}}+\frac{\eta-\alpha s}{\sigma_t\sqrt{s}}}$ and similarly for $\displaystyle{\frac{z-\alpha s}{\sigma_t\sqrt{s}}}$.
Furthermore, we see that
\begin{eqnarray}\label{eqnarray:detinv}
\det \Sigma_{22}=\det \Sigma_t(y,z)&=& \operatorname{Var}(u_t'(y))\operatorname{Var}(u_t'(z))-\operatorname{Cov}^2\left(u_t'(y),u_t'(z)\right)\nonumber\\
&\ge & \operatorname{Var}(u_t'(y))\operatorname{Var}(u_t'(z))-\exp\{-Cr^2\}\operatorname{Var}^2(u_t'(\eta))\nonumber\\
&\ge & \frac{1}{2}\operatorname{Var}(u_t'(y))\operatorname{Var}(u_t'(z))
\end{eqnarray}
for sufficiently large $r$. Using equations (\ref{eqnarray:oneone}), $\dots$,(\ref{eqnarray:detinv}) to estimate the right side of equation (\ref{eqnarray:cov}), we see that there is a $K>0$ for which $$\displaystyle{\operatorname{Cov}\left(|u_t''(y)|,|u_t''(z)| \ \big| u_t'(y)=0, u_t'(z)=0\right) \le C_1\exp\{-C_2r^2\}t^K }.$$ Plugging this into the expression (\ref{eqnarray:cross}), we get
\begin{eqnarray}\label{eqnarray:crossrev}
\mathbb{E}\left(N_t(I_1)N_t(I_2)\right)&\le & C\int_{I_1\times I_2}\frac{C_1\exp\{-C_2r^2\}t^K +\frac{2}{\pi}\sqrt{\operatorname{Var}(u_t''(y))}\sqrt{\operatorname{Var}(u_t''(z))}}{2\pi\sqrt{\operatorname{Var}(u_t'(y))\operatorname{Var}(u_t'(z))}}dydz.
\end{eqnarray}
We know from Lemma \ref{lem:seconddervar} that $$\operatorname{Var}(u_t''(y)) \le Ct^{3/2}.$$ Thus $$\frac{2}{\pi}\sqrt{\operatorname{Var}(u_t''(y))}\sqrt{\operatorname{Var}(u_t''(z))}=O(t^{3/2}).$$ If we choose $r=\sqrt{M\log t}$ for a large enough $M$, then $$\displaystyle{C_1\exp\{-C_2r^2\}t^K <<\frac{2}{\pi}\sqrt{\operatorname{Var}(u_t''(y))}\sqrt{\operatorname{Var}(u_t''(z))}}.$$ Consequently, from (\ref{eqnarray:crossrev}),
\begin{eqnarray*}
\mathbb{E}\left(N_t(I_1)N_t(I_2)\right)&\le & C\int_{I_1\times I_2}\frac{\sqrt{\operatorname{Var}(u_t''(y))}\sqrt{\operatorname{Var}(u_t''(z))}}{\pi^2\sqrt{\operatorname{Var}(u_t'(y))\operatorname{Var}(u_t'(z))}}dydz\\\\
&=&C\left(\int_{I_1}\frac{\sqrt{\operatorname{Var}(u_t''(y))}}{\pi\sqrt{\operatorname{Var}(u_t'(y))}}dy\right)\left(\int_{I_2}\frac{\sqrt{\operatorname{Var}(u_t''(z))}}{\pi\sqrt{\operatorname{Var}(u_t'(z))}}dz\right)\\\\
&=&\left(1+O\left(\sqrt{\frac{\log t}{t}}\right)\right)\mathbb{E}\big(N_t(I_1)\big)\mathbb{E}\big(N_t(I_2)\big),
\end{eqnarray*}
where the last step above follows from (\ref{eqnarray:expmeta}) using Corollary \ref{cor:correlation} (see the proof of Lemma \ref{lem:firstint}).

Thus,
\begin{eqnarray}
\operatorname{Cov}\left(N_t(I_1),N_t(I_2)\right)=O\left(\sqrt{\frac{\log t}{t}}\right)
\end{eqnarray}
for this choice of $r$.\\\\
Now, we have all we need to compute the variance of $N_t\left(\bigcup_{j=1}^{C[\sqrt{t}/r]}I_j\right)$.
\begin{eqnarray*}
\operatorname{Var}N_t\left(\bigcup_{j=1}^{C[\sqrt{t}/r]}I_j\right)&=&\sum_{j=1}^{C[\sqrt{t}/r]}\operatorname{Var}N_t(I_j)+2\sum_{i<j}\operatorname{Cov}\left(N_t(I_i),N_t(I_j)\right)\\
&\le & C_1\frac{\sqrt{t}}{r} +C_2(\sqrt{t\log t})r^{-2},
\end{eqnarray*}
where we used Lemma \ref{lem:secondintensity} crucially in putting the constant bound on $\operatorname{Var}N_t(I_j)$.\\\\
With our choice of $r=\sqrt{M \log t}$, the above becomes
\begin{eqnarray}
\operatorname{Var}N_t\left(\bigcup_{j=1}^{C[\sqrt{t}/r]}I_j\right) \le C\sqrt{\frac{t}{\log t}}.
\end{eqnarray}
Finally, we have, for small $\epsilon >0$,
\begin{eqnarray*}
P\left(\left|\frac{N_t([\delta, \alpha])}{\mathbb{E}N_t([\delta, \alpha])}-1\right| \ge \epsilon\right)&=&P\left(\left| N_t([\delta, \alpha])-\mathbb{E}N_t([\delta, \alpha])\right| \ge \epsilon \mathbb{E}N_t([\delta, \alpha])\right)\\
& \le & \sum_{l=1}^{[\sqrt{M\log t}]}P\left(\left| N_t\left(\bigcup_{j\in \mathcal{C}_l}I_j\right)-\mathbb{E}N_t\left(\bigcup_{j\in \mathcal{C}_l}I_j\right)\right| \ge \epsilon \mathbb{E}N_t\left(\bigcup_{j\in \mathcal{C}_l}I_j\right)\right)\\
&\le & \sum_{l=1}^{[\sqrt{M\log t}]}\frac{\operatorname{Var}N_t\left(\bigcup_{j\in \mathcal{C}_l}I_j\right)}{\epsilon^2 \left(\mathbb{E}N_t\left(\bigcup_{j\in \mathcal{C}_l}I_j\right)\right)^2}\\
& \le &\frac{C}{\epsilon^2}\sqrt{\log t} \frac{\sqrt{\frac{t}{\log t}}}{\left(\sqrt{\frac{t}{\log t}}\right)^2} = \frac{C}{\epsilon^2}\frac{\log t}{\sqrt{t}},
\end{eqnarray*}
which goes to zero as $t \rightarrow \infty$.\qed
\section{\textbf{The continuous two dimensional Conga line}}\label{section:twodimensional}
We will write the (unscaled) continuous two dimensional Conga line as $u(\cdot,t)=(u_1(\cdot,t),u_2(\cdot,t))$, where $u_1$ and $u_2$ are independent and each has the same distribution as the (unscaled) continuous one dimensional Conga line investigated in the previous section. First, we give a fluctuation estimate for $u$ and its first derivative on intervals of the form $[k,k+1]$ for $k \in [\delta t,\alpha t]$ ($\delta \in (0,\alpha)$) for sufficiently large $t$.
\begin{lem}\label{lem:fluctuation}
Take any $\delta \in (0,\alpha)$, $\epsilon>0$ and integer $p \ge 1$. For $t>0$ and $k \in [\delta t,\alpha t]$, there is a constant $C_p$ depending only on $p$ such that
\begin{eqnarray}\label{eqnarray:mombound}
\mathbb{E}\left(\sup_{x \in [k,k+1]}|u(x,t)-u(k,t)|^{2p}\right) &\le& C_pk^{-p/2}\nonumber\\
\mathbb{E}\left(\sup_{x \in [k,k+1]}|\partial_xu(x,t)-\partial_xu(k,t)|^{2p}\right) &\le& C_pk^{-3p/2}
\end{eqnarray}
Consequently, for any $\eta>0$ and $\mu>0$, with probability one, there is a positive integer $N$ such that for all $n \ge N$,
 \begin{eqnarray}\label{eqnarray:supas}
\sup_{x \in [k,k+1]}|u(x,n)-u(k,n)| &\le& \mu k^{-(\frac{1}{4}-\eta)}\nonumber\\
\sup_{x \in [k,k+1]}|\partial_xu(x,n)-\partial_xu(k,n)| &\le& \mu k^{-(\frac{3}{4}-\eta)}
\end{eqnarray}
for all $k \in [\delta n, \alpha n] \cap \Z$.
\end{lem}
\textbf{Proof: }From the scaling relation (\ref{equation:exeq}) and Corollary \ref{cor:ana}, we see that for sufficiently large $t$ and $k \in [\delta t,\alpha t]$, $u(\cdot,t)$ has a power series expansion around $k$ that converges in $[k,k+1]$. Thus, we can write $$\sup_{x \in [k,k+1]}|u(x,t)-u(k,t)| \le \sum_{l=1}^{\infty}\frac{|\partial_x^{l}u(k,t)|}{l!}.$$ The p-th moment can be bounded as
\begin{equation}\label{equation:pmom}
\mathbb{E}\left(\sup_{x \in [k,k+1]}|u(x,t)-u(k,t)|^{2p}\right) \le \left(\sum_{l=1}^{\infty}\frac{[\mathbb{E}|\partial_x^{l}u(k,t)|^{2p}]^{1/2p}}{l!}\right)^{2p}.
\end{equation}
Now we collect the moment estimates of the derivatives from the previous section. In the following, $C_1, C_2, C_1', C_2',\dots$ will be constants depending only on $\alpha$. From Lemmas \ref{lem:derivativecov} and \ref{lem:seconddervar} along with the scaling relation (\ref{equation:exeq}), we have $$\mathbb{E}|\partial_x^{l}u(k,t)|^2 \sim C_lk^{-(l-\frac{1}{2})}$$ for $l=1,2$ and $k \in [\delta t, \alpha t]$. To bound the variances of the higher order derivatives, we use (\ref{eqnarray:higherdervar}) (with $\epsilon$ there taken to be $\frac{y}{2\alpha}$) and (\ref{equation:exeq}) to get $$\mathbb{E}|\partial_x^{l}u(k,t)|^2 \sim C_l\frac{(\sqrt{2\alpha})^{2l-1}}{\sigma^{2l-1}k^{(l-\frac{1}{2})}}(l-1)! + C_l'\frac{8^{l-1}\alpha^{2l-2}}{k^{2(l-1)}}[(l-1)!]^2$$ for $l \ge 3$ and $k \in [\delta t, \alpha t]$. Along with these estimates, and the fact that the derivatives $\partial_x^lu(k,n)$ are Gaussian, we use (\ref{equation:pmom}) to get the first bound in (\ref{eqnarray:mombound}). The second bound follows similarly. Finally, (\ref{eqnarray:supas}) follows from (\ref{eqnarray:mombound}) and the Borel-Cantelli Lemma.\qed

For any integer $n \ge 1$, let $\{X(x,n): x \in [0,n]\}$ denote the \textit{linear interpolation} of $\{X_k(n): k \in [0,n] \cap \Z\}$. Note that $X(\cdot,n)$ is differentiable everywhere except the integers. We set the convention $\partial_xX(k,n)=X_{k+1}(n)-X_k(n)$. Then, Theorem \ref{thm:conga} and Lemma \ref{lem:fluctuation} together imply the following theorem.
\begin{thm}\label{thm:congafl}
For any $\eta>0$ and $\mu>0$, with probability one, there is a positive integer $N$ such that for all $n \ge N$,
 \begin{eqnarray}\label{eqnarray:coneclose}
\sup_{x \in [\delta n, \alpha n]}x^{\frac{1}{4}-\eta}|X(x,n)-u(x,n)| &\le& \mu\nonumber\\
\sup_{x \in [\delta n, \alpha n]}x^{\frac{3}{4}-\eta}|\partial_xX(x,n)-\partial_xu(x,n)| &\le& \mu
\end{eqnarray}
\end{thm}
The above theorem tells us that for large $n$, $\{(X(x,n), \partial_xX(x,n)): x \in [\delta n, \alpha n]\}$ comes uniformly close to $\{(u(x,n), \partial_xu(x,n)): x \in [\delta n, \alpha n]\}$. Furthermore, the distance shrinks as $x$ gets larger. This explains why the discrete Conga line sufficiently far away from the tip \textit{looks so smooth} in Figure 2.

In the rest of this section, we study properties of the continuous \textit{scaled} two dimensional Conga line $u_t=(u_{1,t},u_{2,t})$. 
\subsection{Analyzing length}\label{subsection:lengthanalyze}
The length of the Conga line in the interval $[\delta,\alpha]$ is given by $\displaystyle{l_t=\int_{\delta}^{\alpha}|u_t'(x)|dx}$, where $|\cdot|$ denotes the $L^2$ norm. In this section, we give estimates for the expected length and its concentration about the mean.
\begin{lem}
$\mathbb{E}(l_t) \sim t^{1/4}.$
\end{lem}
\textbf{Proof: }From Lemma \ref{lem:derivativecov}, we have
\begin{eqnarray}\label{eqnarray:sigmalen}
\operatorname{Var}(u_{1,t}'(x))=\frac{\sqrt{t}}{2\sqrt{\pi \alpha}\sigma x^{1/2}}\left(1 + O^{\infty}\left(\sqrt{\frac{\log t}{t}}\right)\right).
\end{eqnarray}
Thus,
\begin{eqnarray*}
\mathbb{E}|u_t'(x)|&=&\mathbb{E}\sqrt{u_{1,t}'^2(x)+u_{2,t}'^2(x)}=\sqrt{\frac{\pi}{2}}\sqrt{\operatorname{Var}(u_{1,t}'(x))}\\
&=&\frac{\pi^{1/4}t^{1/4}}{2\sigma^{1/2}\alpha^{1/4}x^{1/4}}\left(1 + O^{\infty}\left(\sqrt{\frac{\log t}{t}}\right)\right),
\end{eqnarray*}
where the second equality above follows from the fact that $u_t=(u_{1,t},u_{2,t})$ has a bivariate normal distribution. Using this, we get
\begin{eqnarray*}
\mathbb{E}(l_t)=\int_{\delta}^{\alpha}\mathbb{E}|u_t'(x)|dx &=&\int_{\delta}^{\alpha}\left(\frac{\pi^{1/4}t^{1/4}}{2\sigma^{1/2}\alpha^{1/4}x^{1/4}}\right)dx + O\left(\frac{\sqrt{\log t}}{t^{1/4}}\right)\\\\
&=& \frac{2\pi^{1/4}t^{1/4}}{3\sigma^{1/2}\alpha^{1/4}}\left(\alpha^{3/4}-\delta^{3/4}\right) + O\left(\frac{\sqrt{\log t}}{t^{1/4}}\right).
\end{eqnarray*}
\qed\\\\
But this gives us only a rough estimate of the behaviour of length for large $t$. To get a better idea of how the length behaves for large time $t$, we need higher moments and, if possible, some form of concentration about the mean. 
The next lemma gives us an estimate of the variance of $\displaystyle{l_t}$.
\begin{lem}\label{lem:lengthvar}
\begin{eqnarray*}
\operatorname{Var}(l_t)=O(1).
\end{eqnarray*}
\end{lem}
\textbf{Proof: }From Lemma \ref{lem:derivativecov}, we know that for $\delta \le x,y \le \alpha$, $\rho_t(x,y)=\operatorname{Corr}(u_{1,t}'(x), u_{1,t}'(y))$ satisfies
\begin{eqnarray}\label{eqnarray:rhovar}
0 < \rho_t(x,y) \le \exp\left\lbrace-C_2t(x-y)^2\right\rbrace\left(\frac{(xy)^{1/4}}{\left(\frac{x^2+y^2}{2}\right)^{1/4}}\right)\left(1 + O^{\infty}\left(\sqrt{\frac{\log t}{t}}\right)\right)
\end{eqnarray}
and $\sigma_t^2(x)=\operatorname{Var}(u_{1,t}(x))$ satisfies (\ref{eqnarray:sigmalen}).

Let $f$ denote the probability density of $(u_{1,t}'(x),u_{2,t}'(x),u_{1,t}'(y),u_{2,t}'(y))$. We use the fact that this is a Gaussian random vector to write down $f$ explicitly as
\begin{eqnarray*}
f(\mathbf{a},\mathbf{b})&=&\frac{1}{(2\pi)^2\sigma_t^2(x)\sigma_t^2(y)(1-\rho_t^2(x,y))}\\
&\quad & \times \exp\left\lbrace-\frac{1}{2(1-\rho_t^2(x,y))}\left(\frac{|\mathbf{a}|^2}{\sigma_t^2(x)}+\frac{|\mathbf{b}|^2}{\sigma_t^2(y)}-2\rho_t(x,y)\left\langle\frac{\mathbf{a}}{\sigma_t(x)},\frac{\mathbf{b}}{\sigma_t(y)}\right\rangle\right)\right\rbrace
\end{eqnarray*}
where $\mathbf{a}=(a_1,a_2)$, $\mathbf{b}=(b_1,b_2)$, $|\cdot|$ represents the $L^2$ norm and $\langle \cdot,\cdot \rangle$ represents the dot product of two vectors.
Using this expression, we can write down
\begin{eqnarray*}
\operatorname{Cov}(|u_t'(x)|, |u_t'(y)|)=\sigma_t(x)\sigma_t(y)\left(g(\rho_t(x,y))-g(0)\right)
\end{eqnarray*}
where $$g(\rho)=\frac{1}{(2\pi)^2(1-\rho^2)}\int_{\mathbb{R}^2 \times \mathbb{R}^2}|\mathbf{a}||\mathbf{b}|\exp\left\lbrace-\frac{1}{2(1-\rho^2)}\left(|\mathbf{a}|^2+|\mathbf{b}|^2-2\rho\langle\mathbf{a},\mathbf{b}\rangle\right)\right\rbrace \operatorname{d}\mathbf{a} \ \operatorname{d}\mathbf{b}.$$
It is routine to verify from the form of $g$ above that $$\sup_{0 < \rho \le 1/2}\frac{|g(\rho)-g(0)|}{\rho}<\infty.$$
Let $S$ denote the above supremum. From (\ref{eqnarray:rhovar}), it follows that there is $0<A<\infty$ such that if $|x-y|\ge A/ \sqrt{t}$, then $\rho_t(x,y) \le 1/2$ for sufficiently large $t$. On $\{(x,y) \in [\delta, \alpha]^2: |x-y| \ge A/ \sqrt{t}\}$, we bound $\operatorname{Cov}(|u_t'(x)|, |u_t'(y)|)$ by $S \rho_t(x,y)\sigma_t(x)\sigma_t(y)$. On $\{(x,y) \in [\delta, \alpha]^2: |x-y| < A/ \sqrt{t}\}$, we bound $\operatorname{Cov}(|u_t'(x)|, |u_t'(y)|)$ by $\sigma_t(x)\sigma_t(y)$ which follows from Cauchy-Schwarz inequality. From (\ref{eqnarray:sigmalen}), we see that $\sigma_t(x)\sigma_t(y)$ is bounded above by $C\sqrt{t}$ for some finite constant $C$.\\\\
Now, we can combine the estimates above to bound the variance as follows:
\begin{eqnarray*}
\operatorname{Var}(l_t)&=&\int_{(x,y) \in [\delta, \alpha]^2}\operatorname{Cov}(|u_t'(x)|, |u_t'(y)|)dxdy\\
&\le &\int_{\{(x,y) \in [\delta, \alpha]^2: |x-y| \ge A/ \sqrt{t}\}}S \ C \exp\left\lbrace-C_2t(x-y)^2\right\rbrace\left(\frac{\sqrt{t}}{\left(\frac{x^2+y^2}{2}\right)^{1/4}}\right)dxdy\\
&\quad& \quad \quad + \int_{\{(x,y) \in [\delta, \alpha]^2: |x-y| < A/ \sqrt{t}\}}C\sqrt{t} \ dxdy\\
&=& O(1).
\end{eqnarray*}
\qed\\\\
Thus, although the expected length grows like $t^{1/4}$, the variance is bounded. This already tells us that the actual length cannot deviate much from the expected length.

In what we do next, we get \textit{Gaussian concentration} of length about the mean in a window of scale $O(\sqrt{\log t})$.

We know that for most useful concentration results, we need some `independence' in our model. Our strategy here is to construct a new process $\hat{u}_t$ which is very `close' to the original process $u_t$ and is nicer to analyze as $\hat{u}_t(x)$ and $\hat{u}_t(x')$ are \textit{independent} whenever $x$ and $x'$ are sufficiently far apart. As this yields a useful tool which is going to be used in later sections, we give a detailed construction.\\\\
\textbf{Construction of $\hat{u}_t$}\\\\
By Lemma \ref{lem:derivativecov}, we see that the correlation between $u_t'(x)$ and $u_t'(y)$ in the Conga line with $\displaystyle{|x-y|=\frac{\lambda}{\sqrt{t}}}$ decays like $e^{-C\lambda^2}$ as $\lambda$ increases. This indicates that pieces of the Conga line sufficiently far away are `almost independent'. In what follows, we make use of this fact.

Choose and fix a constant $M>2$. Divide the interval $\left[\delta-\sqrt{\frac{M\log t}{t}},\alpha+\sqrt{\frac{M\log t}{t}}\right]$ into subintervals $$I_k=[y_k,y_{k+1}], \ -1 \le k \le \left[(\alpha-\delta)\frac{\sqrt{t}}{\sqrt{M \log t}}\right]+1$$ of length at most $\displaystyle{\frac{\sqrt{M \log t}}{\sqrt{t}}}$. Define the process $\hat{u}_t$ as
\begin{eqnarray}\label{eqnarray:uhat}
\hat{u}_t(x)=\int_{y_{k}/\alpha}^{y_{k+3}/\alpha}\frac{x+\alpha s}{2\sigma_ts^{3/2}}\phi\left(\frac{x-\alpha s}{\sigma_t \sqrt{s}}\right)\left[W(1-s)-W(1-y_{k}/\alpha)\right]ds
\end{eqnarray}
if $\displaystyle{x \in [y_{k+1},y_{k+2}]}$.

In the following, all constants $C, C_1, C_2, \dots$ depend only on $\delta$ and $\alpha$.
\begin{lem}\label{lem:hatuprop}
$\displaystyle{\hat{u}_t}$ satisfies the following properties:\\\\
(i) $\hat{u}_t$ is smooth everywhere except possibly at the points $y_k$.\\\\
(ii) $\{\hat{u}_t(x): x \in I_k\}$ is independent of $\{\hat{u}_t(x): x \in I_{k+3}\}$ for all $k$.\\\\
(iii) For $x \in I_{k+1}$,
\begin{eqnarray*}
\left|\hat{u_t}(x)-u_t(x)+W\left(1-\frac{y_k}{\alpha}\right)\right| \le Ct^{-M/2}||W||
\end{eqnarray*} 
and
\begin{eqnarray*}
\left|u_t'(x)-\hat{u}_t'(x)\right| \le \frac{C}{\sqrt{M}}\left(\frac{x}{t}\right)^{\frac{M}{4}-\frac{1}{2}}||W||.
\end{eqnarray*}
where $\displaystyle{||W||=\sup_{0 \le s \le 1}|W_s|}$.
\end{lem}
\textbf{Proof: }Properties (i) and (ii) follow from the definition of $\hat{u}_t$.

To prove property (iii) notice that, for $\displaystyle{x \in I_{k+1}}$,
\begin{eqnarray*}
u_t(x)-\hat{u}_t(x)&=& \int_{[y_{k}/\alpha,y_{k+3}/\alpha]^c \cap [0,1]}\frac{x+\alpha s}{2\sigma_ts^{3/2}}\phi\left(\frac{x-\alpha s}{\sigma_t \sqrt{s}}\right)W(1-s)ds\\
&\quad & + \ W(1-y_k/\alpha) - W(1-y_k/\alpha)\int_{[y_{k}/\alpha,y_{k+3}/\alpha]^c \cap [0,\infty]}\frac{x+\alpha s}{2\sigma_ts^{3/2}}\phi\left(\frac{x-\alpha s}{\sigma_t \sqrt{s}}\right)ds.
\end{eqnarray*}
From (\ref{eqnarray:derivativeK}) and a similar equation for the derivatives of $\hat{u}_t$, we obtain
\begin{eqnarray*}
u_t'(x)-\hat{u}_t'(x)&=& \int_{[y_{k}/\alpha,y_{k+3}/\alpha]^c \cap [0,1]}\partial_sK_t^0(x,s)[W(1-s)-W(1-y_k/\alpha)]ds\\\\
&\quad & - W(1-y_k/\alpha)\int_1^{\infty}\partial_sK_t^0(x,s)ds,
\end{eqnarray*}
where $K_t^n$ is defined as in (\ref{equation:K}).
If $x \in I_{k+1}$, then $\displaystyle{\frac{x}{\alpha}-L_t^x(M) \ge \frac{y_{k}}{\alpha}}$ and $\displaystyle{\frac{x}{\alpha}+L_t^x(M) \le \frac{y_{k+3}}{\alpha}}$. So, the above differences yield part (iii).\qed

Consequently, if $\hat{l}_t$ is the length of the curve $\hat{u}_t$ restricted to $[\delta,\alpha]$, then
\begin{eqnarray}\label{eqnarray:lengthcompare}
P(|l_t-\hat{l}_t| \ge r) \le C_1\exp\left\lbrace -C_2r^2t^{\frac{M}{2}-1}\right\rbrace
\end{eqnarray}
and
\begin{eqnarray}\label{eqnarray:lengthexp}
E(\hat{l}_t)= \frac{2\pi^{1/4}t^{1/4}}{3\sigma^{1/2}\alpha^{1/4}}\left(\alpha^{3/4}-\delta^{3/4}\right) + O\left(\frac{\sqrt{\log t}}{t^{1/4}}\right)+ O\left(\frac{1}{t^{\frac{M}{4}-\frac{1}{2}}}\right).
\end{eqnarray}
So, to find the concentration of the length around the mean at time $t$, we look at the length of the curve $\hat{u}_t$.
Let $W_k$ be the Brownian motion defined on $\displaystyle{I=\left[0,\frac{3}{\alpha}\sqrt{\frac{M\log t}{t}}\right]}$ by $$\displaystyle{W_k(s)=W\left(1-\frac{y_k}{\alpha}-s\right)-W\left(1-\frac{y_k}{\alpha}\right)}.$$ For each $k$, the Brownian motions $W_k$ and $W_{k+3}$ so defined are clearly independent.

As length is an additive functional, we can find the length on subintervals $I_k$ and add them together. Heuristically, we can see that this gives us concentration as the length of the curve on every third interval is independent of each other, and as these are summed up, the errors get averaged out.\\\\
Now, we give the rigorous arguments. In the following, we fix the probability space $\displaystyle{(\Omega,\mathcal{B}(\Omega), \mathcal{P})}$, where $\displaystyle{\Omega=C\left[0,\frac{3}{\alpha}\sqrt{\frac{M\log t}{t}}\right]}$ denote the set of continuous complex valued functions on $I$ equipped with the sup-norm metric $d$, and $\displaystyle{\mathcal{P}}$ is the Wiener measure.

We need some concepts from \textit{Concentration of Measure} Theory. See \cite{ledoux} for an excellent survey of techniques in this area. We give a very brief outline of the concepts we need.\\\\
\textbf{Transportation Cost Inequalities and Concentration: }Let $(\chi, d)$ be a complete separable metric space equipped with the Borel sigma algebra $\mathcal{B}(\chi)$. Consider the p-th \textit{Wasserstein distance} between two probability measures $P$ and $Q$ on this space, defined as $$\mathcal{W}_p(P,Q)=\inf_{\pi}[\mathbb{E}d(X,X')^p]^{1/p},$$ where the infimum is over all couplings $\pi$ of a pair of random elements $(X,X')$ with the marginal of $X$ being $P$ and that of $X'$ being $Q$.

Now, fix a probability measure $P$. Suppose there is a constant $C>0$ such that for all probability measures $Q << P$, we have $$\mathcal{W}_p(P,Q) \le \sqrt{2CH(Q\mid P)},$$ where $H$ refers to the relative entropy $H(Q\mid P)=\mathbb{E}^Q \log (dQ/dP)$. Then we say that $P$ satisfies the $L^p$ \textit{Transportation Cost Inequality}. In short, we write $P \in \mathcal{T}_p(C)$.

Now, we present one of the key results which connects Transportation Cost Inequalities and Concentration of Measures.
\begin{lem}\label{lem:lipschitz}
Suppose $P$ is a probability measure on $\left(\chi, \mathcal{B} (\chi)\right)$. Suppose further that each $P$ is in $\mathcal{T}_1(C)$. Then, for any $1$-Lipschitz map $F:\chi \rightarrow \mathbb{R}$ and any $r >0$,
\begin{eqnarray*}
P\left(|F-\int FdP|>r\right) \le \exp\left\lbrace-\frac{r^2}{2C}\right\rbrace.
\end{eqnarray*}
\end{lem}
It is easy to see that $\mathcal{T}_2(C)$ implies $\mathcal{T}_1(C)$. But the main advantage in dealing with $\mathcal{T}_2(C)$ comes from its \textit{tensorization property} described in the following lemma.
\begin{lem}\label{lem:tensor}
Suppose $P_i, i=1,2,\ldots n$ are probability measures on $\left(\chi, \mathcal{B} (\chi)\right)$. Suppose further that each $P_i$ is in $\mathcal{T}_2(C)$. On $\chi^n$, define the distance between $x^n=(x_1,x_2,\dots)$ and $y^n=(y_1,y_2,\dots)$ by $$\displaystyle{d^n(x^n,y^n)=\sqrt{\sum_{i=1}^nd^2(x_i,y_i)}}.$$ Then $\displaystyle{\bigotimes_{i=1}^nP_i \in \mathcal{T}_2(C)}$ on $(\chi^n, d^n)$.
\end{lem}
The following lemma, which follows from the developments in \cite{pal}, is of key importance to us.
\begin{lem}\label{lem:transport}
The Wiener measure on $\displaystyle{C\left[0,T\right]}$ satisfies the transportation inequality $\displaystyle{\mathcal{T}_2(T)}$ with respect to the sup-norm metric.
\end{lem}
These tools are all we need to establish a concentration result for $l_t$.

Let us define the function $\displaystyle{T_t^k:\Omega \rightarrow \mathbb{R}}$ as follows:
\begin{eqnarray*}
T_t^k(f)=\int_{y_{k+1}}^{y_{k+2}}\left|\int_{0}^{\frac{y_{k+3}-y_{k}}{\alpha}}\partial_sK_t^0\left(x,s+\frac{y_k}{\alpha}\right)f(s)ds\right|dx.
\end{eqnarray*}
Notice that $\displaystyle{\hat{l}_t=\sum_{k=0}^{\frac{\sqrt{t}}{\sqrt{M \log t}}}T_t^k(W_k)}$. Suppose we prove that $\displaystyle{T_t^k}$ is Lipschitz with respect to $d$ with Lipschitz constant $C_t$. Then, with $\displaystyle{N=\left[\frac{\alpha-\delta}{3}\sqrt{\frac{t}{M\log t}}\right]}$, the functions $\displaystyle{\{T_t^{(i)}:\Omega^N \rightarrow \mathbb{R}:i=0,1,2\}}$ defined by
\begin{eqnarray*}
T_t^{(i)}(f)=\frac{1}{\sqrt{N}}\sum_{k=0}^{N-1}T_t^{3k+i}(f_{k+1})
\end{eqnarray*}
where $f=(f_1,..,f_N)$, is also Lipschitz with respect to $d^N$ with the same constant $C_t$.
\begin{lem}
For each $k$, $\displaystyle{T_t^k}$ is Lipschitz on $(\Omega,d)$ with Lipschitz constant $\displaystyle{C\sqrt{M\log t}}$ where $C$ is a constant depending only on $\delta$ and $\alpha$.
\end{lem}
\textbf{Proof: }For $f_1,f_2$ in $\Omega$, 
\begin{eqnarray*}
|T_t^k(f_1)-T_t^k(f_2)|&\le&\int_{y_{k+1}}^{y_{k+2}}\int_{0}^{\frac{y_{k+3}-y_{k}}{\alpha}}\left|\partial_sK_t^0\left(x,s+\frac{y_k}{\alpha}\right)\right||f_1(s)-f_2(s)|dsdx\\\\
&\le& ||f_1-f_2||\int_{y_{k+1}}^{y_{k+2}}\int_{0}^{\frac{y_{k+3}-y_{k}}{\alpha}}\left|\partial_sK_t^0\left(x,s+\frac{y_k}{\alpha}\right)\right|dsdx.
\end{eqnarray*}
By the estimates obtained in the proof of Lemma \ref{lem:derivative},
\begin{eqnarray*}
\int_{y_{k+1}}^{y_{k+2}}\int_{0}^{\frac{y_{k+3}-y_{k}}{\alpha}}\left|\partial_sK_t^0\left(x,s+\frac{y_k}{\alpha}\right)\right|dsdx \le C\sqrt{t}\int_{y_{k+1}}^{y_{k+2}}dx\le C\sqrt{M\log t}.
\end{eqnarray*}
This proves the lemma.\qed\\\\
Now, for any $f \in \Omega^{3N}$, define the following functions in $\Omega^N$: $f^{(1)}=(f_1,f_4...,f_{3N-2})$, $f^{(2)}=(f_2,f_5...,f_{3N-1})$ and $f^{(3)}=(f_3,f_6...,f_{3N})$. Notice that $$\hat{l}_t=\sqrt{N}\left(T_t^{(1)}(\tilde{W}^{(1)})+T_t^{(2)}(\tilde{W}^{(2)})+T_t^{(3)}(\tilde{W}^{(3)})\right)$$ where $\tilde{W}=(W_{0},W_{1},...,W_{3N-1}) \in \Omega^{3N}$. Using this fact and Lemmas \ref{lem:tensor} and \ref{lem:lipschitz}, we get for any $r>0$,
\begin{eqnarray}\label{eqnarray:cl1}
P(|\hat{l}_t-E\hat{l}_t| \ge r\sqrt{M\log t}) \le C_1 \exp\left\lbrace-C_2r^2\right\rbrace
\end{eqnarray}
As $M>2$ is arbitrary, it can be absorbed in the constant $C_2$. Thus (\ref{eqnarray:cl1}), along with (\ref{eqnarray:lengthcompare}) and (\ref{eqnarray:lengthexp}), gives us our main conclusion:
\begin{thm}\label{thm:lengthconc}
We have
\begin{eqnarray*}
P(|l_t-El_t| \ge r\sqrt{\log t}) \le C_1 \exp\left\lbrace-C_2r^2\right\rbrace,
\end{eqnarray*}
where $C_1, C_2$ are constants depending only on $\delta$ and $\alpha$.
\end{thm}
\subsection{How close is the scaled Conga Line to Brownian motion?}
Though the unscaled Conga line seen far away from the tip `smoothes out' Brownian motion more and more with increasing $t$, we see that in the simulations of the scaled Conga line, \textit{making $t$ larger actually makes the curve rougher} and resemble Brownian motion more and more. Closer analysis reveals that this in fact results from the \textit{scaling}. Again, before we supply the rigorous arguments, we give a heuristic reasoning. Looking at equation (\ref{equation:exeq}), we see that although the scaling takes the Brownian motion $W$ on $[0,t]$ to a Brownian motion $W^{(t)}$ on $[0,1]$, the width of the window on which the smoothing takes place in the unscaled Conga line, which is comparable to $\sqrt{t}$, is taken to $O(t^{-1/2})$ in the scaled version, which shrinks with time $t$.

In the following, we consider the family of two dimensional random curves $u_t(\cdot)$ indexed by $t$, and $L_t^x= \alpha^{-1}\sqrt{-M\sigma_t^2x\log \sigma_t^2x}$.
\begin{thm}\label{thm:Brownianclose}
There exists a deterministic constant $\kappa$ such that, almost surely, there is $T=T(\omega)>0$ for which
\begin{eqnarray}\label{eqnarray:Browclose}
\left| u_t(x)-W\left(1-\frac{x}{\alpha}\right)\right| \le \kappa\sqrt{-\left(\frac{x}{t}\right)^{1/2}\log\left(\frac{x}{t}\right)}
\end{eqnarray}
for all $ x \in (0,\alpha]$ satisfying $x > \alpha L_t^x$ for all $t\ge T$. In particular, for any fixed $\beta \in (0,1)$, the above holds almost surely for $x \in [t^{-\beta},\alpha]$ for all $t \ge T$. Furthermore, $\kappa$ could be chosen appropriately so that the following holds for sufficiently large $t$:
\begin{equation}\label{equation:Browcloseprob}
\P\left(\left| u_t(x)-W\left(1-\frac{x}{\alpha}\right)\right| > \kappa\sqrt{-\left(\frac{x}{t}\right)^{1/2}\log\left(\frac{x}{t}\right)} \text{ for some } x \in [t^{-\beta},\alpha]\right) \le t^{-2}.
\end{equation}
\end{thm}
Thus, the scaled Conga line is close to Brownian motion for large $t$ although the unscaled one is not, as can be seen from the right side of equation (\ref{eqnarray:Browclose}). This subsection is devoted to proving the above theorem.

For any continuous function $f:[0,1] \rightarrow \mathbb{C}$, define
\begin{equation*}
P_tf(x)=\int_0^1\frac{x+\alpha s}{2\sigma_t s^{3/2}}\phi\left(\frac{x-\alpha s}{\sigma_t\sqrt{s}}\right)f(1-s)ds.
\end{equation*} 
Note that the Conga line is given by $\displaystyle{u_t(x)=P_tW(x)}$. $P_tf$ can be thought of as a \textit{smoothing kernel} acting on the function $x \mapsto f(1-x/\alpha)$. The following lemma shows that if $f$ is Lipschitz, then for large $t$, $P_tf(x)$ is close to $f(1-x/\alpha)$.
\begin{lem}\label{lem:Lipclose}
If $f$ is Lipschitz with constant $\mathcal{C}$, then for large enough $t$ and for $ x \in (0,\alpha]$ satisfying $x > \alpha L_t^x$,
\begin{equation}\label{eqn:Lip}
\displaystyle{\left|P_tf(x)-f\left(1-\frac{x}{\alpha}\right)\right| \le \mathcal{C}\sigma_t\sqrt{x}}.
\end{equation}
\end{lem}
Note that
\begin{eqnarray*}
\left|P_tf(x)-f\left(1-\frac{x}{\alpha}\right)\right| &\le & \mathcal{C}\alpha^{-1}\int_0^1\frac{x+\alpha s}{2\sigma_t s^{3/2}}\phi\left(\frac{x-\alpha s}{\sigma_t\sqrt{s}}\right)|x-\alpha s|ds\\
&=& I_t^x + J_t^x + S_t^x
\end{eqnarray*}
where
\begin{eqnarray*}
I_t^x &=&\mathcal{C}\alpha^{-1}\int_0^{(x/\alpha)- L_t^x}\frac{x+\alpha s}{2\sigma_t s^{3/2}}\phi\left(\frac{x-\alpha s}{\sigma_t\sqrt{s}}\right)|x-\alpha s|ds \le \mathcal{C}\alpha^{-1}\int_0^{(x/\alpha)- L_t^x}\frac{x+\alpha s}{2\sigma_t s^{3/2}}\phi\left(\frac{x-\alpha s}{\sigma_t\sqrt{s}}\right)ds\\
&\le & \mathcal{C}\alpha^{-1} \barphi(\sqrt{-M\log \sigma_t^2x}) \le \mathcal{C}\alpha^{-1} (\sigma_t \sqrt{x})^M
\end{eqnarray*}
and
\begin{eqnarray*}
J_t^x &=&\mathcal{C}\alpha^{-1}\int_{(x/\alpha)- L_t^x}^{\min \left((x/\alpha)+ L_t^x,1\right)}\frac{x+\alpha s}{2\sigma_t s^{3/2}}\phi\left(\frac{x-\alpha s}{\sigma_t\sqrt{s}}\right)|x-\alpha s|ds \\
&\le & \mathcal{C}\alpha^{-1}\sigma_t\sqrt{\frac{2x}{\alpha}}\int_{(x/\alpha)- L_t^x}^{\min \left((x/\alpha)+ L_t^x,1\right)}\frac{x+\alpha s}{2\sigma_t s^{3/2}}\phi\left(\frac{x-\alpha s}{\sigma_t\sqrt{s}}\right)\left|\frac{x-\alpha s}{\sigma_t \sqrt{s}} \right|ds\\
&\le & \mathcal{C}\alpha^{-1}\sigma_t\sqrt{\frac{2x}{\alpha}}\int_{-\infty}^{\infty}|s|\phi(s)ds.
\end{eqnarray*}
Similarly as $I_t^x$, $S_t^x$ is small compared to $J_t^x$.\qed

Now, Brownian motion is not Lipschitz, but it can be uniformly approximated on $[0,1]$ by piecewise linear random functions whose Lipschitz constants can be controlled using \textit{Levy's Construction of Brownian motion} which we now briefly describe following \cite{peres}. Define the $n$-th level dyadic partition $\mathcal{D}_n=\{\frac{k}{2^n}: 0 \le k \le 2^n\}$ and let $\mathcal{D}= \cup_{n=0}^{\infty}\mathcal{D}_n$. Let $\{Z_n: n \in \mathbb{N}\}$ be $i.i.d$ standard normal random variables. Define the random piecewise linear functions $F_n$ as follows.

$F_0(x)=xZ_1$ for $x \in [0,1]$. For $n \ge 1$,
\begin{eqnarray*}
F_n(x)=
\left\{
\begin{array}{ll}
2^{-\frac{n+1}{2}}Z_x & x \in \mathcal{D}_n\setminus\mathcal{D}_{n-1}\\
0 & x\in \mathcal{D}_{n-1}\\
\text{linear} & \text{in between}
\end{array} 
\right.
\end{eqnarray*}
With this, Levy's construction says that a Brownian motion $W$ can be constructed via
\begin{eqnarray*}
W(x)=\sum_{n=0}^{\infty}F_n(x).
\end{eqnarray*}
for $x \in [0,1]$.

Let $\displaystyle{W_N(x)=\sum_{n=0}^{N}F_n(x)}$. This function serves as the piecewise linear (hence Lipschitz) approximation to $W$. From Lemma \ref{lem:Lipclose}, for any $N$,
\begin{eqnarray*}
\left|P_tW(x)-W\left(1-\frac{x}{\alpha}\right)\right| \le \sum_{n=0}^N\sigma\sqrt{\frac{x}{t}}||F_n'||_{\infty} + 2\sum_{n=N}^{\infty}||F_n||_{\infty}.
\end{eqnarray*}
Fix $c > \sqrt{2 \log 2}$. Let $\displaystyle{N^*=\inf\{n:|Z_d| \le c\sqrt{n} \ \forall \ d \in \mathcal{D} \setminus \mathcal{D}_n\}}$.
\begin{eqnarray*}
\sum_{n=0}^{\infty}P\left(\text{there exists } d \in \mathcal{D}_n \text{ with } |Z_d| > c\sqrt{n}\right) \le \sum_{n=0}^{\infty}(2^n + 1) \exp\left(\frac{-c^2n}{2}\right) < \infty.
\end{eqnarray*}
So, by Borel-Cantelli Lemma, $P(N^* < \infty)=1$. \\\\
Now, for $n > N^*$, $||F_n||_{\infty} \le c\sqrt{n}2^{-n/2}$ and $\displaystyle{||F_n'||_{\infty} \le \frac{2||F_n||_{\infty}}{2^{-n}} \le 2c\sqrt{n}2^{n/2}}$. So, for $l > N^*$, we get
\begin{eqnarray}\label{eqnarray:Levy}
\left|P_tW(x)-W\left(1-\frac{x}{\alpha}\right)\right| &\le & \sum_{n=0}^{N^*}\sigma\sqrt{\frac{x}{t}}||F_n'||_{\infty} + 2\sum_{n=N^*}^l\sigma\sqrt{\frac{x}{t}}c\sqrt{n}2^{n/2}\nonumber\\
&\quad & \hspace{3cm}+ 2\sum_{n=l}^{\infty}c\sqrt{n}2^{-n/2}.
\end{eqnarray}
Now, take $t$ large enough that, for every $x \in (0,\alpha]$, the first term is less than $\displaystyle{\sqrt{-\left(\frac{x}{t}\right)^{1/2}\log\left(\frac{x}{t}\right)}}$ and $\displaystyle{\sqrt{\frac{x}{t}} \in (2^{-l},2^{-l+1}]}$ for some $l>N^*$. Plugging this $l$ into equation (\ref{eqnarray:Levy}) and using the fact that the second sum above is dominated by its last term, and the third sum is dominated by its leading term, we get (\ref{eqnarray:Browclose}).

To prove (\ref{equation:Browcloseprob}), note that the last two sums in (\ref{eqnarray:Levy}) are bounded above by $\displaystyle{C'\sqrt{-\left(\frac{x}{t}\right)^{1/2}\log\left(\frac{x}{t}\right)}}$ where $C'$ is deterministic (does not depend on $N^*$). So, it suffices to control the first sum. In the remaining part of the proof, $C$ will denote a generic, deterministic constant.

Note that, the first sum is bounded above by $\displaystyle{C\sigma \sqrt{\frac{x}{t}}2^{N^*/2}U_{N^*}}$, where $U_{n}=\sup_{d \in \mathcal{D}_n}|Z_d|$. Thus, the probability in (\ref{equation:Browcloseprob}) is bounded above by $\P(2^{N^*/2}U_{N^*} > t^{1/4}\sqrt{\log t})$. Choose and fix any $\epsilon \in (0,\frac{1}{2})$. Choose $c$ in the definition of $N^*$ above to satisfy $\displaystyle{c>\sqrt{\left(\frac{4}{\epsilon}+2\right)\log 2}}$. Now, $$\P(2^{N^*/2}U_{N^*} > t^{1/4}\sqrt{\log t}) \le \P(N^* > \epsilon \log_2 t) + \P(U_{\lfloor\epsilon \log_2 t\rfloor} > t^{\frac{1}{2}\left(\frac{1}{2}-\epsilon\right)}\sqrt{\log t}).$$ The first probability above can be bounded as $$\P(N^* > \epsilon \log_2 t) \le Ct^{-\epsilon\left(\frac{c^2}{2 \log 2}-1\right)}.$$ The second probability has the following bound: $$\P(U_{\lfloor\epsilon \log_2 t\rfloor} > t^{\frac{1}{2}\left(\frac{1}{2}-\epsilon\right)}\sqrt{\log t}) \le C(1+t^{\epsilon})\exp\left\lbrace-\frac{1}{2}t\left(\frac{1}{2}-\epsilon\right)\log t\right\rbrace.$$ These bounds, along with $\epsilon\left(\frac{c^2}{2 \log 2}-1\right)>2$, give (\ref{equation:Browcloseprob}). 
\subsection{Analyzing number of loops}
A \textbf{loop} $L$ in a continuous curve $f:\mathbb{R} \rightarrow \mathbb{C}$ is defined as a restriction of the form $f|_{[a,b]}$ where $f(a)=f(b)$ and $f$ is injective on $[a,b)$. Note that $L$ divides the plane into a bounded component and an unbounded component. Define the size of the loop 
\begin{eqnarray*}
s(L)=\sup\{R>0: \exists \ x \in \text{ the bounded component } \mathcal{B} \text{ of } L \text{ such that } B(x,R) \subseteq \mathcal{B}\}.
\end{eqnarray*}
It can be shown (the quick way is to look at the expectation meta-theorem from \cite{adlertaylor}) that if $f$ is a continuously differentiable Gaussian process, then with probability one, it has no \textit{singularities} (points where the first derivative of both Re $f$ and Im $f$ vanish). Using this fact, it is easy to see that if $I$ is a compact interval on which $f$ is not injective, then $f|_I$ has at least one loop $L$ of positive size.

As the number of loops is bounded above by the number of critical points of Re $f$ (equivalently Im $f$), we see that by Lemma \ref{lem:critprob}, for a large fixed $t$, the number of loops in the Conga line is bounded above by $C\sqrt{t}$ with very high probability. This section is dedicated to achieving a lower bound. The simulation (Figure 2) shows a number of loops, most of them being small. In the following, we obtain a lower bound for the number of \textit{small loops}, which differs from the upper bound by a logarithmic factor. For this, a key ingredient is the \textit{Support Lemma for Brownian Motion} which we state as the following lemma:
\begin{lem}\label{lem:support}
If $f:[0,1] \rightarrow \mathbb{C}$ is continuous with $f(0)=0$ and $W$ is a complex Brownian motion on $[0,1]$, then for any $\epsilon >0$,
\begin{eqnarray*}
P(||W-f||<\epsilon)>0,
\end{eqnarray*}
where $||g||=\sup_{x\in [0,1]}|g(x)|$.
\end{lem}
The above lemma can be proved either by approximating $f$ by piecewise linear functions and using Levy's construction of Brownian motion, or by an application of the \textit{Girsanov Theorem} (see \cite{revuzyor}).

We also need to exploit the exponentially decaying correlation between $u_t(x)$ and $u_t(x')$ as $|x-x'|$ increases (see Lemma \ref{lem:derivativecov}) by bringing into play the approximation of $u_t$ by the process $\hat{u}_t$ introduced in Subsection \ref{subsection:lengthanalyze}.

Now, we state the main theorem of this section.
\begin{thm}\label{thm:loopcount}
Choose $R > 6 \kappa$, where $\kappa$ is the constant in Theorem \ref{thm:Brownianclose}. Let $N_t^l$ be the number of loops of size less than or equal to $2R\left(\frac{\log t}{t}\right)^{1/4}$ in the (scaled) Conga line $u_t$ in $[\delta,\alpha]$ at time $t$. Then there exist constants $C$ and $C'$ such that
\begin{eqnarray*}
P\left(C\sqrt{\frac{t}{\log t}}\le N_t^l \le C'\sqrt{t}\right) \rightarrow 1
\end{eqnarray*}
as $t \rightarrow \infty$.
\end{thm}
\textbf{Proof: }The upper bound follows from Lemma \ref{lem:critprob}.\\\\
Proving the lower bound is more involved.\\\\
Our strategy is to choose a function $f$ which has a loop and run the Brownian motion $W$ in a narrow tube around $f$, which, by Theorem \ref{lem:support}, we can do with positive probability. Now by Theorem \ref{thm:Brownianclose}, we know that for large $t$, $u_t$ is `close' to the Brownian motion $W$ with very high probability, and thus the curve $u_t$ is forced to run in a narrow sausage around $f$ thereby inducing a loop.

Such a function is $\displaystyle{f(x)=C((4x-2)^3-(4x-2),1-(4x-2)^2)}$ for $x \in [0,1]$, where $C$ is a suitably chosen constant to make the size of the loop in $f$ to be $R$. Let us denote the continuous functions restricted to the $\epsilon$-sausage around $f\mid_{[a,b]}$ as $$S(f;\epsilon,[a,b])=\{g \in C[a,b]: ||f-g|| < \epsilon\}.$$ Fix $\alpha'$ such that $\delta<\alpha'<\alpha$. For $x\in [\delta,\alpha']$, define $$f^{(t)}_x(s)= \left(\frac{M \log t}{t}\right)^{1/4}f\left(\frac{1}{\alpha}\sqrt{\frac{t}{M \log t}}(s-x)\right), \ \ x \le s \le x+\alpha\sqrt{\frac{M \log t}{t}}$$
where $M>2$ is any fixed constant as in Subsection \ref{subsection:lengthanalyze}.
For any continuously differentiable complex-valued Gaussian process $g$ defined on a subset of $[0,1]$ containing $\left[x,x+\alpha\sqrt{\frac{M \log t}{t}}\right]$ and any complex number $c$, if $g \in S\left(c+f^{(t)}_x; \frac{R}{2}\left(\frac{M \log t}{t}\right)^{1/4},\left[x,x+\alpha\sqrt{\frac{M \log t}{t}}\right]\right)$, then $g$ has a self-intersection on $\left[x,x+\alpha\sqrt{\frac{M \log t}{t}}\right]$ and thus, due to absence of singularities with probability one, $g$ has a loop of positive size on this interval.

We break up the proof into parts:
\begin{itemize}
\item[(i)] In Lemma \ref{lem:loopone}, we prove that the probability of $u_t\mid_{\left[x,x + \alpha\sqrt{\frac{M \log t}{t}}\right]}$ having a loop of size comparable to $\left(\frac{\log t}{t}\right)^{1/4}$ is bounded below uniformly for all $x\in [\delta,\alpha']$ by a fixed positive constant $p$ independent of $x$ and $t$.
\item[(ii)] We use part (iii) of Lemma \ref{lem:hatuprop} and Lemma \ref{lem:loopone} to deduce that the probability of $\hat{u}_t$ having a loop of size comparable to $\left(\frac{\log t}{t}\right)^{1/4}$ on each interval $I_{k+1}$ is bounded below by $p/2$.
\item[(iii)] We use the independence of $\displaystyle{\hat{u}|_{I_k}}$ and $\displaystyle{\hat{u}|_{I_{k+3}}}$ for every $k$ to deduce in Lemma \ref{lem:looptwo} that the total number of such loops in $\hat{u}_t$ is bounded below by $\frac{p}{4}\left[(\alpha'-\delta)\sqrt{\frac{t}{M \log t}}\right]$ with very high probability.
\item[(iv)] We finally use part (iii) of Lemma \ref{lem:hatuprop} again to translate the result of Lemma \ref{lem:looptwo} to the original process $u_t$ in Lemma \ref{lem:loopthree}.
\end{itemize}
\begin{lem}\label{lem:loopone}
There is a constant $p>0$ independent of $x$ and $t$ such that
\begin{eqnarray*}
P\left(u_t|_{\left[x,x + \alpha\sqrt{\frac{M \log t}{t}}\right]} \in S\left(W\left(1-\frac{x}{\alpha}\right)+f^{(t)}_x; \frac{R}{2}\left(\frac{M \log t}{t}\right)^{1/4},\left[x,x+\alpha\sqrt{\frac{M \log t}{t}}\right]\right) \right)\ge p>0
\end{eqnarray*}
for all $x\in [\delta, \alpha']$, for all sufficiently large $t$.
\end{lem}
\textbf{Proof: }Choose and fix any $x\in [\delta,\alpha']$. By Theorem \ref{thm:Brownianclose} and by the translation and scaling invariance of Brownian motion, we get for $R >6\kappa$ (here $\kappa$ is the constant in Theorem \ref{thm:Brownianclose}) and large $t$,\\\\
\newcommand{\abs}[1]{\left\vert #1 \right\vert}
$\displaystyle{
P\left(u_t|_{\left[x,x + \alpha\sqrt{\frac{M \log t}{t}}\right]} \in S\left(W\left(1-\frac{x}{\alpha}\right)+f^{(t)}_x; \frac{R}{2}\left(\frac{M \log t}{t}\right)^{1/4},\left[x,x+\alpha\sqrt{\frac{M \log t}{t}}\right]\right) \right)}$\\\\
$\displaystyle{\ge P\left(\sup_{s \in \left[x,x + \alpha\sqrt{\frac{M \log t}{t}}\right]}\abs{u_t(s)-W\left(1-\frac{s}{\alpha}\right)} \le \kappa\left(\frac{M \log t}{t}\right)^{1/4}\right.}\text{ and }$\\
$\displaystyle{\left. \ \ \ \ \hspace{2cm}\sup_{s \in \left[x,x + \alpha\sqrt{\frac{M \log t}{t}}\right]}\abs{W\left(1-\frac{s}{\alpha}\right)-W\left(1-\frac{x}{\alpha}\right)-f^{(t)}_x(s)}\le \frac{R}{3}\left(\frac{M \log t}{t}\right)^{1/4}\right)}$\\\\
$\displaystyle{ \ge P\left(\sup_{s \in \left[x,x+\alpha\sqrt{\frac{M \log t}{t}}\right]}\abs{W\left(\frac{s-x}{\alpha}\right)-f^{(t)}_x(s)}\le \frac{R}{3}\left(\frac{M \log t}{t}\right)^{1/4}\right)}$\\
$\displaystyle{\hspace{4cm}-\ P\left(\sup_{s \in \left[x,x + \alpha\sqrt{\frac{M \log t}{t}}\right]}\abs{u_t(s)-W\left(1-\frac{s}{\alpha}\right)}> \kappa\left(\frac{M \log t}{t}\right)^{1/4}\right)}$\\\\
$\displaystyle{= P\left(\sup_{s \in \left[0,\alpha\sqrt{\frac{M \log t}{t}}\right]}|W^{(t)}_0(s)-f^{(t)}_0(s)|\le \frac{R}{3}\left(\frac{M \log t}{t}\right)^{1/4}\right)}$\\
$\displaystyle{\hspace{4cm}-\ P\left(\sup_{s \in \left[x,x + \alpha\sqrt{\frac{M \log t}{t}}\right]}\abs{u_t(s)-W\left(1-\frac{s}{\alpha}\right)}> \kappa\left(\frac{M \log t}{t}\right)^{1/4}\right)}$\\\\
$\displaystyle{\ge P\left(\sup_{s \in [0,1]}|W(s)-f(s)|\le \frac{R}{3}\right)}$\\
$\displaystyle{\hspace{4cm}-\ P\left(\sup_{s \in \left[\delta, \alpha\right]}\abs{u_t(s)-W\left(1-\frac{s}{\alpha}\right)}> \kappa\left(\frac{M \log t}{t}\right)^{1/4}\right)}$\\
$\displaystyle{\ge p >0}.$\\\\
Here we used Theorem \ref{thm:Brownianclose} and Theorem \ref{lem:support} for the last step. By virtue of the second last step above, we can choose $p$ independent of $x$ and $t$, and the above lower bound works uniformly for all $x \in [\delta,\alpha']$.\qed

Recall that by part (iii) of Lemma \ref{lem:hatuprop}, we know that for $x \in I_{k+1}$,
\begin{eqnarray*}
\left|\hat{u}_t(x)-u_t(x)+W\left(1-\frac{y_k}{\alpha}\right)\right| \le C(\delta)t^{-M/2}||W||.
\end{eqnarray*} 
Define the event
\begin{eqnarray*}
A_k=\left\lbrace\hat{u}_t\mid_{I_{k+1}} \in S\left(W(1-y_{k+1}/\alpha)-W(1-y_k/\alpha) +f^{(t)}_{y_{k+1}}; \frac{R}{2}\left(\frac{M \log t}{t}\right)^{1/4},I_{k+1}\right)\right\rbrace.
\end{eqnarray*}
If $A_k$ holds, then $\hat{u}_t$ has a loop in $I_{k+1}$. Write $$\displaystyle{S_t=\sum_{k=1}^{\left[(\alpha'-\delta)\sqrt{\frac{t}{M \log t}}\right]}\mathbb{I}_{A_k}}.$$ Then the following holds.
\begin{lem}\label{lem:looptwo}
\begin{eqnarray*}
P\left(S_t<\frac{p}{4}\left[(\alpha'-\delta)\sqrt{\frac{t}{M \log t}}\right]\right) &\le 3\exp\left\lbrace-\frac{p^2}{8} \left[(\alpha'-\delta)\sqrt{\frac{t}{M \log t}}\right]\right\rbrace.
\end{eqnarray*}
\end{lem}
\textbf{Proof: }By Lemma \ref{lem:loopone}, it is easy to see that
\begin{eqnarray*}
P\left(A_k \right)&\ge& P\left\lbrace u_t\vert_{I_{k+1}} \in S\left(W(1-y_{k+1}/\alpha) +f^{(t)}_{y_{k+1}}; \left(\frac{R}{2}-\epsilon\right)\left(\frac{M \log t}{t}\right)^{1/4},I_{k+1}\right)\right.\\
&&\hspace{2 in}\left.\text{ and } ||W||\le \frac{\epsilon t^{M/2}}{C(\delta)}\left(\frac{M \log t}{t}\right)^{1/4}\right\rbrace\\
&\ge & \frac{p}{2}>0
\end{eqnarray*}
for large enough $t$ and small enough $\epsilon$. Thus we see that $\displaystyle{ES_t \ge \frac{p}{2}\left[(\alpha'-\delta)\sqrt{\frac{t}{M \log t}}\right]}$.\\\\
Now, as $\hat{u_t}$ is independent on every third interval, so $A_k$ is independent of $A_{k+3}$ for every $k$. The result now follows by Bernstein's Inequality.\qed

The above implies that with very high probability $\displaystyle{S_t\ge\frac{p}{4}\left[(\alpha'-\delta)\sqrt{\frac{t}{M \log t}}\right]}$.\\\\
Define the event $$\displaystyle{B_k=\{u_t \text{ has a loop on } I_{k+1}\}}$$ and the corresponding sum $$\displaystyle{\tilde{N}_t^l=\sum_{k=1}^{\left[(\alpha'-\delta)\sqrt{\frac{t}{M \log t}}\right]}\mathbb{I}_{B_k}}.$$
Our final lemma is the following.
\begin{lem}\label{lem:loopthree}
\begin{eqnarray*}
P\left(N_t^l \ge \frac{p}{4}\left[(\alpha'-\delta)\sqrt{\frac{t}{M \log t}}\right]\right) \rightarrow 1
\end{eqnarray*}
 as $t \rightarrow \infty$.
\end{lem}
\textbf{Proof: }Note that $N_t^l \ge \tilde{N}_t^l$. By part (iii) of Lemma \ref{lem:hatuprop}, we note that for small enough $\epsilon>0$, the events $A_k$ and  $\displaystyle{\left\lbrace||W|| \le\frac{\epsilon t^{M/2}}{C(\delta)}\left(\frac{M \log t}{t}\right)^{1/4}\right\rbrace}$ imply that $B_k$ holds. We see that, for large $t$,
\begin{eqnarray*}
P\left(\tilde{N}_t^l \ge \frac{p}{4}\left[(\alpha'-\delta)\sqrt{\frac{t}{M \log t}}\right]\right) &\ge & P\left(S_t\ge\frac{p}{4}\left[(\alpha'-\delta)\sqrt{\frac{t}{M \log t}}\right] \text{ and } ||W||\le \frac{\epsilon t^{M/2}}{C(\delta)}\left(\frac{M\log t}{t}\right)^{1/4}\right)\\
&\ge & 1-P\left(S_t <\frac{p}{4}\left[(\alpha'-\delta)\sqrt{\frac{t}{M \log t}}\right]\right)\\
&\quad & \hspace{2cm}-P\left(||W||> \frac{\epsilon t^{M/2}}{C(\delta)}\left(\frac{M \log t}{t}\right)^{1/4}\right),
\end{eqnarray*}
which goes to one as $t \rightarrow \infty$ by Lemma \ref{lem:looptwo}.\qed\\\\
The proof of the lower bound in Theorem \ref{thm:loopcount} follows from the above lemmas.\qed
\section{\textbf{Loops and singularities in particle paths}}\label{section:loopevolution}
We start off this section by describing an interesting phenomenon that one notices in simulations of the paths of the individual particles in the discrete Conga line. The leading particle ($k=1$) performs an erratic Gaussian random walk. But as $k$ increases, the successive particles are seen to \textit{cut corners} in the paths of the preceding particles making them smoother (see Figure 5). This can be heuristically explained by the fact that a particle following another one in front directs itself along the shortest path between itself and the preceding particle (see equation (\ref{equation:recursion})), and hence cuts corners. This phenomenon is captured by the process $\baru$ described in Subsection \ref{subsection:related}. So, we use the approximation of the discrete Conga line $X$ by the smooth process $\baru$ in this section.\\\\
Recall that a \textbf{singularity} of a curve $\gamma: \mathbb{R} \rightarrow \mathbb{C}$ is a point $t_0$ at which its speed vanishes, i.e. $|\gamma'(t_0)|=0$. A singularity $t_0$ of a curve $\gamma$ is called a \textbf{cusp singularity} if it is analytic in a neighborhood of $t_0$ and there exists a translation and rotation of co-ordinates taking $\gamma(t_0)$ to the origin, under which, $\gamma$ has the representation $\gamma^*=(\gamma^*_1,\gamma^*_2)$ with the following power series expansions:
\begin{eqnarray}\label{eqnarray:analexp}
\gamma^*_1(t)&=&\sum_{i=2}^{\infty}a_i(t-t_0)^i\nonumber\\
\gamma^*_2(t)&=&\sum_{i=3}^{\infty}b_i(t-t_0)^i
\end{eqnarray}
with $a_2\neq0, \ b_3\neq0$, for $t \in [t_0-\delta, t_0+\delta]$ for some $\delta>0$.\\ 
Intuitively, this means that the graph of $\gamma$ locally around $t_0$ looks like $y=x^{2/3}$ under a rigid motion of co-ordinates taking $\gamma(t_0)$ to the origin. We call these transformed co-ordinates the \textbf{natural co-ordinate frame} based at the cusp singularity $t_0$.\\\\
Making a change of variables $p=t-\frac{x}{\alpha}$ and $\tau=\rho^2 x$, we can rewrite $\baru$ as
\begin{equation}\label{equation:pathfn}
f(p,\tau)=\bar{u}(x,t)=\mathbb{E}_{Z_{\tau}}W(p-Z_{\tau})
\end{equation}
We restrict our attention to $p>0, \tau>0$. Fixing $\tau$ and varying $p$ in the above expression for $f$ yields the path of the particle at distance $x=\tau/\rho^2$ from the tip. This change of variables enables us to write $\baru$ as the solution to the \textit{heat equation} with initial function being the Brownian motion $W$, the space variable represented by $p$ and time by $\tau$. Later, we will see that this makes it easier to write down analytic expansions of $\baru$ around the singularity.

Another interesting observation, which was described briefly in the Introduction, is the \textit{evolution of loops} as we look at the paths of successive particles. If a particle in the (two dimensional) Conga line goes through a loop, the particle following it, which cuts corners and tries to ``catch it", will go through a smaller loop. This is suggested by the simulations, where small loops are seen to `die' (i.e. shrink to a point), and just before death, they look somewhat `elongated', and the death site looks like a \textit{cusp singularity}. Other loops are seen to break after some time, that is, their end points come apart. Figure 5, representing successive particle paths, depicts some loops in various stages of evolution. In this section, we investigate evolving loops in the paths of successive particles, especially the relationship between dying loops and formation of singularities.

Before we can start off, we give some definitions that will be useful in describing the evolution of loops.

We define a metric space $(\mathcal{M},d)$, with a metric similar to the Skorohod metric on RCLL paths, on which we want to study loop evolutions:

$\displaystyle{\mathcal{M}=\{f: [a,b] \rightarrow \mathbb{C}; -\infty < a < b < \infty\}}$\\\\
If $f:[a_1,b_1] \rightarrow \mathbb{C},g:[a_2,b_2] \rightarrow \mathbb{C} \in \mathcal{M}$, define
$$\displaystyle{d(f,g)=\inf\left\lbrace||\lambda_1-\lambda_2|| + ||f\circ \lambda_1-g \circ \lambda_2||\ \  \vert \ \ \lambda_i:[0,1] \rightarrow [a_i,b_i] \text{ is a homeomorphism}\right\rbrace}$$
where $||\cdot||$ denotes the sup-norm metric. It can be easily checked that $(\mathcal{M},d)$ is a metric space.\\\\ Define the \textbf{evolution of a loop $L$} as a continuous function $\displaystyle{L(\cdot):[0,T) \rightarrow \mathcal{M}}$ such that $L(0)=L$ and $L(t)$ is a loop for every $0\le t<T$. If $f:\mathbb{R}^+ \times \mathbb{R}^+ \rightarrow \mathbb{C}$ is a continuous space-time process, and $\displaystyle{L(t)=f(\cdot,t)|_{[a_t,b_t]}}$ is a loop evolution on $0 \le t < T$, we say that $\{L(t): 0\le t < T\}$ is a \textbf{loop evolution of $f$ starting from $L=f(\cdot,0)|_{[a_0,b_0]}$}. Say that the loop $L=f(\cdot,0)|_{[a_0,b_0]}$ \textbf{vanishes after time $T^*$} if $$\displaystyle{T^*=\sup\{T>0: \text{ there exists a loop evolution }L(\cdot):[0,T) \rightarrow \mathcal{M} \text{ of }f \text{ starting from }L\}}.$$ Such a vanishing loop is said to \textbf{die at space-time point $(p,T^*)$} if $a_t \rightarrow p$ and $b_t-a_t \rightarrow 0$ as $t \rightarrow T^*$. Analogous definitions hold for $[0,T)$ replaced by a general interval $[T_1,T_2)$.

We note here that \textit{loops can vanish without dying}. This can happen, for example, when the `end points come apart'. One can check that an instance of this happening for $f$ defined in (\ref{equation:pathfn}) is when the end points of the loop have their velocity vectors parallel, but normal acceleration vectors pointing in opposite directions.\\\\
\textbf{Note: }Although with probability one $f$ has \textit{no singularities for a fixed time $\tau$}, it can be verified by an application of the expectation meta-theorem of \cite{adlertaylor} that the expected number of singularities of $f(\cdot,\tau)$ for $(p,\tau)$ lying in a compact set $K=[a,b] \times [c,d]$ is positive, and thus \textit{singularities do occur with positive probability if we allow both space and time to vary}.

It is easy to see that if a loop dies at a site $(p_0,\tau_0)$, then $p_0$ is a singularity for the curve $f(\cdot,\tau_0)$. We prove in Lemma \ref{lem:singularity} that with probability one, \textit{any singularity is a cusp singularity}. In Theorem \ref{thm:loop} we prove that for any (cusp) singularity $p_0$ in $f(\cdot,\tau_0)$, \textit{there exists a unique loop at each time in some small interval} $[\tau_0-\delta,\tau_0)$. Furthermore, the loop at time $\tau_0-\delta$ \textit{dies} at $(p_0,\tau_0)$ to give birth to the singularity $p_0$ of $f(\cdot,\tau_0)$. Also, the theorem shows that the dying loops, \textit{under some rescaling, converge to a deterministic limiting loop}.
\begin{lem}\label{lem:singularity}
With probability one, any singularity $p_0$ of the curve $f(\cdot,\tau_0)$ is a cusp singularity.
\end{lem}
\textbf{Proof: }We first prove that for every $\tau_0 \in \mathbb{R}^+$, $f(\cdot,\tau_0)$ is analytic. For this, note that\\ $$\displaystyle{\partial_p^nf(p,\tau)=\int_0^{\infty}(-1)^n\tau^{-(n+1)/2}\operatorname{He}_n\left(\frac{y-p}{\sqrt{\tau}}\right)\phi\left(\frac{y-p}{\sqrt{\tau}}\right)W(y)dy}.$$ By using the fact that $\displaystyle{\lim_{y\rightarrow \infty}\frac{W(y)}{y}=0}$ almost surely, we get that with probability one,
\begin{eqnarray}\label{eqnarray:dubla}
\left|\partial_p^nf(p,\tau)\right| \le \frac{C}{\tau^{n/2}}\sqrt{n!}(p^2+\tau)
\end{eqnarray}
for some random constant $C$. This bound implies analyticity of $f(\cdot,\tau_0)$.

Write $f=(f^1,f^2)$. We need to show that if $p_0$ is a singularity of $f(\cdot,\tau_0)$, then there exists a rigid motion of co-ordinates taking $f(p_0,\tau_0)$ to the origin under which (\ref{eqnarray:analexp}) holds (with $\gamma(\cdot)$ replaced by $f(\cdot,\tau_0)$). It suffices to prove the lemma for $(p_0,\tau_0)$ lying in a fixed rectangle $K=[a,b] \times [c,d]$. 
Our first step towards this is to show the following:
\begin{eqnarray}\label{eqnarray:cusp}
P( \exists (p_0,\tau_0) \in K \text{ such that } \partial_p f(p_0,\tau_0)=0 \text{ and the vectors }\left(\partial^2_p f^1(p_0,\tau_0),\partial^3_p f^1(p_0,\tau_0)\right) \text{ and }\nonumber\\
\left(\partial^2_p f^2(p_0,\tau_0),\partial^3_p f^2(p_0,\tau_0)\right) \text{ are linearly dependent})=0
\end{eqnarray}
To show this, define $A_n$ to be the event which holds when all the following are satisfied:
\begin{itemize}
\item[(i)]There exists $(p_0,\tau_0) \in K$ for which $\partial_p f(p_0,\tau_0)=0$ and $$\left(\partial^2_p f^2(p_0,\tau_0),\partial^3_p f^2(p_0,\tau_0)\right)=\lambda \left(\partial^2_p f^1(p_0,\tau_0),\partial^3_p f^1(p_0,\tau_0)\right)$$ for some $\lambda \in [-n,n]$.
\item[(ii)]The Lipschitz constants of the functions $\{\partial^i_p f^j(p,\tau): (p,\tau) \in K; i=1, 2, 3; j=1,2\}$ are less than or equal to $n$.
\end{itemize}
We will show that $P(A_n)=0$ which will yield (\ref{eqnarray:cusp}).\\\\
Partition the rectangle into a grid of sub-rectangles of side length $\le \epsilon$, where $\epsilon$ is small. Call the set of grid points $\hat{K}$.\\ Now, suppose $A_n$ holds. Let $(p_0,\tau_0)$ lie in a sub-rectangle $R$ and let $(p_i,\tau_j) \in \hat{K}$ be a grid point adjacent to $R$. Note that as the Lipschitz constants of the above functions and $\lambda$ are bounded by $n$, the following event holds:
\begin{eqnarray*}
A_n^{ij}&=& \left\lbrace |\partial_p f(p_i,\tau_j)|\le \sqrt{2}n\epsilon \text{ and }\right.\\\\
&\quad & \left.|\left(\partial^2_p f^2(p_i,\tau_j),\partial^3_p f^2(p_i,\tau_j)\right)-\lambda \left(\partial^2_p f^1(p_i,\tau_j),\partial^3_p f^1(p_i,\tau_j)\right)| \le 4n^2 \epsilon \right.\\\\
 &\quad &\left. \text{ for some }\lambda \in [-n,n]\right\rbrace.
\end{eqnarray*}
Thus we have $$\displaystyle{A_n \subseteq \bigcup_{i,j}A_n^{ij}}.$$ We show that there is a constant $C$ depending on $n$ such that $P\left(A_n^{ij}\right) \le C\epsilon$.\\\\
To save us notation, call $$X=(X_1,X_2,X_3)=(\partial_p f^1(p_i,\tau_j),\partial_p^2 f^1(p_i,\tau_j),\partial_p^3 f^1(p_i,\tau_j))$$ and similarly $Y$ for $f^2$. $X$ and $Y$ are independent and each follows a centred trivariate normal distribution. Let us call the density function of $X$ $p_{ij}$ and the distribution of $Y$ as $Q_{ij}$. Then, as $X$ and $Y$ have uniformly bounded densities,
\begin{eqnarray*}
P\left(A_n^{ij}\right)&\le &\int_{x \in [-\sqrt{2}n\epsilon,\sqrt{2}n\epsilon] \times \mathbb{R}^2}p_{ij}(x)\\\\
&\quad & \hspace{1cm} \times Q_{ij}\left(|Y_1| \le \sqrt{2}n\epsilon ; (Y_2,Y_3) \in \text{ the $4n^2\epsilon$ neighbourhood of the linear span of $(x_2,x_3)$}\right)dx\\
& \le & \int_{x \in [-\sqrt{2}n\epsilon,\sqrt{2}n\epsilon] \times \mathbb{R}^2}p_{ij}(x)C_{ij}'\epsilon^2 dx \le C_{ij}\epsilon^3 ,
\end{eqnarray*}
where $C_{ij}' , C_{ij}$ depend on $n$. Note that the determinants of the covariance matrices of $X$ and $Y$ are continuous and do not vanish at any point on the compact set $K$. Thus we can bound $C_{ij}$ by $C$ (which depends on $n$) uniformly over $i,j,\epsilon$. Using these facts, we get
\begin{eqnarray*}
P(A_n) \le \sum_{i,j}P\left(A_n^{ij}\right) \le C\epsilon.
\end{eqnarray*}
As $\epsilon$ is arbitrary, we get $\displaystyle{P(A_n)=0}$.\\\\
Now if $p_0$ is a singularity occurring at time $\tau_0$, i.e. $\exists$ $(p_0,\tau_0) \in K$ for which $\partial_p f(p_0,\tau_0)=0$, we can apply a rigid motion of co-ordinates such that $f(p_0,\tau_0)$ is the new origin and the rotation angle $\theta$ is chosen to satisfy the equation
\begin{eqnarray*}
A_{\theta}
\begin{bmatrix}
\partial_p^2 f^1(p_0,\tau_0) & \partial_p^3 f^1(p_0,\tau_0)\\
\partial_p^2 f^2(p_0,\tau_0) & \partial_p^3 f^2(p_0,\tau_0)
\end{bmatrix}
=
\begin{bmatrix}
a_2 & a_3\\
0 & b_3
\end{bmatrix}
\end{eqnarray*}
where $A_{\theta}$ is the rotation matrix corresponding to $\theta$. By (\ref{eqnarray:cusp}), we see that $a_2$ and $b_3$ are non-zero. Then, the result follows by taking this new co-ordinate frame as the natural co-ordinates.
\qed
\begin{thm}\label{thm:loop}
If $p_0$ is a (cusp) singularity of $f(\cdot, \tau_0)$, then there exists $\delta >0$ such that $f(\cdot,\tau)$ has a unique loop $L(\tau)=f(\cdot,\tau)\mid_{[a_{\tau},b_{\tau}]}$ in an interval $I_{\tau}$ containing $p_0$ for all $\tau \in [\tau_0-\delta, \tau_0)$. Moreover, $I_{\tau}$ can be chosen so that it shrinks to $p_0$ as $\tau \rightarrow \tau_0$, and $L(\cdot): [\tau_0-\delta,\tau_0) \rightarrow \mathcal{M}$ is a loop evolution of $f$ starting from $f(\cdot, \tau_0-\delta)\mid_{[a_{\tau_0-\delta},b_{\tau_0-\delta}]}$.

Furthermore, if $f=(f_1,f_2)$ in natural co-ordinates based at the cusp singularity $p_0$ of $f(\cdot, \tau_0)$, then there is $M>0$ such that the rescaling of $f$ given by
\begin{eqnarray}\label{eqnarray:frescale}
\hat{f}_s(P)=\left(s^{-1}f^1(\sqrt{s}P, \tau_0-s), s^{-3/2}f^2(\sqrt{s}P, \tau_0-s)\right)
\end{eqnarray}
for $s \in (0,\delta]$ has a unique loop $\hat{L}_s$ in $[-M,M]$ which converges to a deterministic loop $\hat{L}_0$ in $(\mathcal{M}, d)$ as $s \rightarrow 0$.
\end{thm}
\textbf{Proof: }First we show that $f$ is jointly analytic in $(p,\tau)$ for $(p,\tau) \in \mathbb{R} \times \mathbb{R}^+$.
From (\ref{eqnarray:dubla}), we know that for any $\tau_0 \in \mathbb{R}^+$, $f(\cdot,\tau_0)$ has an analytic representation in the space variable $p$ given by
\begin{equation}\label{equation:jointanalone}
f(p,\tau_0)=\sum_{i=0}^{\infty}\mathbf{a}_i(p_0,\tau_0)(p-p_0)^i.
\end{equation}
which holds for all $p \in \mathbb{R}$, where $\mathbf{a}_i(p_0,\tau_0) =\partial_p^i f(p_0,\tau_0) \in \mathbb{R}^2$ for each $i$. Note that if we can prove
\begin{equation}\label{equation:dubladee}
\sum_{i=0}^{\infty}\mathbb{E}_{Z_{\delta}}|\mathbf{a}_i(p_0,\tau_0)|\left||p-p_0|+|Z_{\delta}|\right|^i < \infty
\end{equation}
for $p \in \mathbb{R}$ and $0\le \delta< \epsilon$ for some $\epsilon>0$, then we can write
\begin{equation}\label{equation:jointanaltwo}
f(p,\tau)=\sum_{i=0}^{\infty}\mathbb{E}_{Z_{\tau-\tau_0}} \mathbf{a}_i(p_0,\tau_0)\left((p-p_0)-Z_{\tau-\tau_0}\right)^i
\end{equation}
for $p \in \mathbb{R}$ and $\tau \in [\tau_0,\tau_0+\epsilon)$ and it follows that $f$ is jointly analytic on $\mathbb{R}\times [\tau_0,\tau_0+\epsilon)$. Joint analyticity on $\mathbb{R} \times \mathbb{R}^+$ is immediate as a result.

From (\ref{eqnarray:dubla}), we see that (\ref{equation:dubladee}) holds when
\begin{eqnarray*}
\sum_{i=0}^{\infty}\frac{C}{\tau_0^{i/2}\sqrt{i!}}(p_0^2+\tau_0)2^i(|p-p_0|^i+2^{i/2}\sqrt{i!}\delta^{i/2})<\infty
\end{eqnarray*}
which is satisfied for $\displaystyle{\delta < \frac{\tau_0}{8}}$ and all $p \in \mathbb{R}$.

If we set $\zeta_n(p,t)=\mathbb{E}_{Z_t}((p-p_0)-Z_t)^n$ for $p \in \mathbb{R}$ and $t\ge0$, we note that this can be readily computed from the moments of $Z_t$. Further, observe that the polynomial thus obtained for $\zeta_n(p,t)$ is well defined for all $(p,t) \in \mathbb{R}^2$. The joint analyticity of $f$ along with (\ref{equation:jointanaltwo}) tells us that we can write $f(p,\tau)=\sum_{i=0}^{\infty}\mathbf{a}_i(p_0,\tau_0)\zeta_i(p,\tau-\tau_0)$ for $(p,\tau)$ in the region of analyticity around $(p_0,\tau_0)$.

Choose and fix $M \in (\sqrt{3}, 2\sqrt{3}-1)$. If $f(\cdot,\tau_0)$ has a cusp singularity at $p_0$, we can write $f=(f_1,f_2)$ in the natural co-ordinate frame based at $p_0$. We will show that there is $\delta>0$ small enough for which the rescaled process $\hat{f}_s$ given by (\ref{eqnarray:frescale}) in this co-ordinate frame
has a \textit{unique} loop $\hat{L}_s$ in the interval $[-M,M]$ for all $0<s \le \delta$. We will in fact show something stronger: for $0<s \le \delta$, there is a \textit{unique pair} $(\hat{a}_s,\hat{b}_s)$ such that $\hat{a}_s, \hat{b}_s \in [-M,M]$, $\hat{a}_s < \hat{b}_s$ and $\hat{f}_s(\hat{a}_s)=\hat{f}_s(\hat{b}_s)$.

This will imply that for all $0<s\le\delta$, $f(\cdot, \tau_0-s)$ has a unique loop $L(\tau_0-s)=f(\cdot, \tau_0-s)\mid_{[a_{\tau_0-s},b_{\tau_0-s}]}$ in $[-M\sqrt{s}, M\sqrt{s}]$, where $a_{\tau_0-s}=\sqrt{s}\hat{a}_s$ and $b_{\tau_0-s}=\sqrt{s}\hat{b}_s$.\\\\
By scaling properties of Gaussian random variables, $\hat{f}_s$ has the following representation:
\begin{eqnarray}\label{eqnarray:frescanal}
\hat{f}_s^1(P)&=&\sum_{i=2}^{\infty}a_is^{(i-2)/2}\zeta_i(P,-1)\nonumber\\
\hat{f}_s^2(P)&=&\sum_{i=3}^{\infty}b_is^{(i-3)/2}\zeta_i(P,-1).
\end{eqnarray}
Let $$\hat{g}(P)=(a_2\zeta_2(P,-1),b_3\zeta_3(P,-1))=(a_2(P^2-1),b_3(P^3-3P)).$$ Then by (\ref{eqnarray:frescanal}) there is a (random) constant $C$ depending only on $M$ and $\delta$ such that
\begin{equation}\label{equation:hatgone}
|\hat{f}_s(P)-\hat{g}(P)| \le C\sqrt{s}
\end{equation}
and
\begin{equation}\label{equation:hatgtwo}
|\hat{f}_s'(P)-\hat{g}'(P)| \le C\sqrt{s}.
\end{equation}
for all $P \in [-M,M]$ and all $s \in (0,\delta]$.

The bound (\ref{equation:hatgone}) shows that on $[-M,M]$, $\hat{f}_s$ is forced to run in a `tube' of width $C\sqrt{s}$ around $\hat{g}$. Note that $\hat{g}$ has positive speed in $[-M,M]$, and has a unique loop $\hat{g}\mid_{[-\sqrt{3},\sqrt{3}]}$ in this interval. This `forces' $\hat{f}_s$ to intersect itself. More precisely, we can find a (random) constant $A>0$ depending only on $M$ and $\delta$ such that there exist $\hat{a}_s \in [-\sqrt{3} - A\sqrt{s}, -\sqrt{3}+A\sqrt{s}]$ and $\hat{b}_s \in [\sqrt{3}-A\sqrt{s}, \sqrt{3}+A\sqrt{s}]$ satisfying $\hat{f}_s(\hat{a}_s)=\hat{f}_s(\hat{b}_s)$ for $0<s \le \delta$. We will now show that $(\hat{a}_s, \hat{b}_s)$ is the required pair.

Notice that the only point in $[-M,M]$ at which the derivative of $\hat{g}^1$ vanishes is $P=0$, and the only such points for $\hat{g}^2$ are $P=1$ and $P=-1$. Thus, if we choose any $0< \theta <1$, it follows from (\ref{equation:hatgtwo}) that for small enough $\delta$, there is $I>0$ (depending only on $M$ and $\delta$) such that $$\displaystyle{\inf_{P_1,P_2 \in [-M,M], |P_1-P_2| \le \theta}\max\left\lbrace\left|\frac{d}{dP}\hat{f}_s^1(P_1)\right|, \left|\frac{d}{dP}\hat{f}_s^2(P_2)\right|\right\rbrace \ge I}$$ for $0<s \le \delta$. This, in turn, gives us
\begin{equation}\label{equation:derlb}
\inf_{P,Q \in [-M,M], 0<|P-Q| \le \theta}\frac{|\hat{f}_s(P)-\hat{f}_s(Q)|}{|P-Q|} \ge I
\end{equation}
for $0<s \le \delta$. In particular, there does not exist $P,Q \in [-M,M]$ with $0<|P-Q| \le \theta$ such that $\hat{f}_s(P)=\hat{f}_s(Q)$.

Take and fix any $\epsilon \in (0,2\sqrt{3}-1-M)$ and $L \in (1+M+\epsilon,2\sqrt{3})$. As the only self-intersection of $\hat{g}$ in $[-M,M]$ is given by $\hat{g}(-\sqrt{3})=\hat{g}(\sqrt{3})$, $$\displaystyle{m=\inf_{P,Q \in [-M,M], |\theta \le P-Q| \le L}|\hat{g}(P)-\hat{g}(Q)|>0}.$$ thus, for sufficiently small $\delta$,
\begin{equation}\label{equation:dertwo}
\inf_{P,Q \in [-M,M], |\theta \le P-Q| \le L}|\hat{f}_s(P)-\hat{f}_s(Q)| \ge m/2>0
\end{equation}
 for $s \in (0,\delta]$. (\ref{equation:derlb}) and (\ref{equation:dertwo}) together imply that if $P,Q \in [-M,M]$, $P<Q$ such that $\hat{f}_s(P)=\hat{f}_s(Q)$, then $P \in [-M,-1-\epsilon]$ and $Q \in [1+\epsilon,M]$. As $(\hat{g}^1)'<0, (\hat{g}^2)'>0$ on $[-M,-1-\epsilon]$ and $(\hat{g}^1)'>0, (\hat{g}^2)'>0$ on $[1+\epsilon,M]$, therefore for sufficiently small $\delta$, the same relations hold with $\hat{g}$ replaced by $\hat{f}_s$ for $s \in (0,\delta]$. It is then routine to check that such a pair $(P,Q)$ is unique.
 
 Thus, for $\delta$ small enough, the \textit{unique} loop of $\hat{f}_s$ in $[-M,M]$ for $s \in (0,\delta]$ is given by $\hat{L}_s=\hat{f}_s\mid_{[\hat{a}_s,\hat{b}_s]}$. As $\max\{|\hat{a}_s+\sqrt{3}|, |\hat{b}_s-\sqrt{3}|\} \le A\sqrt{s}$, this fact, along with (\ref{equation:hatgone}) and (\ref{equation:hatgtwo}), imply convergence of $\hat{L}_s$ to $\hat{L}_0=\hat{g}\mid_{[-\sqrt{3},\sqrt{3}]}$ in $(\mathcal{M},d)$.
 
 The final thing left to prove is the continuity of the map $L(\cdot):[\tau_0-\delta,\tau_0) \rightarrow \mathcal{M}$. To show this, take any $s_0 \in (0,\delta]$. By joint analyticity of $f$, for small enough $\eta>0$, $\sup_{p \in [-M\sqrt{s_0},M\sqrt{s_0}]}|f(p,\tau_0-s)- f(p,\tau_0-s_0)| \le C|s-s_0|$ whenever $|s-s_0| \le \eta$, for some (random) constant $C$ depending only on $M, s_0$ and $\eta$. This, along with the fact that $f(\cdot, \tau_0-s_0)$ has positive speed in $[-M\sqrt{s_0},M\sqrt{s_0}]$ implies that $\max\{|a_s-a_{s_0}|, |b_s-b_{s_0}|\} < A|s-s_0|$ whenever $|s-s_0| \le \eta$, where $A$ is again a (random) constant depending on $M, s_0$ and $\eta$. This, along with uniform continuity of $f$ in $[-M\sqrt{s_0},M\sqrt{s_0}] \times [s_0-\eta,s_0+\eta]$ implies that $L(\tau_0-s) \rightarrow L(\tau_0-s_0)$ in $(\mathcal{M},d)$ as $s \rightarrow s_0$ proving continuity.
\qed\\\\
From Lemma \ref{lem:singularity} and Theorem \ref{thm:loop}, it is clear that \textit{there is a bijection between dying loops and singularities in the particle path evolution}.\\\\
\textbf{Shape of a dying loop:}
Observe that in the previous theorem, the two coordinates had to be scaled differently to get the limiting loop. The difference in the scaling exponent explains why \textit{the loops look elongated before death} (see Figure 5).
\section{\textbf{Freezing in the tail}}\label{section:freezing}
This section addresses Observation 4 of Burdzy and Pal \cite{burdzypal}. The \textit{tail of the Conga line} refers to the particles at distance $x > \alpha t$ from the leading particle.

Note that, from the scaling relation (\ref{equation:exeq}) and equation (\ref{equation:Browcloseprob}), we have for any fixed $\beta \in (0,1)$,
\begin{equation}\label{equation:headnofreeze}
\P\left(\sup_{[x \in [t^{1-\beta},\alpha t]}\left| u(x,t)-W\left(t-\frac{x}{\alpha}\right)\right| > \kappa t^{\frac{1}{4}}\sqrt{(1+\beta)\log t}\right) \le t^{-2}.
\end{equation}
Thus, almost surely, there is $N$ such that for all $n \ge N$,
\begin{equation}\label{equation:headas}
\sup_{[x \in [n^{1-\beta},\alpha n]}\left| u(x,n)-W\left(n-\frac{x}{\alpha}\right)\right| \le \kappa n^{\frac{1}{4}}\sqrt{(1+\beta)\log n}
\end{equation}
for all $n \ge N$. In particular, for any $\delta \in (0,\alpha)$, if we look at the family of curves $$\mathcal{U}_n(x)=(2n \log \log n)^{-\frac{1}{2}}u(nx,n)$$ for $\delta \le x \le \alpha$, then (\ref{equation:headas}), along with Strassen's Theorem for Brownian motion (see \cite{bill} p. 221, Theorem 22.2), tells us that the set of continuous curves that form the limit points of $\{\mathcal{U}_n(x): \delta \le x \le \alpha\}_{n \ge 1}$ with respect to sup-norm distance on $[\delta, \alpha]$ is precisely given by
\begin{eqnarray*}
\mathcal{K}&=&\left\lbrace \Gamma\left(1-\frac{x}{\alpha}\right), x \in [\delta, \alpha]; \Gamma:[0,1] \rightarrow \mathcal{C} \text{ is absolutely continuous with } \Gamma(0)=0\right.\\
&\quad& \left.\text{ and } \int_0^1|\Gamma'(s)|^2ds \le 1\right\rbrace.
\end{eqnarray*}
Thus, the Conga line $u$ under the above scaling shows appreciable movement in time. In other words, \textit{we do not observe freezing in the scale of the driving Brownian motion $W$}.

On the other hand, if we zoom in on the tail of the Conga line, the particles seem to \textit{freeze in time}. Furthermore, the \textit{direction} in which the particles come out of the origin shows very little change with time after a while (see Figure 2 and Figure 5). The small variance of particles in the tail region compared to that of the driving Brownian motion $W$ for large $t$ necessitates a \textit{rescaling} of the tail to study its properties. Thus, the tail behaves in a very different manner compared to particles near the tip.

To study the phenomenon of `freezing in the tail', we use the continuous version $\baru$ described in Subsection \ref{subsection:related}. For any fixed $\eta \in (0, (1-\alpha)/\alpha)$, we look at the following distances from the tip $$x_t=\alpha(1 + \eta)t$$ and study $$\overline{v}_t(\eta)=\baru(x_t, t).$$ The choice of distances $x_t$ from the tip ensures that these particles remain in the tail region for all $t$. We rescale $\overline{v}_t$ as
\begin{equation}\label{equation:freeze}
\displaystyle{v_t(\eta)=\sqrt{2\pi \rho^2\alpha(1+\eta)t} \ \exp\left\lbrace\frac{\eta^2t}{2\rho^2\alpha(1+\eta)}\right\rbrace\overline{v}_t(\eta)}.
\end{equation}
Also define
\begin{equation}
v(\eta)=\int_0^{\infty}\exp\left\lbrace-\frac{s\eta}{\rho^2\alpha(1+\eta)}\right\rbrace W(s) ds,
\end{equation}
where $W$ is the driving Brownian motion in expression (\ref{equation:baru}).
\begin{thm}\label{thm:freezing}
For any fixed $\eta \in (0, (1-\alpha)/\alpha)$,
\begin{eqnarray*}
\displaystyle{v_t(\eta) \rightarrow v(\eta)}
\end{eqnarray*}
almost surely and in $L^2$ as $t \rightarrow \infty$.
\end{thm} 
\textbf{Proof: }In the following, $C_1, C_2,\dots$ represent finite, positive constants.

From (\ref{equation:baru}) we can write
\begin{eqnarray*}
v_t(\eta)&=&\int_0^{\infty}W(s)\exp\left\lbrace-\frac{s^2}{2\rho^2 t\alpha(1+\eta)}-\frac{s\eta}{\rho^2\alpha(1+\eta)}\right\rbrace ds.\\
\end{eqnarray*}
Almost sure convergence follows from the fact that $$\int_0^{\infty}|W(s)|\exp\left\lbrace-\frac{s\eta}{\rho^2\alpha(1+\eta)}\right\rbrace ds < \infty$$ with probability one, and the Dominated Convergence Theorem.\\\\
To prove $L^2$ convergence, note that
\begin{eqnarray}\label{eqnarray:freezev}
v(\eta)-v_t(\eta)=\int_0^{\infty}g_t(s)W(s)ds=\int_0^{\infty}\left(\int_a^{\infty}g_t(s)ds\right)dW(a),
\end{eqnarray}
where $$g_t(s)=\exp\left\lbrace-\frac{s\eta}{\rho^2\alpha(1+\eta)}\right\rbrace\left(1-\exp\left\lbrace-\frac{s^2}{2\rho^2 t\alpha(1+\eta)}\right\rbrace\right).$$ It is easy to see that
\begin{equation}\label{equation:geebound}
g_t(s) \le \frac{C_1}{\rho^2 \alpha (1+ \eta)t}s^2\exp\left\lbrace-\frac{s\eta}{\rho^2 \alpha(1+\eta)}\right\rbrace \mathbb{I}(s \le \sqrt{t}) + C_2 \exp\left\lbrace-\frac{s\eta}{\rho^2 \alpha(1+\eta)}\right\rbrace \mathbb{I}(s > \sqrt{t}).
\end{equation}
From (\ref{equation:geebound}), we get $$\int_a^{\infty}g_t(s)ds \le \frac{C_3}{t}\mathbb{I}(a < \sqrt{t}) +C_4 \exp\left\lbrace-\frac{\eta(\sqrt{t} \vee a)}{\rho^2\alpha(1+\eta)}\right\rbrace.$$ Thus
\begin{equation}\label{equation:freezeerrone}
\mathbb{E}|v(\eta)-v_t(\eta)|^2=\int_0^{\infty}\left(\int_a^{\infty}g_t(s)ds\right)^2 da \le \frac{C_5}{t^{3/2}}
\end{equation}
giving $L^2$ convergence.\qed\\\\
\textbf{Remark: }Note that the above theorem only gives freezing (under rescaling) in a part of the Conga line that goes to zero exponentially fast with time. We have also seen that freezing does not occur in the scale of the driving Brownian motion. The interesting (and harder) question that remains to be investigated is about possible \textit{freezing in the intermediate regime} (in a window of the form $\left[\alpha t, \alpha t + \sigma\sqrt{\frac{1}{2}t \log t}\right]$ where a sharp transition in variance occurs, see Theorem \ref{thm:variance} and Corollary \ref{cor:vartrans}). We hope to investigate this in a future article.\\\\
\textbf{Conclusion}\\
In this paper, we approximate the \textit{smooth part} of the Conga line, i.e. $X_k(n)$ for $k>>1$, by the smooth random Gaussian process $u$ (equivalently $\baru$), and then investigate the path properties of $u$ and $\baru$. But this sheds little light on the Conga line near the tip $(k=1)$ as well as the motion of particles in time $n$ for smaller values of $k$. Note from (\ref{eqnarray:movtwo}) that for fixed $k$, the increment process $\{X_{k+1}(n)-X_k(n): n \ge 1\}$ behaves like a \textit{Weighted Moving Average} which finds frequent application in time series data analysis to smooth out short term fluctuations and highlight longer term trends and cycles (see \cite{hamilton}). From an analytical point of view, it would be interesting to analyse the Conga line at time $n$ in a suitable window around some sequence $a_n$ which grows to infinity sufficiently slowly so that the approximation by $u$ fails. We will try to address this in a future article.
\\\\ 
\textbf{Acknowledgements}\\
I thank my adviser Krzysztof Burdzy for guiding me through the project and giving interesting ideas to work on. I also thank Shirshendu Ganguly and Soumik Pal for helpful discussions. I thank Krzysztof Burdzy, Shirshendu Ganguly and Mary Solbrig for the figures, and Bharathwaj Palvannan, Yannick Van Huele and Tvrtko Tadic for helping me with computer graphics issues. I am grateful to an anonymous referee for providing valuable comments and suggestions that greatly improved the paper. This research was partly supported by NSF grant number DMS-1206276.
\section{\textbf{Simulations}}\vspace{1cm}
\begin{figure}[H]
\begin{center}
\includegraphics[width=120mm,height=120mm]{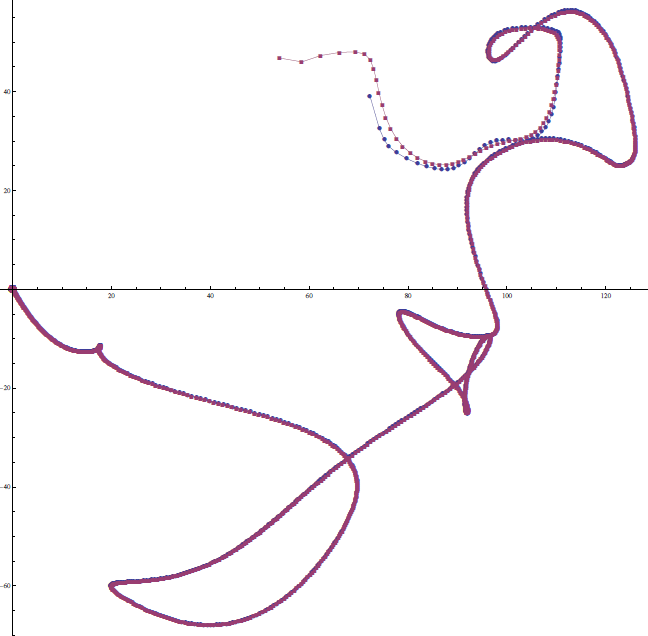}
\caption{\footnotesize \textit{The discrete two dimensional Conga line :} Plot of $\{X_k(10000): k=1,2,\dots,10000\}$ (red) and $\{X_k(9900): k=1,2,\dots,9900\}$ (blue) with $\alpha=0.5$. (courtesy: Krzysztof Burdzy)}
\end{center}
\end{figure}
\begin{figure}
\begin{center}
\includegraphics[width=100mm,height=100mm]{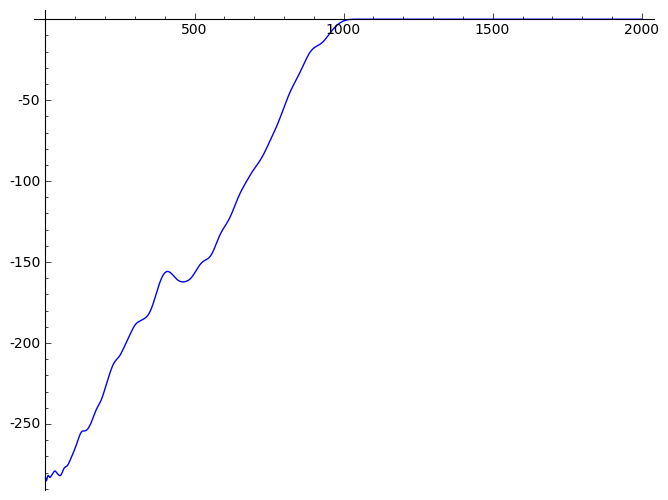}
\caption{\footnotesize \textit{The discrete one dimensional Conga line :} Plot of $\{X_k(2000): k=1,2,\dots,2000\}$ with $\alpha=0.5$. (courtesy: Shirshendu Ganguly)}
 \end{center}
\vspace{1cm}
\begin{center}
\includegraphics[width=90mm,height=90mm]{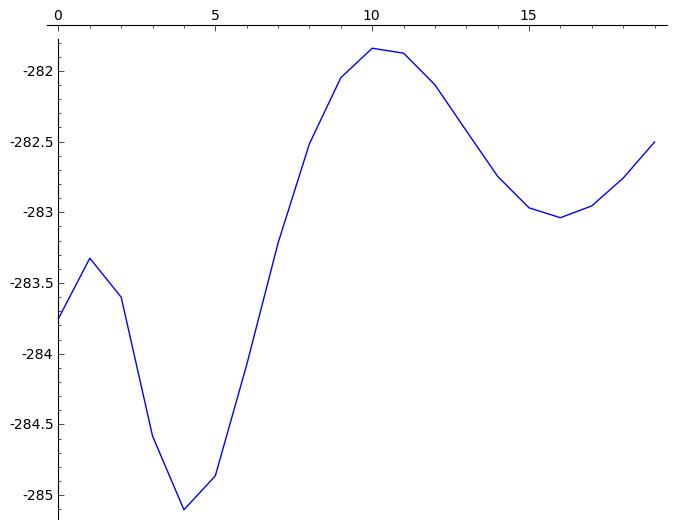}
\caption{\footnotesize \textit{Near the tip of the Conga line :} Plot of $\{X_k(2000): k=1,2,\dots,20\}$ with $\alpha=0.5$. (courtesy: Shirshendu Ganguly)}
 \end{center}
 \end{figure}
 \begin{figure}
\begin{center}
\includegraphics[width=120mm,height=95mm]{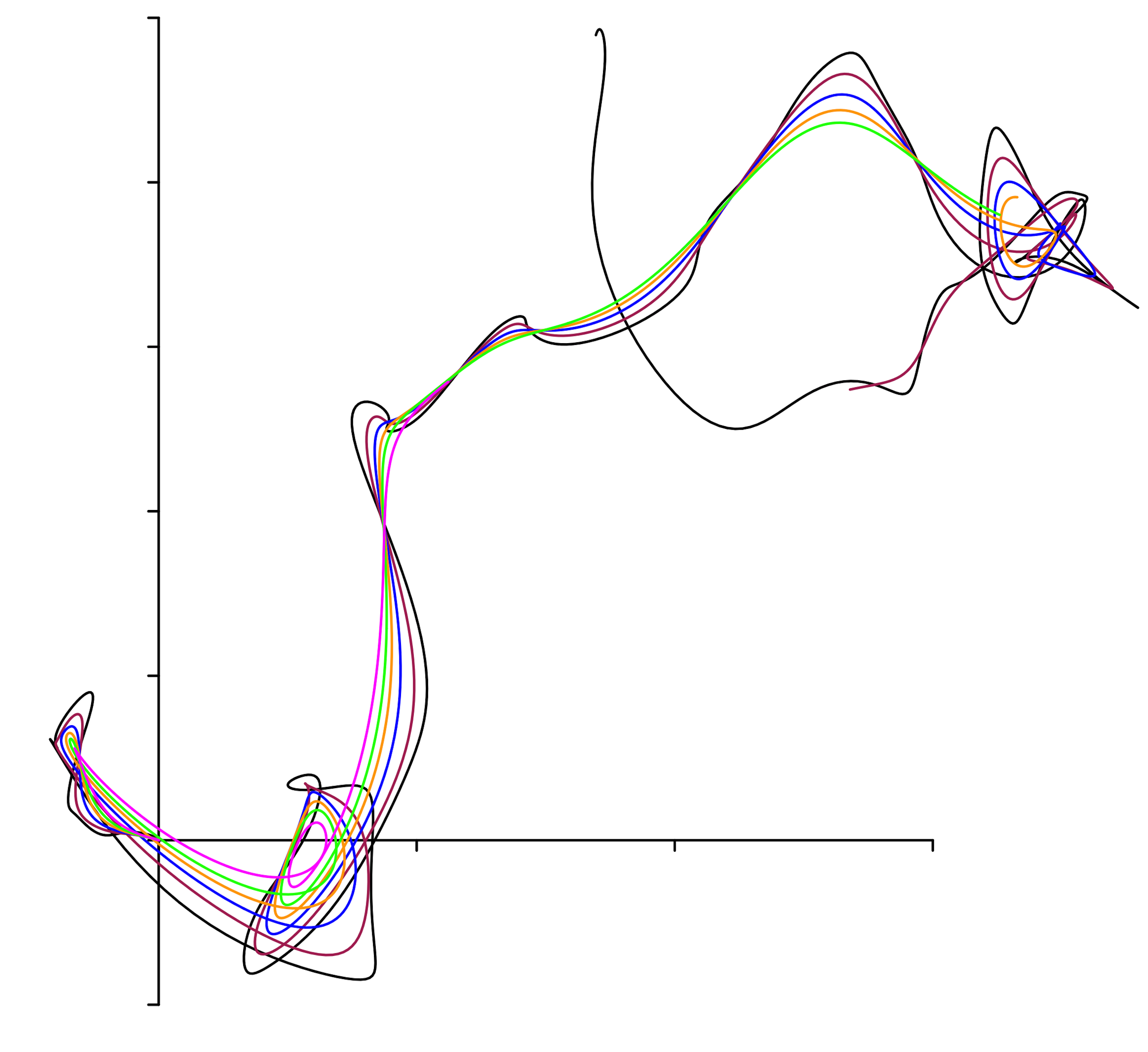}
\caption{\footnotesize \textit{Particle paths :} Plot of $\{X_k(n): n=1,2,\dots,800\}$ for $k=30$ (black), $60$ (maroon), $100$ (blue), $140$ (orange), $180$ (green) and $250$ (magenta).}
 \end{center}
 \end{figure}

\end{document}